\documentclass[12pt]{amsart}
\usepackage{graphicx}
\begin{document}
\newtheorem{theorem}{Theorem}[section]
\newtheorem{prop}[theorem]{Proposition}
\newtheorem{lemma}[theorem]{Lemma}
\newtheorem{claim}[theorem]{Claim}
\newtheorem{cor}[theorem]{Corollary}
\newtheorem{defin}[theorem]{Definition}
\newtheorem{example}[theorem]{Example}
\newtheorem{xca}[theorem]{Exercise}
\newcommand{\map}{\mbox{$\rightarrow$}}
\newcommand{\aaa}{\mbox{$\alpha$}}
\newcommand{\Aaa}{\mbox{$\mathcal A$}}
\newcommand{\bbb}{\mbox{$\beta$}}
\newcommand{\ccc}{\mbox{$\mathcal C$}}
\newcommand{\ddd}{\mbox{$\delta$}} 
\newcommand{\Ddd}{\mbox{$\Delta$}}
\newcommand{\Fff}{\mbox{$\mathcal F$}}  
\newcommand{\Ggg}{\mbox{$\Gamma$}}
\newcommand{\ggg}{\mbox{$\gamma$}}
\newcommand{\kkk}{\mbox{$\kappa$}}
\newcommand{\lll}{\mbox{$\lambda$}}
\newcommand{\Lll}{\mbox{$\Lambda$}}
\newcommand{\mlp}{\mbox{$\mu^{+}_{l}$}}
\newcommand{\ml}{\mbox{$\mu_{l}$}}
\newcommand{\mr}{\mbox{$\mu_{r}$}}
\newcommand{\mlpm}{\mbox{$\mu_{l}^{\pm}$}}
\newcommand{\mrpm}{\mbox{$\mu_{r}^{\pm}$}}
\newcommand{\mlm}{\mbox{$\mu_{l}^{-}$}}
\newcommand{\mrp}{\mbox{$\mu_{r}^{+}$}}
\newcommand{\mrm}{\mbox{$\mu_{r}^{-}$}}
\newcommand{\mm}{\mbox{$\mu^-$}}
\newcommand{\mpm}{\mbox{$\mu^{\pm}$}}
\newcommand{\mpp}{\mbox{$\mu^+$}}
\newcommand{\mt}{\mbox{$\mu^{t}$}}
\newcommand{\mb}{\mbox{$\mu_{b}$}}
\newcommand{\mz}{\mbox{$\mu^{\bot}$}}
\newcommand{\mpq}{\mbox{$\mu^{(p,q)}$}}
\newcommand{\omp}{\mbox{$0_{-}^{+}$}}
\newcommand{\oa}{\mbox{$\overline{a}$}}
\newcommand{\ob}{\mbox{$\overline{b}$}}
\newcommand{\opm}{\mbox{$0_{+}^{-}$}}
\newcommand{\opp}{\mbox{$0_{+}^{+}$}}
\newcommand{\Pt}{\mbox{$\tilde{P}$}}
\newcommand{\rrr}{\mbox{$\rho$}} 
\newcommand{\rz}{\mbox{$\rho^{\bot}$}}
\newcommand{\rp}{\mbox{$\rho^+$}}
\newcommand{\rmm}{\mbox{$\rho^-$}}
\newcommand{\rpq}{\mbox{$\rho^{(p,q)}$}}
\newcommand{\sss}{\mbox{$\sigma$}} 
\newcommand{\Sss}{\mbox{$\mathcal S$}} 
\newcommand{\sm}{\mbox{$\sigma^-$}}
\newcommand{\sz}{\mbox{$\sigma^{\bot}$}} 
\newcommand{\spm}{\mbox{$\sigma^{\pm}$}}
\newcommand{\spp}{\mbox{$\sigma^+$}}
\newcommand{\st}{\mbox{$\sigma^{t}$}}
\newcommand{\Ss}{\mbox{$\Sigma$}}
\newcommand{\Th}{\mbox{$\Theta$}} 
\newcommand{\ttt}{\mbox{$\tau$}} 
\newcommand{\bdd}{\mbox{$\partial$}}
\newcommand{\zzz}{\mbox{$\zeta$}}
\newcommand{\qb} {\mbox{$Q_B$}}
\newcommand{\inter}{\mbox{${\rm int}$}}

\numberwithin{equation}{section}

\title[] {There are no  unexpected tunnel number one knots of genus 
one}

\author{Martin Scharlemann}
\address{\hskip-\parindent
        Mathematics Department\\
        University of California\\
        Santa Barbara, CA 93106\\
        USA}
\email{mgscharl@math.ucsb.edu}

%\author{Martin Scharlemann}
\date{\today} 
\thanks{Research supported in part by an NSF grant, the Miller 
Institute, and RIMS Kyoto}

\begin{abstract} We show that the only knots that are tunnel number 
one and genus one are those that are already known: $2$-bridge knots 
obtained by plumbing together two unknotted annuli and the satellite 
examples classified by Eudave-Mu\~noz and by Morimoto-Sakuma.  This 
confirms a conjecture first made by Goda and Teragaito.
\end{abstract}
\maketitle

\section{Introduction and overview}

There are many useful ways of indexing the complexity of knot types: 
crossing number, bridge number, tunnel number, genus, etc.  
Often the relationship between these indices is unclear, and sometimes 
it is clear that there is no relationship.  Thus, for example, tunnel 
number one knots may be of arbitrarily high genus (e. g. torus knots) 
and genus number one knots may be of arbitarily high tunnel number (e. 
g. doubles of complicated knots.)  Given two indices of complexity it's 
natural to ask a sort of complementary question: how unusual is it for 
a knot (other than the unknot) to be of minimal complexity with 
respect to both indices?  For example, how rare is it that a knot 
have both genus one and tunnel number one?

\begin{figure} 
\centering
\includegraphics[width=0.15\textwidth]{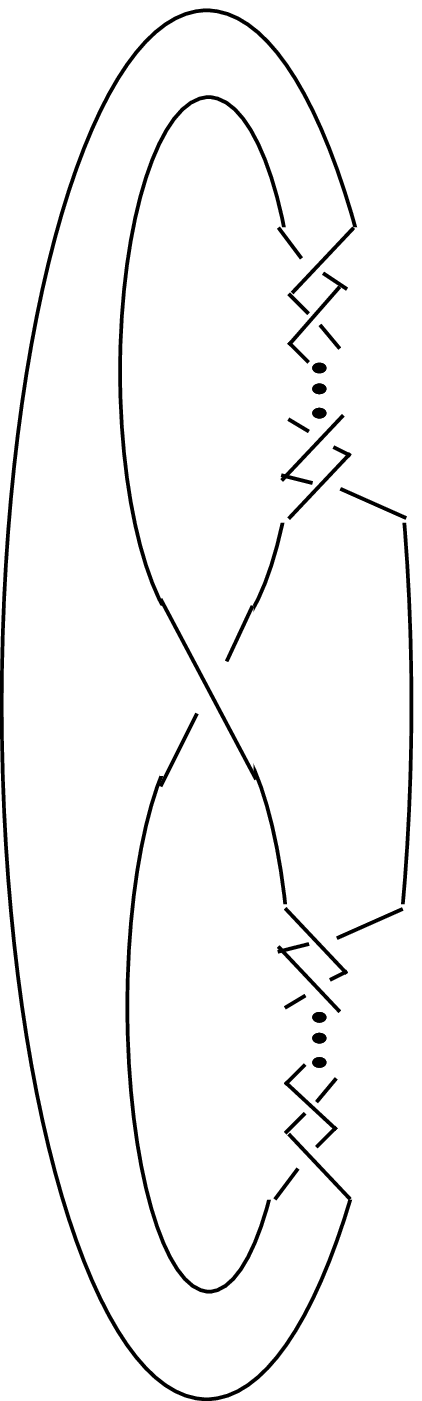}
\caption{} \label{fig:2bridge}
\end{figure}

It's easy to construct examples of knots of this type.  Plumb together 
two twisted unknotted annuli.  The boundary is typically knotted and 
the union of the annuli is visibly a genus one Seifert surface.  If we 
imagine hanging the union of annuli from a single peg we see that its 
boundary is naturally a $2$-bridge knot and therefore has tunnel 
number one.  See Figure \ref{fig:2bridge}.  It is known 
that these are the only $2$-bridge knots of genus one (cf 
\cite[Proposition 12.25]{BZ}).  

Are there other examples of genus one tunnel number one knots?  
Morimoto and Sakuma \cite{MS} and independently Eudave-Mu\~noz 
\cite{EM} classified satellite knots which have tunnel number 
one.  They have a concrete description and can be naturally indexed by 
a $4$-tuple of integers.  In \cite{GT}, Goda and Teragaito determined 
which of these satellite knots have genus one and made the conjecture 
that these knots complete the list of knots that have both genus one and 
tunnel number one.  The conjecture was confirmed by 
Matsuda \cite{Ma} for any knot that admits a $(1,1)$ decomposition; 
that is, for any knot which is $1$-bridge on an unknotted torus.

The central objective of this paper is to prove the Goda-Teragaito 
conjecture in complete generality.  The strategy will be to use thin 
position to show that any tunnel for a genus one knot can either be 
isotoped to lie on a genus one Seifert surface, or isotoped to form an 
unknotted loop.  In the latter case, it is shown that the knot admits 
a $(1,1)$ decomposition and Matsuda's argument applies.  In the 
former case, it follows from work of Eudave-Mu\~noz and Uchida 
\cite{EU} that $K$ is a $2$-bridge knot.

In retrospect, \cite{GST} and \cite{ST1} can be viewed as the first two 
steps of the program.  In \cite{GST} we show how thin position can be 
used to understand unknotting tunnels for unknotting number one 
knots.  In particular, we show that if $K$ and then $\lll$ are put in 
the thinnest possible position, then $\lll$ is level.  If $\lll$ is a 
loop, then Matsuda's theorem applies.  Otherwise, we show in 
\cite{ST1}, either $K$ is $2$-bridge or there is a well-defined 
invariant $\rrr \in \mathbb{Q}/2 \mathbb{Z}$ which, unless it is $1$, 
ensures that the tunnel can be moved onto a minimal genus Seifert 
surface.  So all that remains is to consider the case in which $\rrr = 
1$, which we do here.  

In the case that $\rrr = 1$ and the tunnel $\ggg$ is not an unknotted 
loop, it will be shown that there is another 
useful way of describing $K$ on the boundary of the genus $2$ 
handlebody $H = \eta(K \cup \ggg)$.  That is, there is a different spine for 
$H$, namely a $\Theta$ curve $\theta$, with these properties:

\begin{itemize}

\item The graph $\theta$ can be put in general position in $S^{3}$ in 
such a way that $K \subset \bdd\eta(\theta)$ remains in thin position.

\item $K$ intersects each meridian of each edge of $\theta$ always 
with the same orientation.

\item $\theta \subset S^{3}$ is thinner than the graph $K \cup 
\ggg$. 

\item A minimal genus Seifert surface $F$ for $K \subset H$ 
intersects $H$ only in $K = \bdd F$.

\end{itemize}

A combinatorial argument will show that if $K$ intersects each 
meridian more than once then $genus(F) \geq 2$.  Assuming that $genus(F) 
= 1$ and $\ggg$ cannot be isotoped to $F$, the program then will be to 
find the thinnest spine satisfying the conditions above (plus a more 
technical condition called the ``wave condition'').  For such a graph 
we know that $K$ intersects one of the meridians in only one point.  
We will argue, via thin position, that the cycle obtained by deleting 
this meridian is unknotted.  It will follow that $K$ has a $(1, 1)$ 
decomposition, so Matsuda's result applies.  

\section{Intersecting $(p, q)$ quasi-cables with spheres}

Consider a graph $\theta$ as just described.  Notice that the 
condition on meridians of edges of $\theta$ can be interpreted as 
follows: $K$ and each edge of $\theta$ can be oriented so that $K$ 
always runs along a given edge in the direction consistent with that 
edge's orientation.  In particular, $K$ intersects these meridians 
algebraically as well as geometrically in some $p, q, p+q$ points.  
With this in mind, we establish the more general definition (and 
notation):

\begin{defin}  \label{defin:pq}

Suppose $\theta$ is a $\Theta$-curve in $S^{3}$ with edges $e^{+}, 
e^{-}, e^{\bot}$.  In $H = \eta(\theta)$, denote the corresponding 
meridians by $\mpp, \mm, \mz$.  Suppose $K \subset \bdd H$ is a 
primitive curve in $\bdd H$ (i.  e.  it intersects some essential disk 
in $H$ in a single point) and $K$ intersects each of the meridians $\mpp, 
\mm, \mz$ always with the same orientation and so that some minimal 
genus Seifert surface $F$ for $K$ intersects $H$ only in $K = \bdd F$.  
Arrange the labelling and orientations of the edges and meridians so 
that, geometrically as well as algebraically,

\begin{itemize}

\item $K \cap \mm = q \geq 1$

\item $K \cap \mpp = p \geq q$

\item $K \cap \mz = p + q$.  

\end{itemize}

Then we say that $K$ (or $(K, F)$) is presented on $\theta$ as a $(p, 
q)$ quasi-cable.

\end{defin}

\bigskip

The fact that $K$ is primitive ensures that $p$ and $q$ are relatively 
prime, so $p > q$ unless $p = q = 1$.  Given $p, q$, there is a 
straightforward algorithm to describe the order in which $K$ 
intersects the three meridians (see for example \cite{OZ}): Consider a 
line in $\mathbb{R}^{2}$ of slope $p/q$ that is disjoint from the 
lattice $\mathbb{Z}^{2}$.  Choose a segment $\sss$ that projects to a 
simple closed curve in the torus $\mathbb{R}^{2}/\mathbb{Z}^{2}$.  
Then the order in which $\sss$ intersects respectively lines of the 
form $y \in \mathbb{Z}, x \in \mathbb{Z}, x + y \in \mathbb{Z}$ is the 
order in which $K$ intersects respectively the meridians $\mpp, \mm, 
\mz$.  This has the useful corollary:

\bigskip

\begin{cor}   \label{cor:ziesch}
Suppose $K$ is presented on $\theta$ as a $(p, q)$ quasi-cable, with 
$p > q \geq 2$.  Then there are at least two arcs of $K - \{ \mpp, 
\mm \}$ that are oriented from $\mpp$ to $\mm$ (and of course two then 
oriented from $\mm$ to $\mpp$).  

\end{cor}

\begin{proof} Since $q \geq 2$ the corresponding arc $\sss \subset 
\mathbb{R}^{2}$ crosses at least two vertical lines $x \in 
\mathbb{Z}$.  Since $p > q$, in between such crossings $\sss$ must 
cross at least one horizontal line $y \in \mathbb{Z}$. \end{proof}

\bigskip

In order to appreciate the point of Definition \ref{defin:pq} it's 
useful to observe that any pair $(K,F)$ can be presented as a $(p, 
q)$ quasi-cable for any relatively prime non-negative pair $(p, q)$.  
Consider the following construction.  There is a natural embedding of 
a punctured torus $T_0$ in $S^3$ with two properties:

\begin{itemize}

\item $K \subset T_0$

\item $F$ is transverse to $T_0$ with $F \cap T_0 = \bdd F$

\end{itemize}

For example, $F - \eta(K)$ is a copy of $F$ intersecting $\bdd 
\eta(K)$ in a longitudinal copy of $K$; just let $T_0$ be the 
complement of a disk in $\bdd \eta(K)$.

In the punctured torus $T_0$ choose two non-parallel normally oriented 
essential arcs $\sss^+$ and $\sss^-$.  Once $K \cap (\sss^{+} \cup 
\sss^{-})$ is minimized by isotopy, $K$ will intersect each arc 
$\sss^{\pm}$, always with the same orientation.  Indeed, given $(p, 
q)$ non-negative and relatively prime it's easy to find such arcs and 
to choose their normal orientation so that $\sss^+ \cdot K = p$ and 
$\sss^- \cdot K = q$.  One of the two choices $\sz$ for a third 
essential arc in $T_{0}$ that is not parallel to $\spp$ or $\sm$ will 
have the property that $\sz \cdot K = p + q$.  Let $H$ be the genus 
two handlebody obtained by thickening $T_0 \subset S^3$ slightly, so 
$H \cong (T_{0} \times I)$.  We can then regard $H$ as the 
neighborhood of a $\Theta$-graph $\theta \subset S^3$, with two 
vertices (one for each component of $T_0 - (\spp \cup \sm \cup \sz)$) 
and three edges $e^-, e^+, e^{\bot}$, each dual to its cognate arc.  
The natural meridians for $H$, namely $\mpp \cong \sss^+ \times I$, 
$\mm \cong \sss^- \times I$ and $\mu^{\bot} \cong \sss^{\bot} \times 
I$ are the meridians required to give $\theta$ the structure that 
presents $(K, F)$ as a $(p, q)$ quasi-cable.

So if any pair $(K, F)$ can be presented as a $(p, q)$ quasi-cable, 
what is the point of the construction?  The point will be to use 
general position between $\theta$ and a surface in $S^{3}$ as a 
short-hand way of describing in which meridians we will allow $H$ and 
the surface to intersect. If we view $H$ in this way, then any 
time the $1$-complex $\theta$ is put in general position with respect 
to a surface $S \subset S^3$, it will automatically be true that $S 
\cap H$ is a collection of meridian disks, each parallel to one of 
$\mpp, \mm, \mz$.

The first lemma may clarify the point:

\begin{lemma} \label{lemma:hexagon} Suppose the pair $(K, F)$ 
is presented as a $(p, q)$ quasi-cable on $\theta$ with both $p > q 
\geq 2$.  If there is a sphere $P \subset S^3$ 
so that
\begin{itemize}

\item $\theta$ is in general position with respect to $P$

\item $P$ intersects both of the cycles $e^{\bot} \cup e^+$ and $e^{\bot} 
\cup e^- \subset \theta.$ 

\item each component of the $1$-manifold $P \cap F$ is essential in 
$F$.

\end{itemize}

Then $genus(F) \geq 2$.

\end{lemma}

\bigskip

{\bf Remark:} One might conjecture that this is the simplest case of a 
more general result, perhaps describing how the continued fraction 
expansion of $q/p$ determines a lower bound for the genus of $F$.

\bigskip

\begin{proof} We have seen in Corollary \ref{cor:ziesch} that $K \subset 
\eta(\theta)$ at least twice switches from crossing $\mpp$ (perhaps 
repeatedly) to crossing $\mm$ (and vice versa).

To simplify the number of cases we need to consider, note first that 
we can slit open the part of $H$ corresponding to $e^{\bot} \subset 
\theta$, lengthening $e^{\pm}$ while shortening $e^{\bot}$ until 
$e^{\bot}$ is so short that it is disjoint from $P$.  So, with no loss 
of generality, we may as well assume that $P$ intersects both of the 
segments $e^{\pm}$ but not $e^{\bot}$.  

Once this is done, the components of $\theta - P$ consist of three 
types: subarcs of $e^+$, subarcs of $e^-$, and a single component 
$\theta^{\bot}$ that contains $e^{\bot}$ together with all four ends 
of the two edges $e^{\pm}.$ Then $P$ intersects each of 
the segments $e^{\pm}$ in an even number of points.  Let $\kkk^{+}$ 
denote the component of $e^{+} - P$ that is exactly half-way along 
$e^{+}$ as measured by intersections with $P$.  That is, an arc in $e^{+}$
starting from a point in $\kkk^{+}$ and ending in $\theta^{\bot}$ 
will intersect $P$ in the same number of points no matter which way 
along $e^{+}$ it runs.  Denote by $\kkk^{-}$ the analogous point in 
$e^{-}$.

The knot $K \subset H$ is similarly split 
up into segments by $P$, some parallel to segments of $e^{+} - P$, 
some to segments of $e^{-} - P$ and some lying on 
$\eta(\theta^{\bot})$.  Any subsegment of $K$ that is a union of 
components of the first (resp.  second) type will be said to be {\em 
colored $+$} (resp.  $-$.)  Remembering that $K$ is oriented, each 
segment of $K \cap \eta(\theta^{\bot})$ can be described as one of 
three types:

\begin{enumerate}

\item components of $K - P$ that run from an end of $e^+$ to an end 
of $e^-$

\item components of $K - P$ that run from an end of $e^-$ to an end of $e^+$

\item components of $K - P$ that run from an end of $e^+$ to the other 
end of $e^+$

\end{enumerate}

(There are no components of $K \cap \eta(\theta^{\bot})$ that run from 
$e^-$ to $e^-$ since $p > q$.)  We will say that these three types of 
segments of $K - P$ are colored $\omp, \opm, \opp$ respectively. The 
notation is meant to suggest that in clockwise rotation, the color 
changes from $+$ to $-$, etc.

Since all components of $F \cap P$ are essential in $F$, $F$ cannot be 
a disk.  We will suppose $genus(F) = 1$ and arrive at a contradiction.  

There are various ways that the $1$-manifold $F \cap P$ can lie in 
$F$.  If any component is a closed curve, then the complement of the 
closed curve is a simple pair of pants (i.  e.  a 3-punctured sphere) 
so all arc components of intersection must be parallel to each other.  
If there are no closed curves, the arcs of $F \cap P$ fall into (at 
most) three classes of parallel arcs in $F$.  We will assume for the 
purposes of this argument that $F \cap P$ consists of three such 
classes of parallel arcs in $F$; if there are fewer classes of 
parallel arcs, the same method works, but more easily.

Abstractly, the three families of parallel arcs of $F \cap P$ in $F$ 
give $K = \bdd F$ the structure of a hexagon $R$, in which opposite 
sides are connected via arcs of intersection that are parallel in $F$.  
See Figure \ref{fig:hexagon}.  Each end of such an arc of intersection lies in a 
meridian disk of $H$, corresponding to a point of $\theta \cap P$; if 
two ends of arcs of intersection lie in the same meridian, we 
say that the ends {\em have the same label}.  

\bigskip

\begin{figure} 
\centering
\includegraphics[width=.7\textwidth]{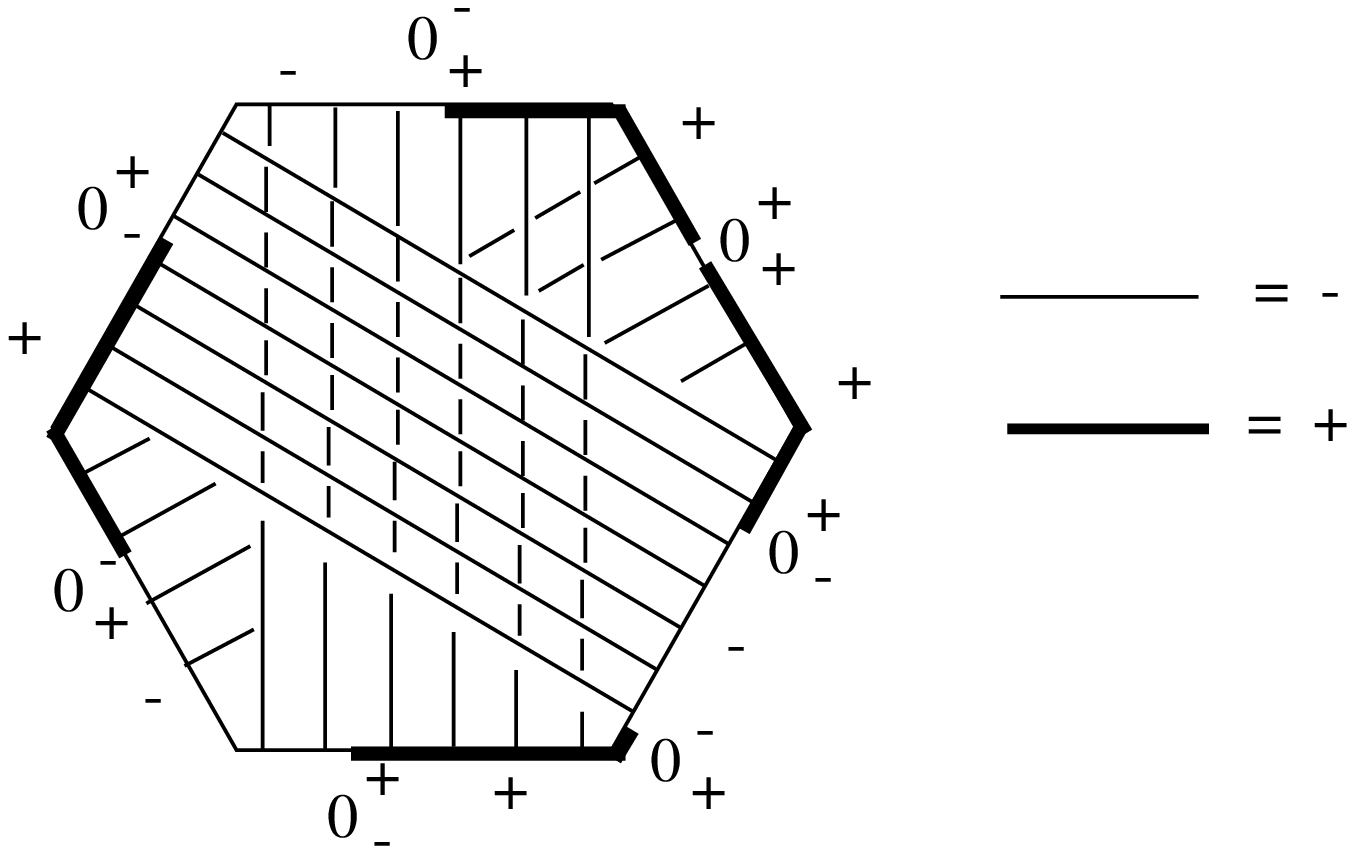}
\caption{} \label{fig:hexagon}
\end{figure}  

{\bf Claim 1:}  Opposite ends of the same arc of $F \cap P$ cannot 
have the same label.

{\bf Proof of Claim 1:} Since $K = \bdd F$ always crosses each 
meridian with the same orientation, a normal orientation induced on 
the intersection arc by a normal orientation of $P \subset S^3$ would 
have to have opposite direction at each end of the intersection arc.

\bigskip

{\bf Claim 2:} Suppose $\aaa_1$ and $\aaa_2$ are intersection arcs 
parallel in $F$, connecting opposite sides $s_1$ and $s_2$ of the 
hexagon $R$.  Suppose further that the labels of $\aaa_1$ at $s_1$ and 
$\aaa_2$ at $s_2$ are the same.  Then both are $+$-labels and all 
labels lying between the ends of the $\aaa_{i}$ on one of the $s_{i}$ 
are $+$-labels.  On the other side, between the ends of the 
$\aaa_{i}$, there is exactly one subsegment of $-$-labels.

{\bf Proof of Claim 2:} If there were a counterexample, choose 
$\aaa_1$ and $\aaa_2$ to be as close as possible (among parallel arcs 
of intersection in $F$) among all such counterexamples.  They cannot 
be the same intersection arc, by Claim 1.  We now show they cannot be 
adjacent intersection arcs in $F$.  For if they were, then the 
segments of $s_1$ and $s_2$ that lie between them would correspond to 
parallel segments of $K - P$ on $H$.  This is obvious unless the 
component of $\theta - P$ on which the segments of $s_i$ lie is 
$\theta^{\bot}$, i.  e.  the segments are colored $0$.  But even if 
the segments are colored $0$, then the fact that this is a 
counterexample forces both segments to be colored $\opm$, or both 
$\omp$ or both $\opp$, so they are in fact parallel on $\bdd H - P$.  Now 
follow a standard argument that traces its origins to \cite{GL} or 
\cite{Sc}: Consider the union of the sphere $P$, the single $1$-handle 
subsection of $H$ corresponding to the segments in the $s_i$ between 
the intersection arcs, and a $2$-handle whose core is the rectangle in 
$F$ lying between the two intersection arcs.  (See Figure \ref{fig:xcycle}.)  
This defines a Heegaard splitting of a punctured Lens space $L(2,1) = 
RP^3$, lying in $S^3$, clearly an impossibility.

\begin{figure} 
\centering
\includegraphics[width=.7\textwidth]{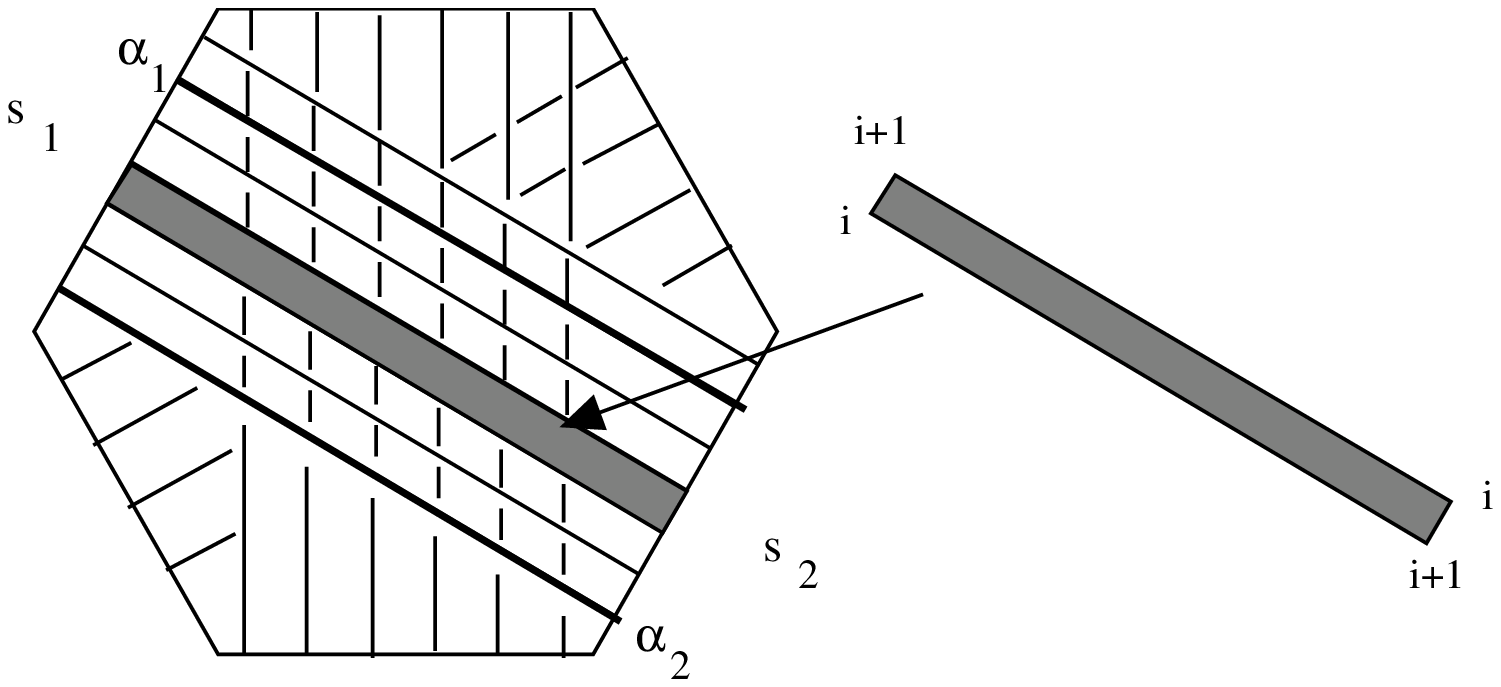}
\caption{} \label{fig:xcycle}
\end{figure} 

Since $\aaa_1$ and $\aaa_2$ are not adjacent, we can consider the 
intersection arcs $\aaa_1'$ and $\aaa_2'$ adjacent to $\aaa_1$ and 
$\aaa_2$ but closer together.  Since the intersection arcs $\aaa_1'$ 
and $\aaa_2'$ are not a counterexample, either the $\aaa_i'$ have 
different labels at $s_{i}, i = 1, 2$ or they have the same labels but 
the count of segments colored $\opm$ or $\omp$ between them changes.  
Either outcome is only possible if the segments $k_i$ between $\aaa_i$ 
and $\aaa_i'$ on $s_i$ are colored $0$, for both $i = 1, 2$ and at 
least one, say $k_1$, is colored $\opm$ or $\omp$, say $\opm$.  Then 
$k_2$ is colored either $\opp$ (and the labels of the $\aaa_i'$ at 
$s_i$ are different for $i = 1, 2$) or $k_2$ is also colored $\opm$.  
Consider first the former case, $k_2$ is colored $\opp$.  This 
immediately implies that both of the initial labels were $+$-labels.  
If another segment colored $\opm$ lies between $\aaa_{1}$ and 
$\aaa_{2}$ on $s_{2}$ then an intersection arc adjacent to it would be 
a counterexample closer to $\aaa_1$ than $\aaa_2$ is.  So we deduce 
that all labels on $s_{2}$ between the $\aaa_{i}$ are $+$-labels.  
Because $\aaa_1$ and $\aaa_2$ are a counterexample, another segment 
colored $\opm$ must lie between $\aaa_1$ and $\aaa_2$ on $s_{1}$.  
Then there is also one colored $\omp$ between $\aaa_1'$ and $\aaa_2'$.  
In order for an intersection arc adjacent to it and the intersection 
arc $\aaa_2'$ not to be a counterexample, there must be a segment on 
$s_1$ even closer to $\aaa_2'$ that is colored $\opm$.  Between it and 
the segment colored $\omp$ lies every label colored $+$.  Opposite to 
this segment on $s_2$ every label is colored $+$, since no $\opm$ or 
$\omp$ color appears on $s_2$ between $\aaa_1$ and $\aaa_2$.  So for 
any of these intersection arcs, whatever the label is on $s_2$, 
there's a parallel intersection arc with that label on $s_1$ and 
between them lies no label $\opm$ or $\omp$.  This creates a closer 
together pair of intersection arcs that are a counterexample, a 
contradiction.  See Figure \ref{fig:side1}.

If instead $k_2$, like $k_1$, is colored $\opm$, then the labels of 
the $\aaa_i'$ at $s_i$ are the same for $i = 1, 2$ and, since these 
arcs are not a counterexample, they must be $+$-labels, all the 
labels on one side between the $\aaa_{i}'$ must be $+$-labels, but on 
the other side there must be a switch to $-$-labels.  But this 
produces exactly the same contradiction as before: all the labels on 
one side are $+$-labels and an entire sequence of $+$-labels appears 
on the other side.  Hence there can be no counterexamples, 
completing the proof of claim 2.

\begin{figure} 
\centering
\includegraphics[width=.5\textwidth]{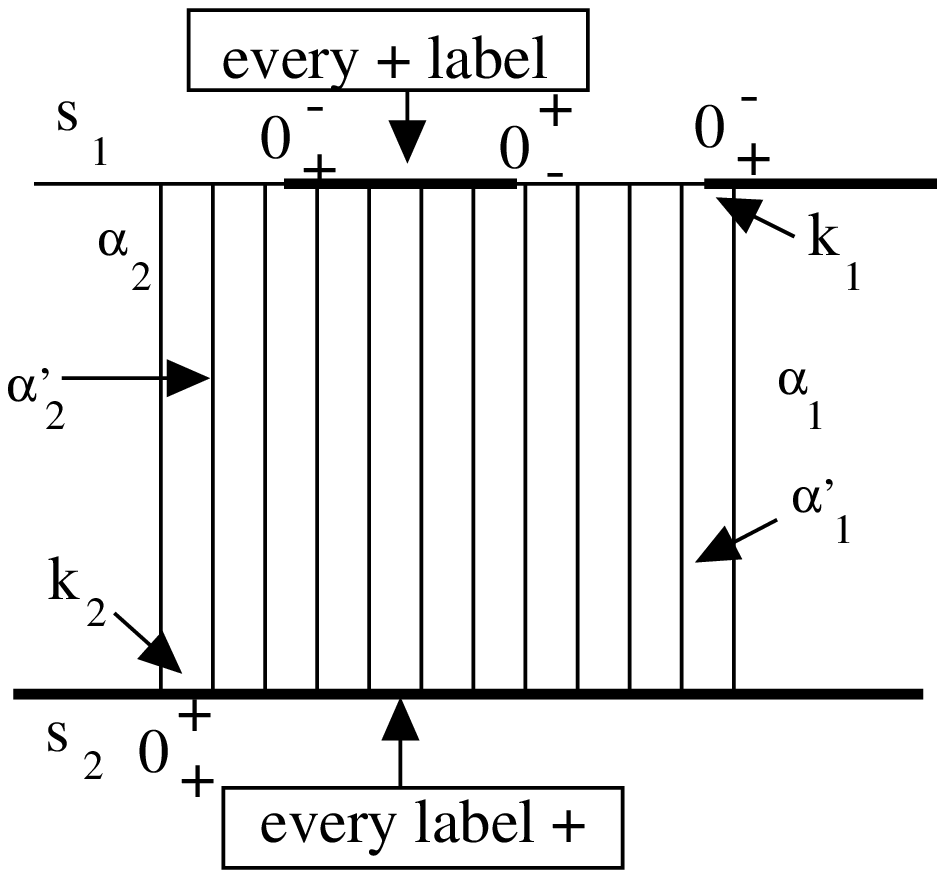}
\caption{} \label{fig:side1}
\end{figure}  

\bigskip

{\bf Claim 3:} Segments corresponding to $\kkk^{+}$ (or $\kkk^{-}$) cannot appear 
on opposite sides of $R$ (including corners).

{\bf Proof of Claim 3:} If $\kkk^{-}$ appeared on opposite sides, it 
would immediately contradict Claim 2.  Suppose $\kkk^{+}$ appeared on 
opposite sides of $R$.  Choose those occurences that are closest 
together (as measured by arcs between them) and let $\aaa_{1}$ and 
$\aaa_{2}$ denote the arcs adjacent to those occurences of $\kkk^{+}$ 
that are closer together.  Let $m = |e^{+} \cap P|$.  If the number of 
arcs between the $\aaa_{i}$ is less than $m/2$ it would contradict 
Claim 2, since neither side could then have any $-$-labels.  If the 
number of arcs between them is no less than $m-2$, then the fact, from 
Claim 2, that one side consists entirely of $+$-labels would ensure 
that another label $\kkk^{+}$ occurs even more closely to the opposite 
$\kkk^{+}$, contradicting the choice of $\kkk^{+}$ segments that are 
closest together.  Suppose finally that there are between $m/2$ and 
$m-2$ arcs between them.  Then, by Claim 2, on one side between them 
will be a segment colored $\opm$ or $\omp$ and on the other a segment 
colored $\opp$.  Moreover, the arcs $\aaa_{i}'$ adjacent to these 
segments and closer together would have the same $+$-label but would 
have no $-$-labels between them on either side, again contradicting 
Claim 2.

\bigskip

{\bf Claim 4:}  Two segments, one corresponding to each of 
$\kkk^{\pm}$ cannot occur on the same side of $R$ (including corners.)

{\bf Proof of Claim 4:} Let $n = |e^{-} \cap P|$.  If both types of 
$\kkk$ occur on one side, say $s_{1}$ then that side is incident to at 
least $(m+n)/2$ arcs and contains a segment of type $\opm$ or $\omp$.  
Also, from Claim 3, the opposite side $s_{2}$ can be incident to 
neither type of $\kkk$-interval.  Since the $s_{i}$ are incident to 
the same number of arcs, this means that $s_{2}$ contains some segment 
of the form $\opm, \omp$ or $\opp$, say $\opm$.  In fact, it must be 
of the form $\opp$ since otherwise its having length $> (m+n)/2$ would 
force a $\kkk$ segment to appear.  For the same reason, the $\opp$ 
segment in $s_{2}$ must be opposite a segment in $s_{1}$ that lies 
between a $\kkk_{+}$ label and the $\opm$ label.  The arcs $\aaa_{i}$ 
adjacent to the $\opm$ and $\opp$ segments and closer together would 
have the same $+$-label but would have no $-$-labels between them on 
either side, again contradicting Claim 2.

\begin{figure}
\centering
\includegraphics[width=.4\textwidth]{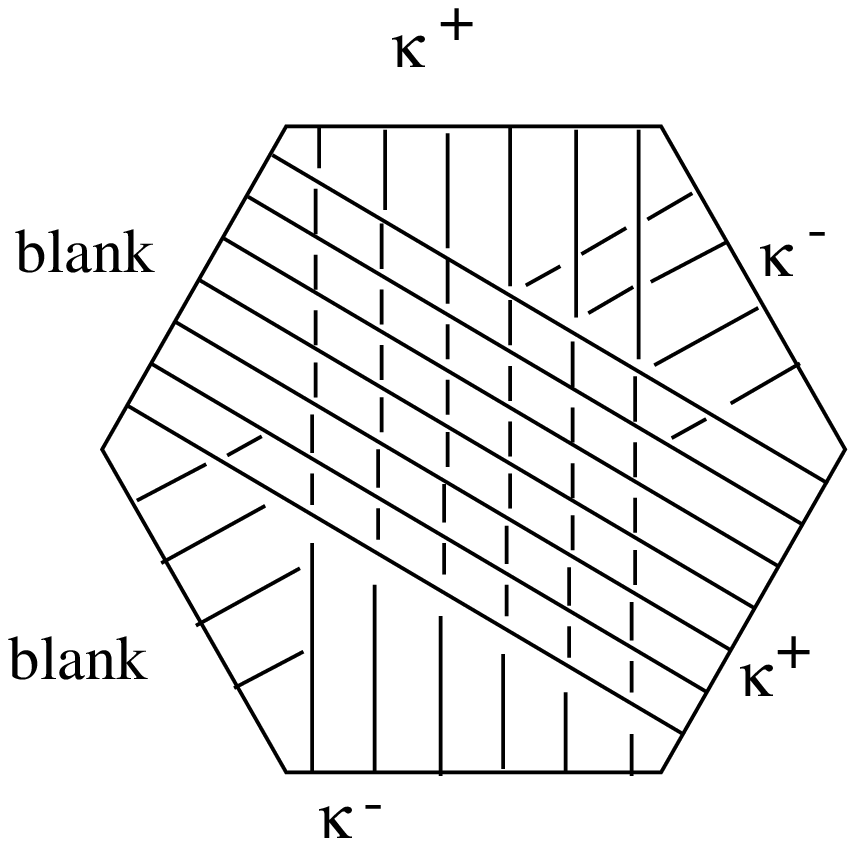}
\caption{} \label{fig:kappa}
\end{figure}

Following Claim 4, we can think of each side of $R$ as either a 
$+$-side, a $-$-side or blank, depending on whether a copy of 
$\kkk_{+}, \kkk_{-}$ or no $\kkk$ at all appears..  Moreover, by 
Claim 3 opposite 
sides can't both be $+$-sides or both be $-$-sides and (a crucial 
point) following Corollary \ref{cor:ziesch}, at least $4$-sides have 
signs, alternating around $R$ as $+,-,+,-$.  Combining these facts, the 
only possible signing of the sides of $R$ is, in order, $blank, 
blank,+,-,+,-$ (with some orientation of $\bdd R$).  Now the fact 
that the two adjacent blank sides have no sign whereas their adjacent sides 
have different signs means that the total number of arcs intersecting 
those two blank sides must be less than $(m + n)/2$.  On the other 
hand, the adjacent sides opposite these blank sides have signs $+$ and $-$ 
which guarantees that their combined length is greater than 
$(m+n)/2$.  The contradiction proves the lemma.  See Figure 
\ref{fig:kappa}.
\end{proof}

\bigskip

The requirement that $p, q \geq 2$ in the above lemma is central to the 
proof, of course, since it guarantees the repetitions in patterns 
around $\bdd F$ that lead to the combinatorial contradiction.  
Nonetheless, there are important situations in which the results of 
Lemma \ref{lemma:hexagon} hold true even when $q = 1$.

\begin{figure}
\centering
\includegraphics[width=.6\textwidth]{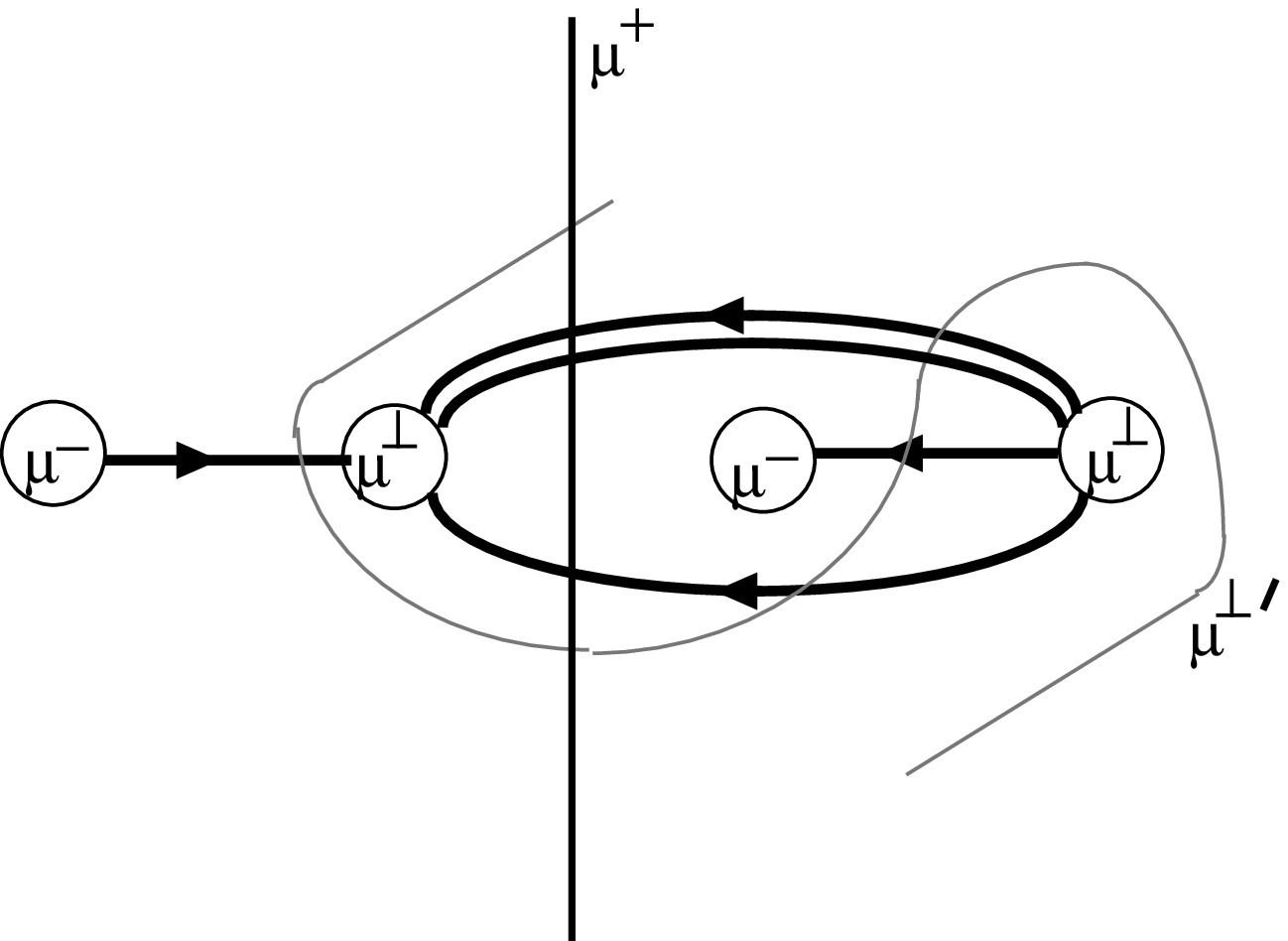}
\caption{} \label{fig:Whitney}
\end{figure}

The easiest example to see requires a preliminary construction.  
Suppose $(K, F)$ is presented as a $(p, q)$ quasi-cable on $\theta$ 
with regular neighborhood $H$, and and consider the $4$-punctured 
sphere $\Ss$ obtained from $\bdd H$ by removing copies of the 
meridians $\mz$ and $\mm$.  Then the boundary of $\Ss$ consists of two 
copies of $\bdd \mz$ and two copies of $\bdd \mm$, one on each side of 
the circle $\bdd \mpp \subset \Ss$.  The arcs of $K \cap \Ss$ are 
oriented to flow from one side of $\bdd \mpp$ to the other.  The 
slopes of these arcs naturally define another circle $\mz' \subset 
\Ss$ with the property that each arc of $K \cap P$ intersects $\mz'$ 
exactly once and $|\mz' \cap \mpp| = 2$.  See figure 
\ref{fig:Whitney}.  (Actually, there are two candidates for $\mz'$; 
the other is obtained by vertical reflection.)  If we discard the 
meridian $\mpp$ of $H$ and replace it with $\mz'$ then the associated 
$\Theta$-graph $\theta'$ is one obtained from $\theta$ by a Whitney 
move.  $(K, F)$ is still presented as a quasi-cable on $\theta'$ but 
now of type $(p + q, q)$.  With the new structure, the old $\mp$ is 
discarded, the old $\mm$ becomes also the new $\mm'$ and the old 
$\mz$ becomes the new $\mpp'$.

There is a similar move in which $H$ is cut up along $\mpp$ and 
$\mz$, the meridian $\mm$ is discarded and the new graph $\theta'$ 
presents $(K, F)$ as a $(p + q, p)$ quasi-cable.

\begin{defin} \label{defin:Whitney} The moves on $\theta$ just 
described, which change the presentation of $(K, F)$ from that of a 
$(p, q)$ quasi-cable to, respectively, a $(p + q, q)$ quasi-cable or a 
$(p + q, p)$ quasi-cable are called {\em standard Whitney moves} on 
$e^{+}$ and $e^{-}$ respectively.
\end{defin}

\begin{lemma} \label{lemma:hexagon1} Suppose the pair $(K, F)$ 
is presented as a $(p, q)$ quasi-cable on $\theta$ with $p > q 
\geq 1$.  If there is a sphere $P \subset S^3$ 
so that
\begin{itemize}

\item $\theta$ is in general position with respect to $P$

\item $P$ intersects both of the edges $e^{\bot}$ and $e^+$ but is 
disjoint from $e^-$ 

\item each component of the $1$-manifold $P \cap F$ is essential in 
$F$.

\end{itemize}

Then $genus(F) \geq 2$.

\end{lemma}

\begin{proof} Perform a standard Whitney move on $e^{-}$ so that 
afterwards the new $\Theta$-graph $\theta'$ presents $(K, F)$ as a $(p 
+ q, p)$ quasi-cable.  Since $P$ was disjoint from $e^{-}$ this has 
no effect on $P \cap F$ and, since $p+q \geq p \geq 2$ Lemma 
\ref{lemma:hexagon} applies.  
\end{proof}

In a related but more complicated case, the combinatorics 
is so close to that of an actual cable knot that the arguments are 
considerably easier than those above.

\begin{lemma} \label{lemma:hexagon2} Suppose the pair $(K, 
F)$ is presented as a $(p, q)$ quasi-cable on $\theta$, $p \geq q \geq 
1$.  If there is a 
sphere $P \subset S^3$ so that
\begin{itemize}
\item $\theta$ is in general position with respect to $P$
\item $P$ intersects $e^{\bot}$ but is disjoint from $e^{\pm}$ and
\item each component of the $1$-manifold $P \cap F$ is essential in 
$F$.
\end{itemize}
Then $genus(F) \geq 2$.
\end{lemma}

\begin{proof} Suppose, as in the proof of Lemma \ref{lemma:hexagon}, 
$genus(F) = 1$.  The case $q \geq 2$ is already settled, so we assume 
$q=1$.  Set $n = |e^{\bot} \cap P|$, necessarily even; then $|K \cap 
P| = n(p+1) \geq 4$ and the number of edges of $F \cap P$ is 
$n(p+1)/2$.  Of the intervals in $K - P$, $(p+1)(n-1)$ lie on 
segments of $e^{\bot} - P$, $p$ contain $e^+$ in their interior, and 
exactly one contains $e^-$.  Call the last the ``special'' component.  
We now consider how these intervals are distributed around the hexagon 
$R$ described in the proof of Lemma \ref{lemma:hexagon2} above.  
First note that he 
two vertices corresponding to the corners of the hexagon (when 
reidentified to give $F$) either lie on the same or opposite sides of 
$P$, depending on whether $n(p+1)/2$ is even or odd, and this in turn 
determines whether the number of intersections of any edge of $R$ with 
$P$ is even or odd.  The upshot is that every edge of $R$ intersects $P$ with 
the same parity and that parity is determined by the parity of 
$n(p+1)/2$.  

The combinatorial argument gets easier as $p$ gets larger, since 
intervals cooresponding to the same components of $\theta - P$ appear 
more often.  So for brevity we just do the case $p = 1, 2$ (and of 
course $q = 1$) and merely outline the argument.  (In fact the 
argument for these cases also instantly gives as well all cases in 
which $p+1$ is a multiple of $2$ or $3$.)

When $p = 1$ there are $2n$ segments in $\bdd F$ and it follows that, 
except perhaps for the special component, centers of opposite edges in 
$R$ represent the same segment of $e^{\bot}$, for there are exactly as 
many intersection points with $P$ (namely $n$) going one way around 
$K$ between them as the other.  At most one of these three opposite 
pairs contains the special component, so the other two display lens 
spaces $L(2, 1)$ in $S^3$, a contradiction.  The only way to avoid 
this contradiction is if only two (opposite) sides of $R$ intersect 
$K$, so all arcs in $F \cap P$ are parallel in $F$, and the special 
component appears in the center of this band of parallel arcs.  This 
only transfers the contradiction: if we let $A$ be the annulus which 
is the complement in $F$ of the single band containing all the 
parallel arcs of $P \cap F$ then it is easy to see that $\bdd A$ lies 
on the boundary $T$ of the punctured solid torus obtained by attaching 
a single component of $\eta(e^{\bot}) - P$ to $P$.  Moreover, the 
components of $\bdd A$ are oriented so that they form parallel circles 
in $T$.  This implies that $A$, together with an annulus in $T$, form 
a Klein bottle in $S^3$, another (related) contradiction.  

When $p = 2$ then $|K \cap P| = 3n$ so the number of edges in $K \cap 
P$ is $3n/2 \geq 3$.  There are two cases:  

If each edge of $R$ 
intersects $P$ the same number of times (as our parity discussion
guarantees will happen when $n = 2$), then each triple of corners 
of the hexagon (corresponding to a vertex of the punctured torus when 
it's reassembled) represent the same interval of $e^{\bot}$, or 
perhaps the special segment.  At least one of the triples doesn't 
contain the special segment, and the hexagon in $F$ cut off from $F$ by the 
arcs of $P \cap F$ adjacent to these corners, together with $P$ and 
the segment of $e^{\bot}$ the hexagon is incident to describe the 
spine of a Lens space $L(3,1) \subset S^{3}$, a contradiction.  

If some edges of $R$ intersect $P$ more often than others (so $n \geq 
4$ and so $K \cap P \geq 12$), then 
consider a longest pair $\rrr_{i}, i = 1, 2$ of opposite edges of $R$ 
(length here is shorthand for the number of points of intersection 
with $P$).  Picturing these opposite sides as the top and bottom of 
the hexagon, consider the distance in $R$ between the left ends of 
the $\rrr_{i}$.  Our choice of edge guarantees that it is less than 
$1/3$ the circumference of $R$, that is $K \cap P$.  Similarly for the 
right hand ends.  It follows readily that there are at least two 
rectangles in $F$, cut off by a pair of adjacent arcs of $P \cap F$ 
running between the $\rrr_{i}$ that either make up part of the spine 
of an $L(2,1) \subset S^{3}$ or contain the special component.  At 
most one can contain the special component, leading to the same 
contradiction.  
\end{proof}

The same argument, with the roles of $e^+$ and $e^{\bot}$ 
switched and $p$ replacing $p+1$, shows

\begin{lemma} \label{lemma:hexagon2b} Suppose the pair $(K, 
F)$ is presented as a $(p, q)$ quasi-cable on $\theta$, $p > q \geq 
1$.  If there is a sphere $P \subset S^3$ so that
\begin{itemize}
\item $\theta$ is in general position with respect to $P$
\item $P$ intersects $e^{+}$ but is disjoint from 
$e^{\bot}$ and $e^-$ and
\item each component of the $1$-manifold $P \cap F$ is essential in 
$F$.
\end{itemize}
Then $genus(F) \geq 2$.
\end{lemma}

Essentially the same argument applies in a slightly different setting:

\begin{lemma} \label{lemma:hexagon3}  Suppose the pair $(K, 
F)$ is presented as a $(p, q)$ quasi-cable on $\theta$.  If there is a 
sphere $P \subset S^3$ so that
\begin{itemize}
\item $\theta$ is in general position with respect to $P$
\item $P$ intersects each of $e^{\pm}$ in a single point and
\item each component of the $1$-manifold $P \cap F$ is essential in 
$F$.
\end{itemize}
Then either also $|e^{\bot} \cap P| = 1$ or $genus(F) \geq 2$.
\end{lemma}

\begin{proof} Here the intervals of $K - P$ that are incident to $P 
\cap e^{-}$ constitute two adjacent ``special'' components.  Since we 
can take $|e^{\bot} \cap P| \geq 3$ (it's necessarily odd, since the 
hypothesis guarantees that its ends lie on opposite sides of $P$) then 
$K \cap P \geq 8$.  Then the combinatorial arguments of Lemma 
\ref{lemma:hexagon2} applies with little change.
\end{proof}

We can switch the roles of $e^+$ and $e^{\bot}$ in the above 
proofs, at the cost of raising $p$ to $2$:

\begin{lemma} \label{lemma:hexagon4} Suppose the pair $(K, F)$ is
presented as a $(p, q)$ quasi-cable on $\theta$ with 
$p > q \geq 1$.  If there is a sphere $P
\subset S^3$ so that

\begin{itemize}

\item $\theta$ is in general position with respect to $P$

\item $P$ intersects $e^{\bot}$ and $e^-$ in a single point and

\item each component of the $1$-manifold $P \cap F$ is essential in 
$F$.

\end{itemize}
Then either also $|e^{+} \cap P| = 1$ or $genus(F) \geq 2$.
\end{lemma}

\begin{proof} We need only consider the case $q = 1$.  Let $n = |e^{+}
\cap P|$, necessarily odd.  Then $|K \cap P| = p(n+1) + 2$.  In fact,
$p(n - 1)$ segments lie parallel to segments of $e^+ - P$, $2p+2$
segments (in sequential pairs) are incident to the point $e^{\bot}
\cap P$ and one of these sequential pairs (called the ``special
pair'') are incident to the point $e^-$.  Much as before, we set $p
= 2, 3$ as the most difficult but also roughly representative cases.

When $p = 2$ there are an odd number of arcs in $P \cap F$, since $|K 
\cap P|/2$ is odd and so there are an odd number of arcs incident to 
each edge.  The extra ``special'' segments in $K$ mean that the 
intersection arc of $K \cap F$ that lies in the center of each pair of 
oppsoite sides of $R$ does not have its ends at the same point of 
intersection of $\theta$ with $P$ (an immediate contradiction via its 
normal orientation) but it does mean that adjacent to each such 
central arc, at opposite ends, are segments of $K - P$ that are 
parallel to the same segment of $\theta - P$.  This would exhibit, as 
usual, the absurd $L(2, 1) \subset S^3$.  This contradiction is only 
avoided if each side of $R$ intersects $P$ in exactly one point.  But 
this means that $|K \cap P| = 6$, so $n = 1$ as required.

When $p = 3$ then, since the number of arcs in $P \cap F$ is $1 mod 
3$, not all edges of $R$ intersect $P$ the same number of times.  We 
may as well restrict to the case $n \geq 3$ so there are at least $7$ 
arcs in $F \cap P$.  If there is a single pair of longest edges, each 
must then intersect $P$ at least $5$ times.  It's easy to see in this 
case, that wherever the adjacent special edges lie, they cannot 
disrupt the existence of at least one pair of opposite intervals of 
$K \cap P$ (among those lying in these longest sides) that correspond 
to the same component of $e^{+} - P$ and so constitute part of a Lens 
space $L(2, 1) \subset S^{3}$.  Similarly, if there are two pairs of 
opposite longest sides, then each must be of length at least $3$ and, 
wherever the adjacent special components lie, they cannot disrupt the 
existence of a similar Lens space contradiction from at least one 
pair of opposite longest sides.
\end{proof}

\section{Knots thinly presented on handlebodies}

In \cite[Section 2]{GST}, we extended Gabai's notion of thin position 
for knots to include also certain types of graphs in $3$-space.  We 
briefly review (and incidentally somewhat extend) that development 
here, since it will be an important ingredient of our argument.

Choose a height function $h:S^3-\{x,y\}=S^2\times \mathbb R 
\rightarrow\mathbb R$ and let $P(t)=h^{-1}(t)$.

\begin{defin}
A finite trivalent graph $\Ggg \cup S^3-\{x,y\}$ is in {\it normal 
form} with respect to $h$ if
\begin{enumerate} 
\renewcommand{\labelenumi}{(\alph{enumi})} 

\item For each edge $e 
\subset \Ggg$ the critical points of $h|e$ are nondegenerate and lie 
in the interior of $e$, 

\item The critical points of $h|edges$, and the vertices of 
$\Ggg$, all occur at different heights.

\item At each (trivalent) vertex $v$ of $\Ggg$ either two ends of 
incident edges lie above $v$ (we say $v$ is a {\em $Y$-vertex}) or two 
ends of incident edges lie below $v$ (we say $v$ is a {\em $\lll$-vertex})
\end{enumerate}
\end{defin}

Standard Morse theory shows that any finite trivalent graph in $S^3$ can be 
infinitesimally isotoped so that it is in normal form.

\begin{defin} The {\em maxima} of $\Ggg$ consist of all local 
maxima of $h|edges$ and all $\lambda$-vertices.  Similarly, the {\em 
minima} of $\Ggg$ consist of all local minima of $h|edges$ and all 
$Y$-vertices.  A maximum (resp.  minimum) that is not a 
$\lambda$-vertex (resp.  $Y$-vertex) will be called a {\em regular} 
maximum (resp.  minimum).  The union of the maxima and minima (hence 
including the vertices) are called the {\em critical points} of $\Ggg$ 
and their heights the {\em critical values} or {\em critical heights}.
\end{defin}

\begin{defin} Let $t_0 < \ldots < t_n$ be the successive critical 
heights of $\Ggg$ and suppose $t_{v_1}, \ldots, t_{v_j}$ are that 
subset of levels at which vertices occur.  Let $s_i, 1 \leq i \leq n$ 
be generic levels chosen so that $t_{i-1} < s_i < t_i$.  Define the 
width of $\Ggg$ to be $$W(\Ggg) = 2(\Sigma_{i \notin v_1, \ldots, 
v_j}|P(s_i)\cap (\Ggg)|) + (\Sigma_{i \in v_1, \ldots, v_j}|P(s_i) 
\cap (\Ggg)|).$$
\end{defin}

\begin{defin} 
A {\it thin position} of a graph $\Ggg \subset S^3$ is a normal form 
(with respect to $h$) which minimizes the width of $\Ggg$.
\end{defin}

{\bf Remark:}  In practice, the chief property of a thin positioning 
of a graph $\Ggg$ that we will need is this:  The positioning becomes thinner 
if a maximum is pushed below a minimum, but the width is unaffected 
by pushing one maximum above or below another maximum, or one minimum above 
or below another minimum.  See \cite[Section 3]{GST} for details.

A graph $\Ggg \subset S^3$ in normal form with respect to $h$ can be 
thickened slightly to give a solid handlebody $\eta(\Ggg) \subset 
S^3$.  Standard techniques allow us to take a neighborhood so thin 
that the height function $h|\bdd(\eta(\Ggg))$ has the obvious Morse 
structure: very near any regular maximum (resp.  minimum) of $\Ggg$ 
there are two non-degenerate critical points of $h|\bdd(\eta(\Ggg))$, 
one a saddle just below (resp.  above) and one a maximum (resp.  
minimum) just above (resp.  below).  Similarly, just above (resp.  
below) a $Y$-vertex (resp.  $\lll$-vertex) there is a single saddle 
singularity.  When we refer to a regular neighborhood of $\eta(\Ggg)$ 
of $\Ggg$ we will always mean a thickening with this property.  
Slightly abusing notation, $S^3 - \eta(\Ggg)$ will denote the closed 
complement of $\eta(\Ggg)$.  We will be concerned with simple closed 
curves on $\bdd\eta(\Ggg)$ and with properly 
imbedded surfaces in $S^3 - \eta(\Ggg)$.

\begin{defin} Suppose $\Ggg$ is a graph, in normal form with respect 
to $h$, and $K \subset \bdd\eta(\Ggg)$ is a simple closed curve.  Then 
$K$ is in normal form on $\bdd\eta(\Ggg)$ if each critical point of 
$h$ on $K$ is non-degenerate, and occurs near an associated critical 
point of $\Ggg$ in $\bdd\eta(\Ggg)$.  Furthermore, the number of 
critical points of $K$ has been minimized via isotopy of $K$ in 
$\bdd\eta(\Ggg)$.
\end{defin}

\begin{defin} \label{def:normal} A properly imbedded surface 
$$(F,\bdd F) \subset (S^3 - \eta(\Ggg), \bdd\eta(\Ggg))$$ is in {\it 
normal form} if
\begin{enumerate}

\item each critical point of $h$ on $F$ is nondegenerate, 

\item $\bdd F$ is in normal form with respect to $h$

\item no critical point of $h$ on int$(F)$ occurs near a 
critical height of $h$ on $\Ggg$, \label{near}

\item no two critical points of $h$ on int$(F)$ occur at the same height, 

\item the minima (resp.  maxima) 
of $h|\partial F$ at the minima (resp.  maxima) of $\Ggg$ are 
also local extrema of $h$ on $F$, i.e., `half-center' singularities, 

\item the maxima of $h|\partial F$ at $Y$-vertices and the minima of 
$h|\partial F$ at $\lambda$-vertices are, on the contrary, 
`half-saddle' singularities of $h$ on $F$.
\end{enumerate}
\end{defin}

{\bf Remark:} The meaning of ``near'' in (\ref{near}.) is probably best thought 
of informally, but the technical requirement (for, say, the 
critical height of a maximum $v$ of $\Ggg$) is this: No critical point of $h$ 
on the interior of $F$ occurs at a height between the levels of the 
maxima of $\bdd F$ (if any) near $v$ and the level of the saddle point 
of $\bdd\eta(\Ggg)$ near $v$.  Standard Morse theory ensures that, 
for $\Ggg$ in normal form, any properly imbedded surface $(F, \partial 
F)$ can be put in normal form.

\begin{defin} \label{defin:bridge}
$\Ggg$ is in {\it bridge position} if there is a level sphere, called 
a {\em dividing sphere} for the bridge position, that lies above all 
minima of $\Ggg$ and below all maxima.  \end{defin}

\begin{defin} 
Given $\Ggg$ in normal form and $P$ a level sphere for $h$ at a 
generic height, let $B_u$ and $B_l$ denote the balls which are the 
closures of the region above $P$ and below $P$ respectively.  An {\it 
upper disk} (resp.  {\it lower disk}) for $P$ is a disk $D \subset S^3 
- \eta(\Ggg)$ transverse to $P$ such that $\partial D = 
\alpha\cup\beta$, where $\alpha$  is an arc imbedded on 
$\partial\eta(\Ggg)$, $\beta = \bdd D \cap P$ is an arc properly 
imbedded in $P - \eta(\Ggg)$, $\partial\alpha=\partial\beta$ and a small 
product neighborhood of $\bdd D$ lies in $B_u$ (resp.  $B_l$) i.e., it 
lies {\it above} (resp.  {\it below}) $P$.  \end{defin} 

Note that $\inter(D)$ may intersect $P$ in simple closed curves.  An 
innermost such simple closed curve cuts off a disk that lies either 
above or below $P$.  Such a disk is called an {\em upper cap} or {\em 
lower cap}.  For the moment, these caps will be unimportant.

A natural occurence of upper (or, symmetrically, lower) disks is this:  
According to Definition \ref{def:normal}, a maximum of $\bdd F$ 
near a maximum of $\Ggg$ is a half-center singularity on $\bdd F$.  In 
particular, a sphere $P$ just below this maximum will cut off an upper 
disk from $F$.

\begin{defin} Suppose $\Ggg$ is a graph, in normal form with respect 
to the height function $h$, and $K$ is a normal form simple closed 
curve on $\bdd\eta(\Ggg)$.  If $K$ is also in thin position (as a knot 
in $S^3$) with respect to $h$, then we say that $K$ is {\em thinly 
presented on $\eta(\Ggg)$}.
\end{defin}

\begin{lemma} \label{lemma:essential}
Suppose $\Ggg$ is a graph, in normal form with respect to the height 
function $h$, $F \subset S^3 - \eta(\Ggg)$ is an incompressible 
surface in normal form and $K = \bdd F$ is thinly presented on 
$\eta(\Ggg)$.  Suppose a maximum and a minimum of $K$ occur 
respectively at heights $u$ and $l$ with $l < u$.  Then there is a 
generic level sphere $P = P(t)$, $l < t < u$ so that every arc 
component of $P \cap F$ is essential in $F$.
\end{lemma}  

\begin{proof} We have seen that a sphere just below $P(u)$ cuts off an 
upper disk from $F$ and that a sphere just above $P(l)$ cuts off a 
lower disk from $F$.  The seminal point of thin position (see 
\cite{G}) is that there cannot simultaneously (even at a critical 
point of $h$ on the interior of $F$) be both an upper and a lower 
disk, for these disks could be used to push a maximum of $K$ below a 
minimum, thinning $K$.  Hence there is a generic height $t$ between 
$l$ and $u$ for which the level sphere $P = P(t)$ cuts off neither an 
upper nor a lower disk from $F$.  But this means there can be no arcs 
of $P \cap F$ which are inessential, for an outermost such inessential 
arc would cut off either an upper or a lower disk from $F$.  
\end{proof}

Combining Lemma \ref{lemma:essential} with the central lemma of the 
previous section, we have this corollary:

\begin{cor} \label{cor:thinpq} Suppose the pair $(K, F)$ is 
thinly presented 
as a $(p, q)$ quasi-cable on $\theta$ with both $p, q 
\geq 2$. 

Then $genus(F) \geq 2$.
\end{cor}

\begin{proof} Let $M^+$ and $M^-$ be the highest maxima of, 
respectively, the cycles $e^{\bot} \cup e^+$ and $e^{\bot} \cup e^-$.  
Similarly, let $m^+$ and $m^-$ be the respective lowest minima of 
these cycles.  Let $u = min \{M^+, M^-\}$ and $l = max\{m^+, m^-\}$.  
Since $e^{\bot}$ is in both cycles, we know that any point on 
$e^{\bot}$ lies below $u$ and above $l$, so $l < u$.  Since $p, q > 
0$, $u$ and $l$ are (near) the heights of, respectively, maxima and 
minima of $K$.  Choose a level sphere as in Lemma 
\ref{lemma:essential} between $u$ and $l$.  By construction, such a 
level sphere lies at a height between the maximum and minimum of each 
of the cycles $e^{\bot} \cup e^+$ and $e^{\bot} \cup e^-$ and so 
intersects both of them.  Now apply Lemma \ref{lemma:hexagon}.
\end{proof}

\begin{lemma} \label{lemma:bridgepq} Suppose $K$ is a knot 
which is thinly presented on a as a $(p, q)$ quasi-cable on $\theta \subset S^3$.  
Suppose $\theta$ has been made as thin as possible subject to this 
condition.  Suppose furthermore that the complement $S^3 - 
\eta(\theta)$ is a genus two handlebody.

Then $\theta$ is in bridge position.
\end{lemma}

\begin{proof} The proof is analogous to that of \cite[Proposition 
4.4]{GST}.  Suppose that $\theta$ is not in bridge position.  Then 
there is a level sphere $P$ that lies between a sequential pair of 
critical levels for $\theta$, a maximum just below $P$ and a minimum just above 
$P$.  Maximally compress $P - \theta$ in the complement of $\theta$.  
The resulting meridional planar surface $\Pt$ is incompressible in the 
handlebody $S^3 - \eta(\theta)$ so each component is parallel to a 
subsurface of $\bdd\eta(\theta)$ (see \cite{Mo}).  Since the boundary 
components of $\Pt$ are meridians of $\eta(\theta)$, each component of 
$\Pt$ can be completed to a sphere in $S^3$, and the piece of $\theta$ 
lying in one of the balls bounded by that sphere is an unknotted tree 
(possibly just an arc) in the ball.  

For concreteness, choose an innermost such ball $B$ and suppose it 
lies above $\Pt$.  Then every arc of $K$ lying in $B$ has at least one 
maximum in $B$.  If $\theta \cap B$ is a single arc $\aaa$, then 
isotope that arc to lie in $\Pt$.  The arc in $\Pt$ can be chosen to 
be disjoint from those disks which are the results of the compressions 
that created $\Pt$ from $P$, so in fact then $\aaa$ lies in $P$.  
After this isotopy, the width of $\theta$ is reduced (since a maximum 
has been pushed below whatever minimum lay just above $P$) and that of 
$K$ is not increased (since each arc of $K \cap B$ still has at most 
one maximum.)  This argument shows more: any arc of $K \cap B$ must 
have {\em exactly} one maximum in $B$, for otherwise a disk of 
parallelism between that arc and $\Pt$ (guaranteed by \cite{Mo}) could 
be used to reduce the number of critical points on $K$, contradicting 
the assumption that $K$ is in thin position, hence in minimal bridge 
position.  

If the tree $\theta \cap B$ contains only maxima 
(including perhaps $\lll$-vertices), then 
pushing one to $\Pt$ would push it below the minimum that we know 
lies (elsewhere) just above $\Pt$, again thinning $\theta$.  So we 
know that $\Pt$ contains at least one minimum and, since it can't be a 
minimum of $K$, it must be a $Y$-vertex, with ends of $e^{\pm}$ 
descending into it.  Similarly, if both vertices are in $B$ then both 
must be $Y$-vertices. Since $\theta \cap B$ is a tree, at most one of 
$e^{\pm}$ lies in $B$ so we can assume that, say, $e^{-}$ intersects 
$\bdd B$.  Also $e^{\bot}$ intersects $\bdd B$ since otherwise $K$ would 
contain a minimum in $B$.  Let $\aaa$ be the arc in $\theta \cap B$ 
consisting of the end of of $e_{-}$ and the end of $e_{\bot}$ at the 
$Y$-vertex.  Then the parallelism between $\aaa$ and an arc on $\Pt$, 
guaranteed by Morimoto's theorem \cite{Mo}) describes how to pull the 
end of $e_{-}$ down to change the vertex into a $\lll$-vertex.  This 
thins $\theta$.  
\end{proof}

\section{Knots presented as $p$-eyeglasses}

We will need a second way in which $K$ can be viewed as lying on a 
neighborhood of a normal form graph in $S^3$.  Let $\bowtie$ be the 
``eyeglass'' graph, obtained from two circles $e_l$ and $e_r$ by 
attaching an edge $e_b$ running between them.  If $\bowtie$ is 
imbedded in $S^3$ then a regular neighborhood $\eta(\bowtie)$ of 
$\bowtie$ is a genus two handlebody that can be 
described as follows:  Take two solid tori $T_l$ 
and $T_r$ with cores the loops $e_l$ and $e_r$ and meridian disks 
$\ml$ and $\mr$ respectively, and join them together by a $1$-handle with 
core $e_b$ and meridian $\mb$.  

We will sometimes refer to $e_b$ as the bridge between the cycles 
$e_{l}$ and $e_{r}$.  

\begin{defin} \label{defin:eye}
 Given a normal form eyeglass $\bowtie \subset S^3$, a normal form 
 knot $K \subset \bdd\eta(\bowtie)$ is {\em presented as a 
 $p$-eyeglass on $\bowtie$} if

\begin{itemize}

\item $|K \cap \mr| = 1$

\item $|K \cap \mb| = 2$

\item $|K \cap \ml| = p \geq 1$

\item $K$ always intersects $\ml$ with the same orientation.

\end{itemize}

\end{defin}

Less formally, $K$ can be described as the band-sum, via a band 
running once along the bridge $1$-handle, of a longitude of $T_r$ and a 
$(p, q)$ cable of $T_l$.  Note that if $K$ is presented as a 
$p$-eyeglass on $\bowtie$ and $S^3 - \eta(\bowtie)$ is a handlebody, 
then $K$ is tunnel number one:  Since $K$ goes just once through a 
meridian of $e_r$ we can isotope $H = \eta(\bowtie)$ in $S^3$, 
``vacuuming'' up $K$ with $e_r$ until $K$ is simply a longitude of 
$e_r$.  

It's easy to see that any knot presented as a $p$-eyeglass is also 
presented as a $(p-1, 1)$ quasi-cable on the same underlying 
handlebody $H$.  The difference is in how meridians are chosen to 
define the graph that $H$ is a neighborhood of.  The correspondence is 
given by $\mr = \mm$, $\ml = \mz$ while the meridian disks $\mb$ and 
$\mpp$ intersect in a single arc.  That is, the difference of the two 
graphs is a simple Whitney move.  The $\theta$ graph obtained from 
$\bowtie$ this way is called the {\em associated $(p-1,1)$ 
quasi-cable} and the disk in $\bowtie$ just described that becomes the 
meridian $\mpp$ will be called the pre-cable disk in $\bowtie$.

\bigskip

Recall that a graph is in bridge position if every maximum lies above 
every minimum.  Following \cite{ST2}, we will extend this notion in 
the case of an eyeglass graph.  

\begin{figure}
\centering
\includegraphics[width=.6\textwidth]{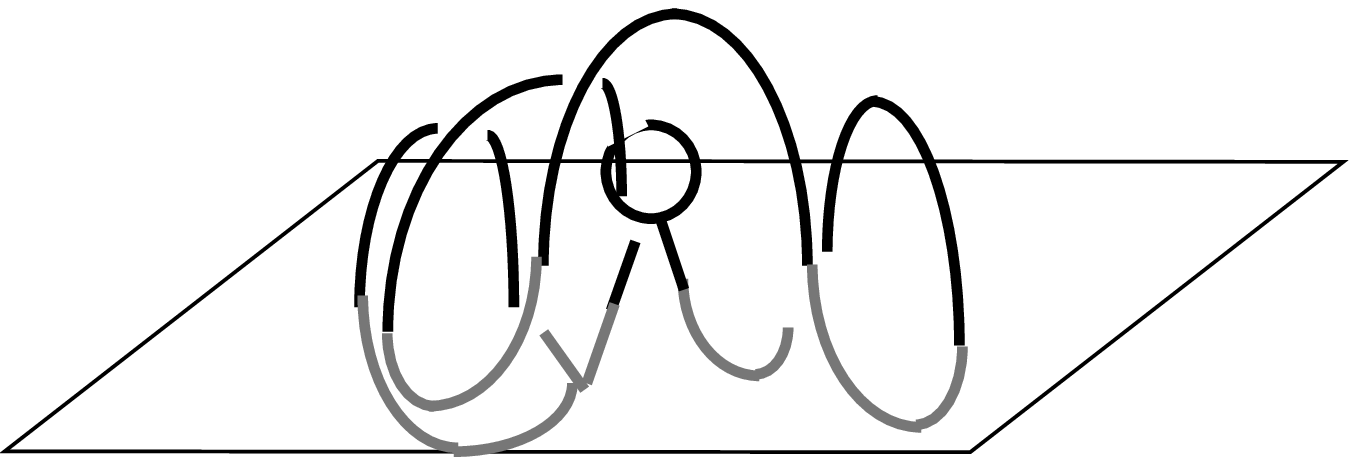}
\caption{} \label{fig:extbridge}
\end{figure}

\begin{defin}  \label{defin:extend}
Suppose a height function is defind on $S^{3}$.  A cycle in $S^{3}$ is 
{\em vertical} if it has exactly one minimum and one maximum.  An 
eyeglass graph is in {\em extended bridge position} if any minimum 
that does lie above a regular maximum (resp.  maximum that lies below 
a regular minimum) is a $Y$-vertex at the minimum (resp.  
$\lll$-vertex at the maximum) of a vertical cycle.  A vertical cycle 
whose minimum is a $Y$-vertex is called an extended maximum.  One 
whose maximum is a $\lll$-vertex is called an extended minimum.  Such 
a $Y$-vertex or $\lll$-vertex is called a {\em base vertex} of the 
extended maximum (resp.  minimum).  A level sphere that lies above all 
minima and extended minima (except perhaps a base $Y$-vertex) and 
below all maxima and extended maxima (except perhaps a base 
$\lll$-vertex) is called a {\em dividing sphere} for the 
extended bridge position.  See Figure \ref{fig:extbridge}
\end{defin}

\begin{prop} \label{prop:bridgeye}  
Suppose $K$ is thinly presented as a $p$-eyeglass on $\bowtie$, whose 
closed complement is also a genus two handlebody.  Suppose that 
$\bowtie$ is made as thin as possible, subject to the condition that 
it thinly presents $K$.  Then either 
\begin{enumerate}

\item $e_{l}$ is vertical, 

\item  $p = 1$  and $e_{r}$ is vertical, or 

\item $\bowtie$ is in extended bridge position.

\end{enumerate}

Moreover, in the last case, $e_{b}$ is disjoint from some dividing sphere 
for $\bowtie$.
\end{prop}

\begin{proof} The proof follows the same line of argument as the proof 
of the main theorem of \cite{ST2}.  We only need to verify that the 
argument there does not interfere with the thin presentation of $K$.  
In fact, the argument here is simpler because some of the more 
complicated steps in \cite{ST2} are required only after a step that, 
in our case, clearly thins $K$ or shows that e.  g.  $e_{l}$ is 
vertical.

So suppose $\bowtie$ thinly presents $K$ but is not in extended bridge 
position.  We've noted above that $K$ is a tunnel number one knot so 
we know that $K$ is in bridge position.  So if a maximum of $\bowtie$ 
lies below a minimum, either the maximum is a $\lll$-vertex or the 
minimum is a $Y$-vertex, or both.  So there are at most two level 
spheres with the property that each lies just below a minimum and just 
above a maximum.  Let $Q$ be the sphere or pair of spheres with this 
property.  Compress $Q$ as much as possible in the complement of 
$\bowtie$ and call the result $Q'$.

A path in $\bdd H = \eta(\bowtie)$ between meridians is {\em regular} 
if the corresponding path in $\bowtie$ is embedded.  It is shown in 
\cite{ST2} that there is a disk $F$ in $S^3$ whose boundary is the 
union of a path $\aaa$ in $\bowtie$ and an arc $\bbb$ in $Q'$.  
Moreover the interior of $F$ is disjoint from $Q'$ and either $\aaa$ 
is a regular path that is disjoint from some meridian of 
$e_{b}$ or $\aaa$ has both its ends at the same point $p$ of $e_{b} 
\cap Q'$ and runs once around either $e_{r}$ or $e_{l}$.  Consider 
each possibility in turn.

{\bf Case 1:} The path $\aaa$ is a regular path that is
disjoint from some meridian of $e_b \subset \bowtie$.

Say $F$ lies below $\bbb \subset Q'$.  By general position we can 
assume that $\bbb$ is disjoint from the disks in $Q'$ which are the 
remains of the compressing disks of $Q$, so in fact $\bbb$ lies on 
$Q$.  If $\aaa$ does not pass through a vertex of $\bowtie$ then just 
use $F$ to isotope the arc of $\bowtie - Q$ that contains $\aaa$ to 
$\bbb \subset Q$.  During the isotopy, as $\aaa$ perhaps passes 
through $Q$ (though not through $Q'$), $\bowtie$ may get thicker, but 
once it reaches $\bbb$ it will have been thinned, since all that 
remains of its internal critical points is one minimum, which will 
have been brought up above the level of the maximum just below $Q$.  
The move similarly cannot thicken $K$.

Essentially the same argument applies even when $\aaa$ passes through 
a vertex.  $F$ is used to slide an end of one of the edges incident to 
the vertex down to $\bbb$ perhaps thereby just extending $e^{b}$ and 
not affecting the bridge structure of $K$.  In any case $\bowtie$ is 
thinned and $K$ is not thickened.  

{\bf Case 2:} $\aaa$ has both its ends at the same point $q$ of $e_{b} 
\cap Q'$ and runs once around either $e_{r}$ or $e_{l}$.

Much as in the previous case, $F$ can be used to move the cycle 
$e_{r}$ or $e_{l}$ together with the end of $e_b$ between $q$ and the 
cycle to $Q$.  Unless the cycle was already vertical, this move (once 
the cycle is tilted again to restore genericity) will thin $\bowtie$ 
and will not thicken $K$.  So we can assume the cycle is vertical.  
This means we are done, unless in fact the cycle is $e_{r}$ and $p 
\geq 2$.  This case only arises if $Q$ intersects $e_{l}$ but not 
$e_{r}$ since if it is disjoint from both, we can appeal to \cite{Mo} 
directly to get a disk as in Case 1.  (Unless, of course, some 
component of $Q$ intersects $\bowtie$ in exactly one point of $e_{b}$.  
But then $\bowtie$ would be planar.)

Now note that, unless the maximum just below $Q$ is a $\lll$-vertex 
maximum of $e_{r}$ the move on $e_{r}$ just described pushes a minimum 
of $K$ (on $e_{r}$) past a maximum of $K$, contradicting the thin 
position of $K$.  We deduce that the one and only maximum just below 
$Q$ is in fact a $\lll$-vertex maximum of $e_{r}$.  If $Q$ consists of 
more than one sphere, we could repeat the same argument, just using 
the component to which we have not just pushed $e_{r}$.  But that 
would lead to the contradiction that the $\lll$-vertex maximum of 
$e_{r}$ also lies just below the other plane.  We deduce that $Q$ is a 
single plane.  We have shown then that, aside from the base 
$\lll$-vertex of $e_{r}$, only minima lie below $Q$.  Moreover, just 
above $Q$ is a minimum and at no other level does a maximum lie below 
a minimum.  It follows that $\bowtie$ is in extended bridge position.

It remains to show that $e_{b}$ is disjoint from some dividing sphere.  
Much of the proof mimics \cite{ST2}.  We suppress most of 
the technical details, except to note that many of the technical 
problems do not arise in our context.  Most importantly, if $P$ is a 
dividing sphere and there are disjoint lower and upper caps, then 
pushing a vertical cycle which is, say, an extended maximum down past 
a minimum would immediately thin $K$, even if (as discussed in 
\cite{ST2}) passing other maxima might thicken $\bowtie$.  We deduce 
that in our context disjoint lower and upper caps cannot arise for 
elementary reasons.

In any case, the upshot of the argument in \cite{ST2} is 
that there is a dividing sphere $P$ that cuts off from a meridian 
disk $E$ of $S^{3} - \eta(\bowtie)$ both an upper disk $D_{u}$ and a 
lower disk $D_{l}$.  Moreover the interior of each is disjoint from 
$P$.

Consider the components $C_u$ and $C_l$ of $\eta(\bowtie) - P$ to which 
$D_u$ and $D_l$ are incident.  If neither $C_u$ nor $C_l$ contain 
vertices, or if they have no ends in common, or if together they 
contain at most one vertex and they have a single end in common, then 
it is easy to use $D_u$ and $D_l$ to push a maximum down past a 
minimum, contradicting the thinness of $\bowtie$ or $K$.  Note in particular 
that if the boundary of one of the disks, say $D_{l}$, goes once 
around $e_{r}$, although the move described may thicken $\bowtie$ (cf 
\cite{ST2}) it does push a minimum of $K$ (namely the minimum of 
$e_{r}$) past a maximum of $K$ and so would violate the assumption 
that $\bowtie$ thinly presents $K$.  

We now proceed to dispose of the other cases.  Suppose that, say, 
$C_u$, contains a vertex and that $C_u$ and $C_l$ have two end 
meridians in common.  We can assume $C_{u}$ contains only one vertex, 
else $e_{b}$ is disjoint from $P$ and we are done.  Then $C_{u} \cup 
C_{l}$ contains a vertical cycle.  If that cycle is $e_{l}$ or $p = 1$ 
we are done, so we'll assume it's $e_r$ and that $p \geq 2$.  Either 
$D_u$ and $D_l$ can be used to make $\bowtie$ (indeed $K$!)  thinner 
or they can be used to isotope $e_r$ into $P$.  This last move not 
only makes $\bowtie$ non-generic, but it may thicken $K$ if $K$ winds 
around $e_r$.  Nonetheless, we persist, inspired by the proof of 
\cite[Theorem 5.14]{GST}.  That argument shows that, once $e_r$ is 
level, so the solid torus neighborhood $T_r$ divides $P$ into two 
disks, an innermost disk component of $E \cap P$ in $E$ or a disk cut 
off by an outermost arc of $E \cap P$ in $E$ can be used to push a 
maximum (resp.  minimum) of $\bowtie$ (possibly the maximum near the 
end of $e_b$ at $e_r$) down (resp.  up) through the level of $P$.  
Afterwards, $e_r$ can be tilted slightly to restore genericity and 
thereby to remove the extra bridges of $K$ that may have been 
introduced when $e_{r}$ was made perfectly level.  Since a 
maximum has been pushed down (or a minimum up) past $e_r$, it 
follows that $\bowtie$ (indeed $K$, since $p \geq 2$) has been 
thinned, the usual contradiction.  

The possibility remains that $C_u$ and $C_l$ each have a single vertex 
and they also have a single end in common (they can't have two ends in 
common since the result would be a cycle in $\bowtie$ containing both 
vertices.)  Their common end must be a point of $e_b \cap P$, since in 
$\bowtie$ that is the only arc that connects the two vertices; in 
particular, $e_b$ is monotonic.  In this case, the disks $D_u$ and 
$D_l$ either could be used to thin $\bowtie$ (an immediate 
contradiction) or they can be used to make $e_b$ level.    If 
neither component $C_{u}$ or $C_{l}$ contains all of $e_{r}$, then the move simply levels 
$e_{b}$.  Once again, this move makes $\bowtie$ no longer generic and 
may also thicken $K$, for $K$ may wind many times around the edge $e_b$.  But 
we continue anyway, inspired this time by the proof of \cite[Theorem 
6.1, Subcase 3b]{GST}.  The argument has a number of subcases, but all 
result in the following conclusion: a maximum (say) of $e_l$ or $e_r$ 
(possibly contiguous to an end of $e_b$) can be pushed down to the 
level of $e_b$ (or below, if it is not contiguous).  Once this is 
achieved, tilt $e_b$ slightly to restore genericity, but leave the 
pushed down maximum at (or below) the lower end of $e_b$.  The result 
is a thinning of $\bowtie$ (indeed $K$).

The final possibility is that the move just described levels all of 
$e^{b}\cup e^{r}$ because, say, $C_{u}$ contains $e_{r}$.  In this 
case it seems that the move might pull $C_{u}$ past other maxima lying 
below it, thickening $\bowtie$ so that, after $e^{b}\cup 
e^{r}$ is tilted to restore genericity, $\bowtie$ actually ends up 
thicker.  Nevertheless, it is argued in \cite{ST2} that in fact this 
does not happen, or at least, if it does, the extra thickness (and 
more) can immediately be removed by a move analogous to that described 
above when $e^{r}$ was levelled.   This is established by a somewhat 
complicated combinatorial argument on $\bdd E$.  We won't repeat the 
argument here. The upshot is that either $\bowtie$ ends up thinner, a 
contradiction, or there were in fact no maxima between $e_{r}$ and 
$P$.  But in this last case, all the minima (including the base 
$Y$-vertex of $e_{r}$) lie below all the maxima, so $\bowtie$ is in 
non-extended bridge position.  (Only the base vertex of $e_{r}$ lies 
between $P$ and the level plane $P'$ for this bridge presentation.)  
Moreover, all of $e_{b}$ lies below $P'$, verifying the proposition in 
this case as well.
\end{proof}

\begin{cor} \label{cor:eyetocable} 
Suppose $K$ is thinly presented as a $p$-eyeglass on $\bowtie$, $S^3 - 
\eta(\bowtie)$ is a handlebody, and $\bowtie$ has been made as thin as 
possible.  Then either

\begin{enumerate}

\item $e_l$ is vertical 

\item $p = 1$ and $e_r$ is vertical 

\item $genus(K) \geq 2$

\item $p \geq 2$ and $K$ is also thinly presented as a $(p-1, 1)$ 
quasi-cable on the associated $\Theta$-graph $\theta$ and $\theta$ is 
no thicker than $\bowtie$.

\item $p = 1$ and the graph $K \cup \ggg$ is no thicker 
than $\bowtie$.

\end{enumerate}
\end{cor}

{\bf Remark:} The last two possibilities are essentially the same and 
are only distinguished by the value of $p$.  It's convenient to 
restrict the terminology $(p, q)$ quasi-cable to the case $p, q \geq 
1$.  If we were to extend that definition to $(1, 0)$ quasi-cable, 
then the original graph $K \cup \ggg$ would be its natural meaning.  
See the beginning of Section \ref{section:thinquasi}.

\begin{proof} Following Proposition \ref{prop:bridgeye} we only need 
to consider the case in which $\bowtie$ is in possibly extended bridge 
position and $e_{b}$ is disjoint from a dividing sphere.  Suppose 
first that $\bowtie$ is in fact in (non-extended) bridge position and, 
with no loss of generality, suppose $e_{b}$ lies above the dividing 
sphere.  Then $e_{b}$ ascends from the lower $\lll$-vertex and either 
is monotonic or it has one internal minimum and descends into the 
other $\lll$-vertex as well.  It's easy to move from one position to 
the other without affecting the width of either $\bowtie$ or $K$, so 
we'll assume for concreteness that $e_b$ is monotonic.  By (perhaps) 
twisting around the other two ends at the lower $\lll$-vertex we can 
ensure that the meridian disk for the associated $\Theta$-graph, 
namely the pre-cable disk in $\bowtie$ that runs the length of 
$e_{b}$, is disjoint from the descending disk incident to $e_{b}$ 
given by the bridge structure.  See Figure \ref{fig:descend2}. Then the 
Whitney move has no effect on the bridge structure (hence the width) 
of $K$, nor the width of $\bowtie$: a pair of $\lll$-vertices with an 
edge between them is replaced by exactly the same thing.

\begin{figure}
\centering
\includegraphics[width=.8\textwidth]{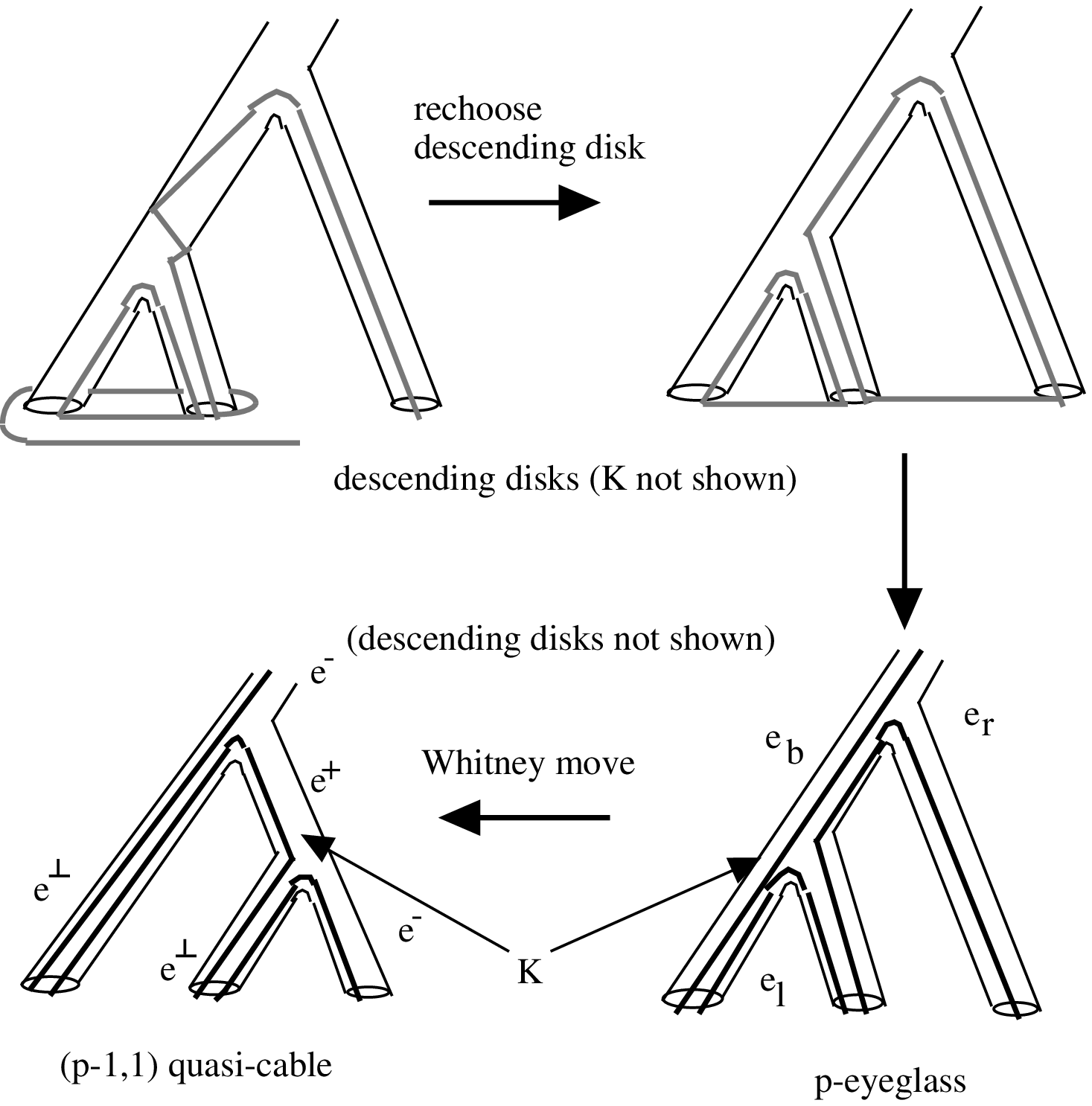}
\caption{} \label{fig:descend2}
\end{figure}	 
\bigskip

Now assume that $\bowtie$ has an extended maximum, say.  Since the 
extended maximum contains a vertical cycle, we are done immediately 
unless the vertical cycle is $e_{r}$ and $p \geq 2$.  Since $e_{b}$ is 
disjoint from the dividing sphere (say it lies above), it runs 
monotonically from a $\lll$-vertex to the $Y$-vertex base of $e_{r}$.  
Since $p \geq 2$ the knot $K$ has a maximum at the $\lll$-vertex.  Now 
find, somewhere below the $\lll$-vertex, a level sphere $P$ as in 
Lemma \ref{lemma:essential}, so every component of $P \cap F$ is 
essential in the Seifert surface $F$.  Then Lemma \ref{lemma:hexagon2} 
applied to the associated $(p-1, 1)$ quasi-cable shows that $genus(K) 
\geq 2$.
\end{proof}

\section{Thinning quasi-cables} \label{section:thinquasi}

Suppose $K$ is a tunnel number one knot with $\ggg$ an unknotting 
tunnel.  According to \cite[Proposition 4.2]{ST1} there is a minimal 
genus Seifert surface $F$ for $K$ that is disjoint from $\ggg$.  We 
now pursue the line of argument used in \cite{GST} and \cite{ST1} to 
analyze the relation between $F$ and the pair $(K, \ggg)$.  The 
philosophy will be to view the pair $(K, \ggg)$ as an incipient case 
of $K$ being thinly presented as a $(p, q)$ quasi-cable, though with 
$(p, q) = (1, 0)$.  Here the graph $\theta$ would be $K \cup \ggg$, 
with $\ggg$ playing the role of $e^{-}$ and the two segments of $K$ 
into which the ends of $\ggg$ divide $K$ playing (interchangeably) the 
roles of $e^{\bot}, e^{+}$.  Roughly, the idea is this: inspired by 
\cite{GST} we will consider the thinnest graph $\theta \subset S^{3}$ 
that thinly presents $K$ as a $(p, q)$ quasi-cable and, inspired by 
\cite{ST1}, ask how $F$ and the height function $h$ interact with 
``splitting'' spheres for the handlebody $\eta(\theta) \subset S^{3}$.

We begin by setting some terminology and notation.  In analogy to the 
notation used for quasi-cables, denote the meridian of $H = \eta(K 
\cup \ggg)$ corresponding to a point of $\ggg$ by $\mu^{-}$ 
and meridians corresponding to points in the two edges of $K - \ggg$ 
by $\mu^{+}$ and $\mu^{\bot}$.  In general, for $H$ a 
handlebody in $S^{3}$ whose closed complement is also a handlebody, 
a {\em splitting sphere} $S$ for $H$ is a sphere that intersects $\bdd H$ 
in a single essential circle.  In other words, it is a reducing sphere 
for the Heegaard splitting $S^3 = H \cup_{\bdd H} (S^3 - int(H))$.  A splitting 
sphere $S$ is best viewed as the union of two disks, $D = S \cap H$ 
and $E = S - int(H)$ that have a common boundary in $\bdd H$.  The exterior 
disk $E$ will, much as above, suggest possible thinning moves for 
$\theta \supset K$.  The interior disk $D$ will give useful 
information about how $\bdd E$ can behave, since $\bdd D = \bdd E$.

\begin{defin} \label{defin:wave}
Let $\{ \mu_i \}$ be a (not necessarily complete) family of pairwise 
disjoint meridian disks for $H$ and $c$ be a simple closed curve on 
$\bdd H$ isotoped so as to minimize $|c \cap (\cup_i \mu_i)|$.  
Suppose an arc component $c_{0}$ of $c - \{ \mu_i \}$ has both its 
ends on a single meridian $\mu_{i}$ in the family and the union of 
$c_{0}$ and a subarc of $\bdd \mu_{i}$ bound a disk in $H$.  Then 
$c_{0}$ is a {\em wave} of $c$ with 
respect to $\{ \mu_i \}$. The wave 
is said to be {\em based} at the meridian $\mu_{i}$.

In particular, if $D$ is an essential disk in $H$, isotoped so as to 
minimize $|D \cap (\cup_i \mu_i)|$, then all components of this
intersection are arcs, and an outermost arc of intersection in $D$ 
cuts off a disk $D_0 \subset D$ so that $\bdd D_0 \cap \bdd H$ is a 
wave of $\bdd D$ with respect to $\{ \mu_i \}$.  The disk $D_0$ is 
called a {\em wave disk}.
\end{defin}

\begin{defin}
Suppose $\theta$ thinly presents the pair $(K, F)$ as a $(p, q)$ 
quasi-cable with $q \leq p$.  Consider the family of meridians $\{ 
\mm, \mpp, \mz \}$ for the handlebody $H = \eta(\theta)$.  Then an 
essential disk $D$ in $H$ {\em satisfies the wave condition} if there 
is a wave of $\bdd D$ based at either $\mm$ or $\mpp$ (so in 
particular the wave is disjoint from $\mz$).  

Similarly, if $S$ is a splitting sphere for $H$, then 
$S$ satisfies the wave condition if $D = S \cap H$ does.
\end{defin}

Note that $D$ satisfies the wave condition if and only if some 
outermost disk $D_{0}$ of $D$ (hence all outermost disks) cut off by 
$\{ \mm, \mpp, \mz \}$ is cut off by an arc lying in either $\mm$ or 
$\mp$.

With this terminology, \cite[Corollary 5.3]{ST1} can be reinterpreted 
as follows (noting that ``$\rho$ is finite'' translates to ``$S$ has a 
wave at $\mm$''):

\begin{lemma}  \label{lemma:onezero}
Suppose $K$ is a tunnel number one knot with $\ggg$ an unknotting 
tunnel, and $F$ is a minimal genus Seifert surface $F$ for $K$ that is 
disjoint from $\ggg$.  As described above, let $\theta = K \cup \ggg$ 
present $(K, F)$ as a $(1,0)$ quasi-cable.  Then $\ggg$ may be slid 
and isotoped so that either 

\begin{itemize}

\item $\ggg$ lies on $F$ or 

\item there is a splitting sphere for $\eta(\theta)$ so that, with 
respect to the set of meridians $\{ \mm, \mpp, \mz \}$, a wave of $D$ 
is based at $\mm$.

\end{itemize}

\end{lemma}

In fact more of \cite{ST1} can be reinterpreted in this setting.  If 
$K$ is not a $2$-bridge knot, an invariant $\rho(K, \ggg) \in 
\mathbb{Q}/2 \mathbb{Z}$ is defined and, if $\rho \neq 1$ it is shown 
that $\ggg$ can be isotoped onto $F$.  The case $\rho = 1$ translates 
to this statement: For the pair of meridians $\{ \mpp, \mz \}$ (which 
are interchangeable in this context) and any splitting sphere $S$ 
there are waves of $D = S \cap H$, with the property that each wave 
disk intersects $\mm$ in a single arc.  In particular, a wave disk 
at $\mpp$, say, can be glued to a subdisk of $\mpp$ to get a 
non-separating meridian $\mu^{wave}$ of $H$ that is disjoint from the 
wave disk and intersects $\mm$ in a single arc.  Then the 
meridians $\{ \mpp, \mz, \mu^{wave} \}$ give $H$ the structure of 
a $\Th$-graph $\theta'$ which presents $(K, F)$ as a $(1, 1)$ 
quasi-cable satisfying the wave condition.  In $\theta'$, the meridian 
$\mpp$ has become the meridian of the edge $e^{\prime -}$, $\mz$ the 
meridian of the edge $e^{\prime +}$, and $\mu^{wave}$ the meridian of 
the edge $e^{\prime \bot}$.  The two graphs $\theta$ and $\theta'$ 
differ by a standard Whitney move on $\ggg$.  So we have:

\begin{lemma} \label{lemma:oneone}
Suppose $K$ is a tunnel number one knot with $\ggg$ an unknotting 
tunnel, and $F$ is a minimal genus Seifert surface $F$ for $K$ that is 
disjoint from $\ggg$.  Then either there is a graph $\theta$ that 
presents $(K, F)$ as a $(1,1)$ quasi-cable and a splitting sphere that 
satisfies the wave condition, or $\ggg$ may be isotoped to lie on $F$.
\end{lemma}

The next ingredient to throw into the mix is thin position.  Of course 
in Lemma \ref{lemma:onezero} there is no obstacle to having $\theta$ 
thinly present $K$ - just begin with $K$ in thin position.  It's not 
obvious that the process that leads from Lemma \ref{lemma:onezero} to 
Lemma \ref{lemma:oneone} preserves the property that $K$ is thinly 
presented, but we now show that it does, (or $\ggg$ can be made into  
an unknotted loop).  Or, more accurately, we show that this follows 
immediately from the results of \cite{GST} and \cite{ST1}.  

\begin{lemma}  \label{lemma:oneonethin}
Suppose $K$ is a tunnel number one knot with $\ggg$ an unknotting 
tunnel, and $F$ is a minimal genus Seifert surface $F$ for $K$ that is 
disjoint from $\ggg$.  Then either

\begin{enumerate}

\item there is a graph $\theta$, no thicker than the graph $K \cup 
\ggg$, that {\em thinly} presents $(K, F)$ as a $(1,1)$ quasi-cable 
and a splitting sphere that satisfies the wave condition, or

\item $\ggg$ may be isotoped to lie on $F$, or 

\item $\ggg$ can be slid and isotoped to form an unknotted loop with 
its ends at the same point of $K$.

\end{enumerate}
\end{lemma}

\begin{proof} Following \cite[Theorem 3.5]{GST}, $\ggg$ can be slid to 
become either an unknotted loop, and we are done, or $\ggg$ is a level 
edge, with its ends incident to the top two maxima (say) of the thinly 
presented $K$.  (In the latter case, regain a generic positioning by 
slightly perturbing $\ggg$ from its level position, changing it to a 
monotone edge connecting two $\lll$-vertices.)  If $K$ is $2$-bridge, 
it's easy to see that $F$ contains an isotopic copy of $\ggg$ (cf 
\cite[Remark 12.26]{BZ}), and we are done.  If $K$ 
is not $2$-bridge and the invariant $\rrr(K, \ggg) \in \mathbb{Q}/2 
\mathbb{Z}$ then defined in \cite{ST1} is not $1$, it follows from 
\cite[Theorem 5.2]{ST1} that $\ggg$ may be isotoped into $F$ and we are 
done. 

Consider finally the case $\rrr(K, \ggg, S) = 1$ for $S$ a splitting 
sphere as in Lemma \ref{lemma:onezero}.  The fact that $\rrr = 1$ 
means that, the pair of meridians $\mz, \mpp$ cuts off a wave of $D = 
S \cap H$ based at the meridian $\mpp$, say, and that wave intersects 
$\mm$ in a single arc.  In particular (as above) the wave disk can be 
glued to a subdisk of $\mpp$ to get a non-separating meridian 
$\mu^{wave}$ of $H$ that is disjoint from the wave disk and intersects 
$\mm$ in a single arc.  Appropriately twist  
the two arcs of $K$ that descend from the bottom $\lll$-vertex of 
$\ggg$ (equivalently, choose an appropriate set of descending disks) 
as discussed in Corollary \ref{cor:eyetocable} (see Figure 
\ref{fig:descend2}) so that the descending disk is disjoint from 
$\mu^{wave}$.  Then the standard Whitney move on $\ggg$, using 
$\mu^{wave}$ not only converts $\theta$ to a graph $\theta'$ that 
thinly presents $K$ as a $(1,1)$ quasi-cable, it does it without 
thickening $\theta$, for one pair of $\lll$-vertices is just replaced 
with another.  (See again Figure \ref{fig:descend2}). \end{proof}

We consolidate our results a bit more:

\begin{defin} Suppose $K$ is a knot, $F$ is a minimal genus Seifert 
surface for $K$ and $\theta$ is a $\Theta$-graph such that $\theta$ 
thinly presents $(K, F)$ as a $(p, q)$ quasi-cable, $p \geq q \geq 1$ 
and such that there a splitting sphere for $\eta(\theta)$ that 
satisfies the wave condition.  We say that $\theta$ is an {\em 
appropriate} $\Th$-graph (for the pair $(K, F)$ and splitting sphere 
$S$.)
\end{defin}

\begin{prop}  \label{prop:oneonethin}
Suppose $K$ is a tunnel number one knot and $\ggg$ is an unknotting 
tunnel for $K$.  Then either 

\begin{enumerate}

\item $\ggg$ can be slid and isotoped to form 
an unknotted loop with its ends at the same point of $K$ or 

\item there is a 
minimal genus Seifert surface $F$ for $K$ with this property: either 

\begin{itemize}

\item $\ggg$ can be slid and isotoped to lie on $F$, or 

\item there is an appropriate $\Theta$-graph for $(K, F)$ (and some 
splitting sphere).

\end{itemize}
\end{enumerate} 

In the last case, if the $\Theta$-graph is put in thin position then 
it is in bridge position.

\end{prop}

\begin{proof} As noted above, according to \cite[Proposition 4.2]{ST1} 
there is a minimal genus Seifert surface $F$ for $K$ that is disjoint
from $\ggg$.  The rest follows from Lemmas
\ref{lemma:oneonethin} and \ref{lemma:bridgepq}.
\end{proof}

We will expand on the last possibility, but it will be useful to have 
the following general lemma:

\begin{lemma} \label{lemma:bridgehelp}
Let $\theta \subset S^{3}$ be a $\Theta$-graph with edges $e_{1}, 
e_{2}, e_{3}$.  Suppose $\theta$ is in bridge position, $P$ is a 
dividing sphere, and the cycle $e_{1} \cup e_{2}$ is vertical.  Then 
either

\begin{enumerate}

\item all cycles in $\theta$ are vertical

\item there are bridges above and below $P$ made up of interior 
subarcs of $e_{3}$ or

\item one of $e_{1}, e_{2}$, say $e_{1}$, is disjoint from $P$ and 
$e_{1} \cup e_{3}$ is vertical.

\end{enumerate}
\end{lemma}

\begin{proof} Suppose first that the vertices lie on the same side of 
a dividing sphere $P$, so, with no loss of generality, both are 
$Y$-vertices, say.  Then each edge intersects $P$ in an even number of 
points.  Since $e_{1} \cup e_{2}$ is vertical, together the edges 
intersect $P$ in $2$ points.  Hence one of these two edges, say 
$e_{1}$, is disjoint from $P$.  This means that $e_{1}$ is the edge 
that descends from the higher $Y$-vertex.  $e_{3}$ can't also descend 
from this vertex, so $e_{3}$ intersects $P$.  If $e_{3}$ intersects 
$P$ in two points then $e_{1} \cup e_{3}$ is vertical.  If $e_{3}$ 
intersects $P$ in four or more points then there are bridges above and 
below $P$ made up of interior subarcs of $e_{3}$, as required.  

Suppose next that the vertices lie on opposite sides of $P$, so one 
vertex is a $\lll$-vertex and the other is a $Y$-vertex, and each edge 
intersects $P$ in an odd number of points.  Since $e_{1} \cup e_{2}$ 
is vertical, each of these edges intersects $P$ in a single point.  If 
$e_{3}$ also intersects $P$ in a single point then every cycle is 
vertical.  If it intersects $P$ in three or more points, then there 
are bridges above and below $P$ made up of interior subarcs of 
$e_{3}$, as required.
\end{proof}

\begin{prop}  \label{prop:oneonethin2}

Suppose $K$ is a knot with Seifert surface $F$ and there is an 
appropriate $\Theta$-graph for $(K, F)$ and some splitting sphere.  
Let $\theta$ be a thinnest such $\Theta$-graph.  Suppose that 
$genus(F) = 1$.  Then $\theta$ is in bridge position and (thinly) 
presents $K$ as a $(p, 1)$ quasi-cable, for some $p \geq 1$.

If furthermore one of the cycles $e^- \cup e^+ $ or $e^- \cup 
e^{\bot}$ is vertical and a dividing sphere intersects both edges of 
the vertical cycle then 
either
\begin{itemize}
\item $e^+ \cup e^{\bot}$ is vertical or
\item $p = 1$ and the cycle $e^- \cup e^{\bot}$ is vertical.
\end{itemize}
\end{prop}

\begin{proof} Corollary \ref{cor:thinpq} notes that $q = 1$. 
Proposition \ref{prop:oneonethin} shows that $\theta$ is in bridge position.

Suppose first that the cycle $e^- \cup e^+$ is vertical.  We apply 
Lemma \ref{lemma:bridgehelp}, using $e^-, e^{+}$ for $e_{1}, e_{2}$.  
If all cycles are vertical then of course we are done.  By hypothesis, 
neither $e^-$ nor $e^{+}$ is disjoint from a dividing sphere.  So we 
may assume, following Lemma \ref{lemma:bridgehelp}, that there are 
bridges above and below a dividing sphere made up entirely of interior 
subarcs of $e_{b}$.  We can of course arrange that the lowest maximum 
and highest minimum are these bridges.  Now choose a level sphere $P$ 
as in Lemma \ref{lemma:essential}.  Since $e^- \cup e^+$ is vertical 
and $P$ intersects both edges, $P$ intersects each edge $e^-, e^+$ in 
a single point.  The result then follows from Lemma 
\ref{lemma:hexagon3}.

Similarly, suppose the cycle $e^- \cup e^{\bot}$ is vertical.  Again 
apply Lemma \ref{lemma:bridgehelp} this time using $e^-, e^{\bot}$ for 
$e_{1}, e_{2}$.  If all cycles are vertical we are done.  By 
hypothesis, neither $e^-$ nor $e^{\bot}$ is disjoint from a dividing 
sphere so we may arrange that the lowest maximum and highest minimum 
are from bridges that lie entirely in $e^{+}$.  Now choose a level 
sphere $P$ as in Lemma \ref{lemma:essential}.  Since $e^- \cup 
e^{\bot}$ is vertical and $P$ intersects both edges, $P$ intersects 
each edge $e^-, e^{\bot}$ in a single point.  The result then follows 
from Lemma \ref{lemma:hexagon4}.
\end{proof}

\bigskip

\begin{defin} Suppose $\theta$ is an appropriate $\Th$-graph for the 
pair $(K, F)$ and splitting sphere $S$.  Suppose there is a level 
sphere $P$ at a generic height for $\theta$ such that $P$ cuts off 
both an upper and a lower disk from $E = S - int(H)$.  (As usual, $P$ 
may lie at a critical height of $interior(E)$).  Then $P$ is called a 
{\em critical sphere} for $\theta$.
\end{defin} 

\begin{lemma}  \label{lemma:levelsphere}
Suppose $\theta$ is a thinnest $\Th$-graph appropriate for $(K, F)$.  
Then either
\begin{enumerate}
\item Some dividing sphere for $\theta$ is also a critical sphere 

\item all the cycles in $\theta$ are vertical or

\item the edge $e^w$ (one of $e^{\pm}$) on whose meridian a wave is 
based is disjoint from a dividing sphere and both of the cycles 
containing $e^w$ are vertical.
\end{enumerate}
\end{lemma}

\begin{figure}
\centering
\includegraphics[width=.8\textwidth]{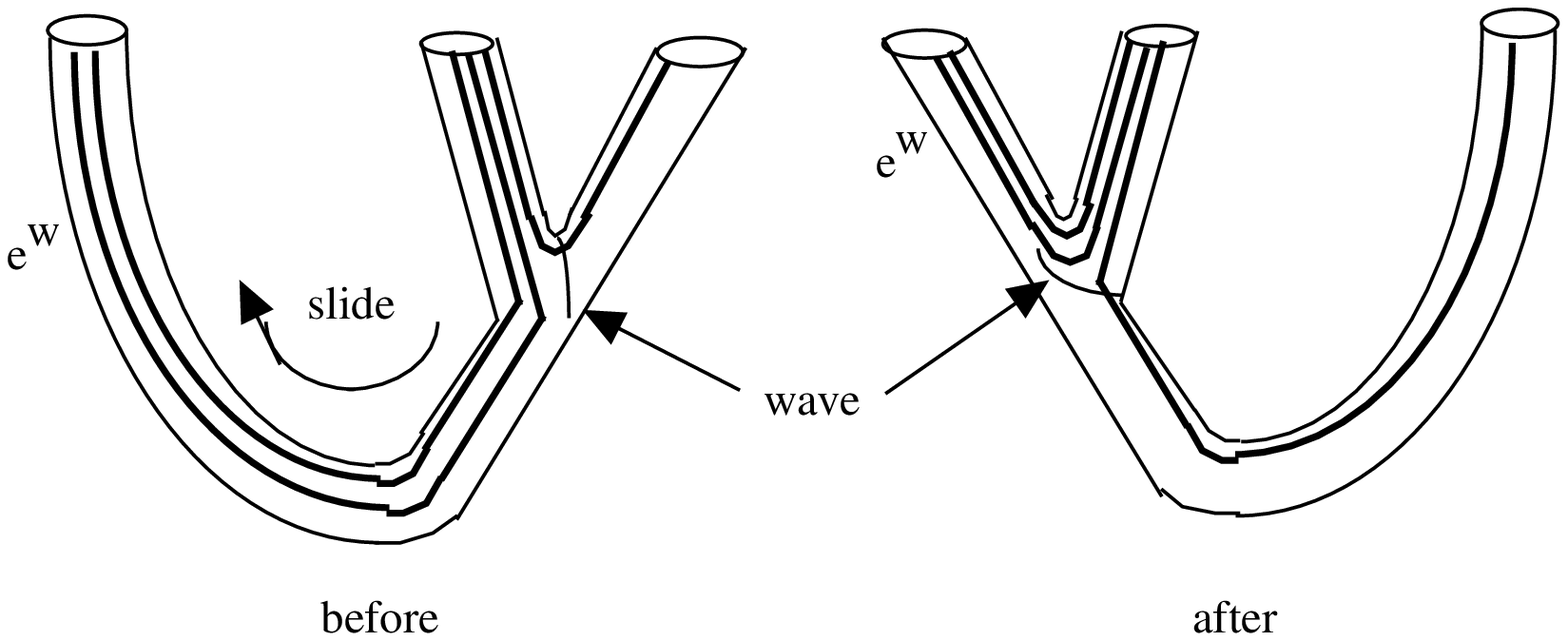}
\caption{} \label{fig:levelsphere}
\end{figure}

\begin{proof} The conclusions make sense, since we know from 
Proposition \ref{prop:oneonethin2} that $\theta$ must be in bridge 
position.  If a sphere just above the highest minimum cuts off a lower 
disk, and a sphere just below the lowest maximum cuts off an upper 
disk, then some level sphere between them is a critical sphere.  This 
condition is guaranteed unless the highest minimum (or the lowest 
maximum) is a $Y$-vertex (resp.  $\lll$ vertex) with, via the wave 
condition, an end of $e^w$ descending (resp.  ascending) from the 
vertex.

So suppose the highest minimum, say, is a $Y$-vertex with an end of 
$e^w$ descending.  Let $P$ be any dividing sphere for $\theta$.  If 
any component of $\theta - P$ below $P$ is a simple arc, its minimum 
could be pushed higher than the $Y$-vertex, eliminating the problem, 
so we can assume that all components of $\theta - P$ below $P$ contain 
vertices.  If there is only one such component and it contains a 
single vertex, then each of the cycles has at most one minimum and so 
each is vertical.  If there is only one such component and it contains 
both vertices, then the edge between them must be incident to the 
higher vertex from below, hence that edge is $e^w$.  Moreover, both 
cycles containing $e^{w}$ have exactly one minimum, and so both are 
vertical.

The remaining case is when the vertices are in separate components, 
each lying below $P$, i.  e.  both of them $Y$-vertices.  Let $v$ 
denote the higher $Y$-vertex.  Push the regular minimum on the 
component containing $v$ up to a height just below $v$.  Now slide the 
end of $e^{\bot}$ ascending from $v$ down to the regular minimum and 
back up the other side.  This has no effect on the width of $K$ 
(since, for example, it doesn't change the number of bridges) but it 
alters the arrangement of the edges around $v$.  In particular, 
afterwards the end of $e^w$ ascends from $v$ so, by the wave 
condition, we can assume that just above $v$, a level sphere cuts off 
a lower disk from $E$, and so somewhere between that level and that of 
the lowest maximum there is a critical sphere as required.  See 
Figure \ref{fig:levelsphere}.
\end{proof}

It will be useful to assume that we always perform the move described 
at the end of the proof of Lemma \ref{lemma:levelsphere} when it is 
possible.  That is, if 
\begin{itemize}

\item $\theta$ is in bridge position

\item a dividing sphere $P$ cuts off a component of $\theta$ 
containing a single vertex (say a $Y$-vertex) and

\item the end of $e^{w}$ descends from that  $Y$-vertex

\end{itemize} then we slide the end of $e^{\bot}$ at that vertex 
over the contiguous regular minimum of $e^{w}$.  This ensures that $e^{w}$ 
ascends from every such $Y$-vertex (and, symmetrically, descends from 
any such $\lll$-vertex) and has no effect on the width of anything, 
since we are just rearranging the heights of minima (or the heights of 
maxima).

\begin{lemma} \label{lemma:edgeabove}
Suppose $\theta$ is a thinnest $\Th$-graph appropriate for $(K, F)$ 
and $\theta$ has a dividing sphere that is a critical sphere.  Then 
either
\begin{enumerate}
\item one of the edges of $\theta$ is disjoint from the critical 
sphere (hence from a dividing sphere) or

\item the edge $e^w$ of $\theta$ on whose meridian a wave is based is 
monotonic and one of the two circuits in $\theta$ containing $e^{w}$ 
is vertical.
\end{enumerate}
\end{lemma}

\begin{proof} Let $P$ be the critical sphere.  Let $C_l$ and $C_u$ be 
the components of $\theta - P$ to which the lower and upper disks 
$D_l$ and $D_u$ are incident.  Suppose first that one of these 
components, say $C_l$, has no vertex (so $C_l$ is a regular minimum).  
If $C_u$ also has no vertex, then $K$ could be thinned, a 
contradiction.  If $C_u$ has two vertices, then it contains an edge 
disjoint from $P$ and we are done.  If $C_u$ has one vertex then, 
since every circuit in $\theta$ has two vertices, $C_l$ and $C_u$ have 
at most one common end point on $P$.  Moreover, by the wave condition, 
the path $\bbb_{u} = \bdd D_{u} \cap P$ has at least one end incident 
to an end of $e^{w}$ in $C_{u}$.  Then $D_u$ and $D_l$ describe how to 
slide the end of $e^w$ in $C_u$ down to $P$ while simultaneously 
isotoping all of $C_l$ up to $P$.  If the slide of the end of $e^w$ is 
down an end of $e^{\bot}$ that would thin $K$; otherwise the effect of 
the slide is to sew together the ends of $e^{\pm}$ extending the end 
of $e^{\bot}$ down to or below the level of the minimum at $D_l$.  
(See Figure \ref{fig:sew}.)
This thins $\theta$ (though not necessarily $K$), possibly by changing 
a $\lll$-vertex into a $Y$-vertex.  This contradicts the 
assumption that $\theta$ is a thinnest such graph.

\begin{figure}
\centering
\includegraphics[width=.8\textwidth]{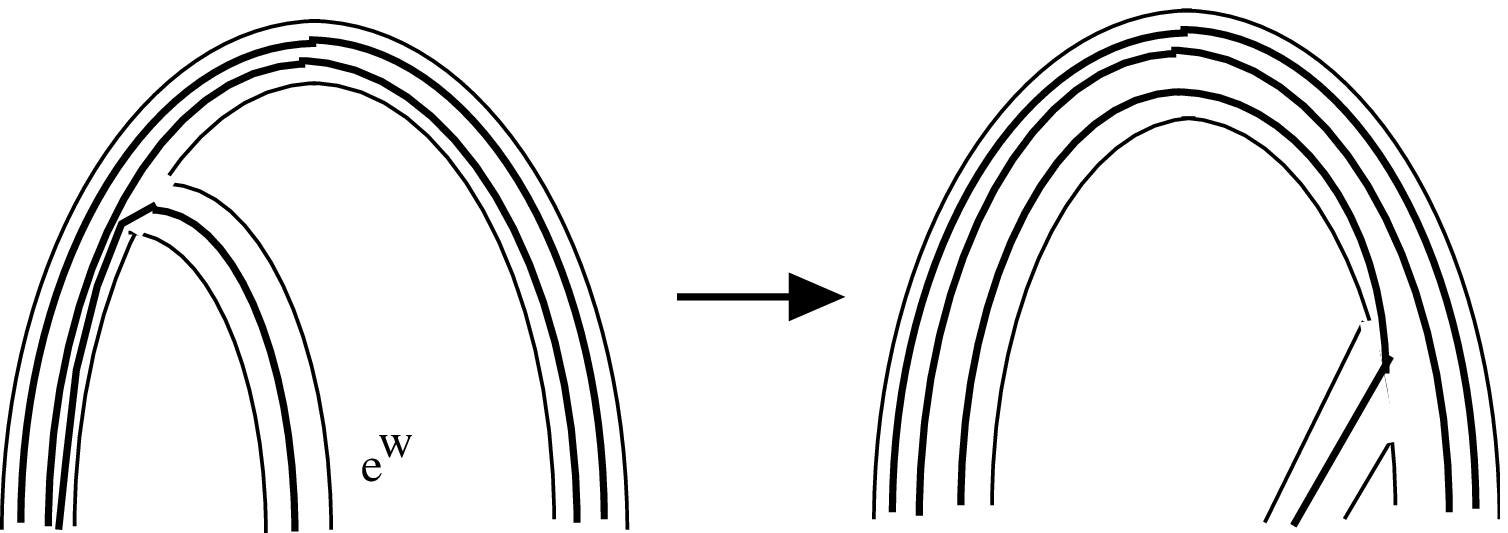}
\caption{} \label{fig:sew}
\end{figure}

The remaining case is that $C_l$ and $C_u$ both contain a single 
vertex.  As above, $D_l$ and $D_u$ can be used to slide the ends of 
$e^w$ in $C_l$ and $C_u$ up and, respectively, down, until they lie in 
$P$.  This would thin $\theta$ unless $e^w$ crossed $P$ in exactly one 
point (an end of both $C_u$ and $C_l$).  In the latter case, $e^{w}$ 
was monotonic and the slide would level $e^w$ by isotoping it into 
$P$.  If the levelled $e^w$ becomes a loop, then originally it was 
part of a vertical cycle, and we are done.  This is also true if one 
of the other edges of $\theta$ has a single interior critical point.  
If neither of these cases occur, then the argument of 
\cite[Theorem 
6.1, Subcase 3b]{GST} can be used (as it was in Proposition 
\ref{prop:bridgeye}) to move the levelled 
$e^w$ up (or down) to connect two maxima or two minima.  Afterwards 
$K$ is no wider but $\theta$ is thinned, a contradiction completing 
the proof.
\end{proof}

\begin{cor} \label{cor:edgeabove} Suppose that $\theta$ is a thinnest 
$\Th$-graph appropriate for $(K, F)$.  Then either

\begin{enumerate}

\item some dividing sphere is a critical sphere and the critical 
sphere is disjoint from one of the edges of $\theta$ or

\item $e^+ \cup e^{\bot}$ is vertical or

\item $p = 1$ and the cycle $e^- \cup e^{\bot}$ is vertical

\item $genus(F) \geq 2$.

\end{enumerate}
\end{cor}

\begin{proof}  Suppose first that some dividing sphere is a critical sphere. 
Then following Lemma \ref{lemma:edgeabove}, if the critical sphere 
intersects all the edges, then some cycle involving one of $e^{\pm}$ is 
vertical.  If it's the cycle $e^+ \cup e^{\bot}$ we're done.  If it's 
either of the other two cycles, apply Proposition 
\ref{prop:oneonethin2}, observing that either some edge is above a 
dividing sphere or the dividing sphere will intersect both edges of 
any vertical cycle. 

If no dividing sphere is a critical sphere then apply Lemma 
\ref{lemma:levelsphere}.  If all the cycles in $\theta$ are vertical 
then of course we are done.  If the wave is based on $e^{+}$ so $e^{w} 
= e^{+}$ then Lemma \ref{lemma:levelsphere} says $e^+ \cup e^{\bot}$ 
is vertical, as required.  Finally, suppose $e^{w} = e^{-}$.  If $p = 
1$ then we are done, since $e^- \cup e^{\bot}$ is vertical.  So 
suppose $p = 2$ and, following Lemma \ref{lemma:levelsphere}, $e^{-}$ 
is disjoint from a dividing sphere.  Find a dividing sphere $P$ as in 
Lemma \ref{lemma:essential} and apply Lemma \ref{lemma:hexagon1}.  
\end{proof}

\section{Regular annuli in handlebody complements}

\begin{defin} Let $\Ggg$ be a trivalent graph in $S^3$, in normal form 
with respect to a height function $h$.  Then a normal form simple 
closed curve $c$ on $\bdd\eta(\Ggg)$ is {\em regular} if $c$ never 
``back-tracks'' along an edge, traversing the edge twice in opposite 
directions, after looping around a vertex at the end of the edge.  
\end{defin}

More precisely, if $\{ \mu_i \}$ is a collection of meridian disks 
in $\eta(\Ggg)$, one for each $1$-handle corresponding to an edge of 
$\Ggg$, then no subsegment of $c$ with interior disjoint from these 
meridians, has its ends incident to the same side of the same 
meridian.  In particular, no minimum of $c$ occurs near a 
$\lll$-vertex of $\Ggg$ and no maximum of $c$ occurs near a 
$Y$-vertex of $\Ggg$.

\begin{lemma}  \label{lemma:unbalanced}
Suppose $\Ggg$ is a trivalent graph in $S^3$ in bridge position with 
respect to a height function $h$ and $A$ is a properly imbedded 
annulus in $S^3 - \eta(\Ggg)$ in normal form with respect to $h$.  
Suppose some dividing sphere $P_0 = P(t_0)$ intersects $\bdd_+ A$ more 
often than it intersects $\bdd_- A$.  Suppose finally that the lowest 
maximum and the highest minimum of $\Ggg$ are regular critical points, 
on edges that are incident to $\bdd_{+} A$.  Then some dividing sphere 
cuts off from $A$ both an upper disk $D_u$ and a lower disk $D_l$, 
both of which are incident to $\bdd_+ A$ (as opposed to $\bdd_{-} A$).
\end{lemma}

\begin{proof} Since the lowest maximum is a regular maximum and $\bdd_{+} 
A$ runs along the edge that contains it, a dividing sphere $P(y)$ just 
below the lowest maximum cuts off an upper disk that is incident to 
$\bdd_{+} A$.  Similarly, a dividing sphere $P(x)$ just above the 
highest minimum cuts off a lower disk that is incident to $\bdd_{+} 
A$.  Since some (hence every) dividing sphere intersects $\bdd_{+} A$ 
more often than it intersects $\bdd_{-} A$, every dividing sphere cuts 
off some disk, either upper or lower, incident to $\bdd_+ A$.  Since 
at $x$ there is a lower one, and at $y > x$ there is an 
upper one, and at every height between there is one or the other, it 
follows that at some height there is both an upper and a lower disk 
incident to $\bdd_{+} A$.  (As usual, this height may be at a saddle 
tangency of $P$ with an interior point of $A$.)  \end{proof}

In the next section we will see that such a useful annulus $A$ can 
often be found.  In this section we examine how, by exploiting upper 
disks and lower disks lying in $A$, we can thin a $\Theta$ curve that 
presents $K$ as a $(p, 1)$ quasi-cable.  Until we begin to use Lemma 
\ref{lemma:unbalanced} there is nothing special about using $A$; any 
properly embedded surface in the graph complement would do, though it 
is important that the upper and lower disks themselves are incident to 
$\theta$ only along regular curves disjoint from $K$.  For example, 
the external disk $E$ used above cannot generally be used for the 
purposes of this section because its boundary is typically not a 
regular curve.

\begin{lemma}  \label{lemma:reguplow+}
Suppose there is a $\Th$-graph appropriate for $(K, F)$ and, for a 
thinnest one $\theta$, there is a dividing sphere $P$ that cuts off 
from $A$ disjoint upper and lower disks $D_u$ and $D_l$ so that the 
arcs $\aaa_u = D_u \cap \bdd\eta(\theta)$ and $\aaa_l = D_l \cap 
\bdd\eta(\theta)$ are regular curves disjoint from $K$.

If the edge $e^+$ is disjoint from $P$ then either

\begin{enumerate}
\item $\aaa_u$ runs from a point of $e^- \cap P$ to a point of $e^{\bot}
\cap P$, traversing $e^+$ once,  

\item $genus(K) \geq 2$ or

\item $e^+ \cup e^{\bot}$ is unknotted 
%(or, if $p = 1, e^- \cup e^{\bot}$ is unknotted.)
\end{enumerate}

\end{lemma}

\begin{proof} We may as well assume $e^+$ lies above $P$, so the
vertices are $\lll$ vertices.  By thinness we can assume that $\aaa_u$
lies on the $4$-punctured sphere component $\Ss$ of $\bdd\eta(\theta)
- P$ lying above $P$. We can think of the components $K
\cap \Ss$ as edges of a graph on $\Ss$, with vertices the four meridian 
boundary components.  That is, if we
label the meridian components of $\bdd \Ss$ as $\mlm$, $\mrm$,
$\mz_l$, $\mz_r$ in the obvious way, then the components of $K \cap
\Ss$ consist of a single arc connecting $\mlm$ to $\mz_l$, a single
arc connecting $\mrm$ to $\mz_r$ and $p$ arcs connecting $\mz_l$ to
$\mz_r$.  The $p$ arcs are either all parallel or comprise two
families of parallel arcs, separating in $\Ss$ the points $\mlm,
\mrm$.  See Figure \ref{fig:reguplow}.

\begin{figure}
\centering
\includegraphics[width=.6\textwidth]{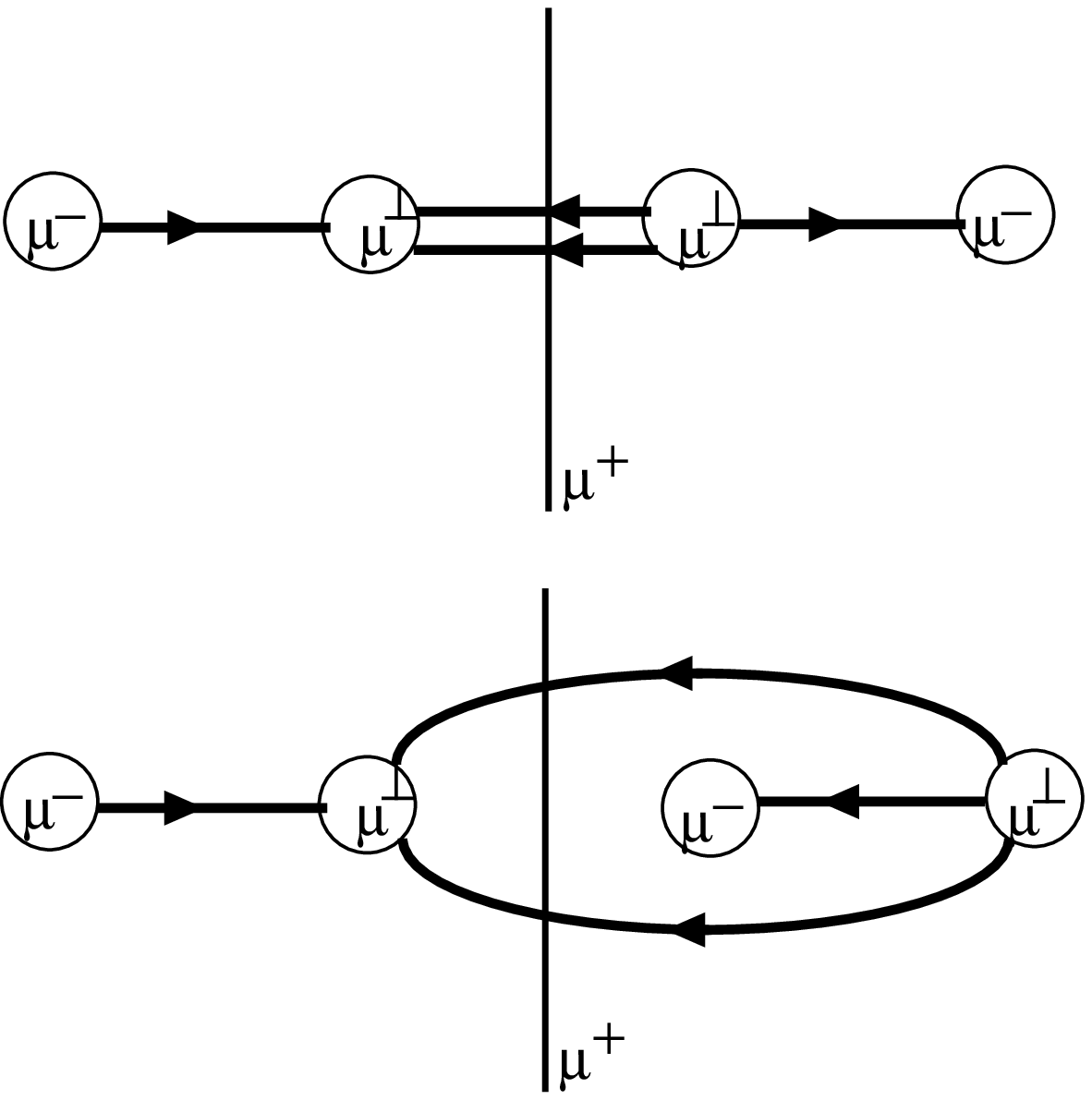}
\caption{} \label{fig:reguplow}
\end{figure} 

Consider the possibilities for $\aaa_{u}$: There is only one path in 
$\Ss$ that is disjoint from $K$ and has ends at meridians $\mlm$ and 
$\mz_l$.  And it's parallel in $\Ss$ to an arc of $K$.  Hence if 
$\aaa_u$ were that path, $K$ could be thinned, using $D_{u}$ and 
$D_{l}$.  Similarly for paths from $\mrm$ to $\mz_r$.  There are (at 
most) two paths from $\mlm$ to $\mz_r$ disjoint from $K$, each 
traversing $e^+$ once, and similarly two paths from $\mrm$ to $\mz_l$.  
If $\aaa_u$ is any of these paths, then the first conclusion of the 
lemma follows.

Now suppose $\aaa_{u}$ is the path in $\Ss$ (available only if all $p$ 
arcs are parallel) that is disjoint from $K$ and runs from $\mlm$ to 
$\mrm$.  Suppose, to begin with, that $\aaa_{l}$ does not also run 
from $\mlm$ to  $\mrm$ but rather has at least one end at another 
point of $\theta \cap P$.  Consider the $p+1$ eyeglass graph 
$\bowtie$ obtained from $\theta$ by a Whitney move along $e^{+}$, 
using as the new meridian a neighborhood of $\mlm \cup \aaa_{u} \cup 
\mrm$ (see Figure \ref{fig:theta2eye}).  Then $D_{u}$ describes a 
descent of the new bridge edge $e_{b}$ for $\bowtie$ down to $P$.  
Simultaneously, $D_{l}$ describes how to move a minimum of (now) 
$\bowtie$ above or at least to the level of $P$.  In particular, 
$\bowtie$ is in bridge position, and a dividing sphere necessarily 
intersects $e_{b}$.  Moreover, $\bowtie$ is thinner than $\theta$ 
since, in effect, a maximum has been pushed below a minimum.  Now 
put $\bowtie$ in thin position.  According to Corollary 
\ref{cor:eyetocable} (exploiting the fact here that $\bowtie$ is a $p+1 > 
1$ eyeglass) either $e_{l}$ becomes vertical (which implies 
that $e^{+} \cup e^{\bot}$ was unknotted) or $genus(K) \geq 2$ or the 
associated $\Theta$-graph, namely $\theta$ can be made as thin as 
$\bowtie$.  The last contradicts the hypothesis and the first two 
possibilities are the conclusions we seek.

\begin{figure}
\centering
\includegraphics[width=.8\textwidth]{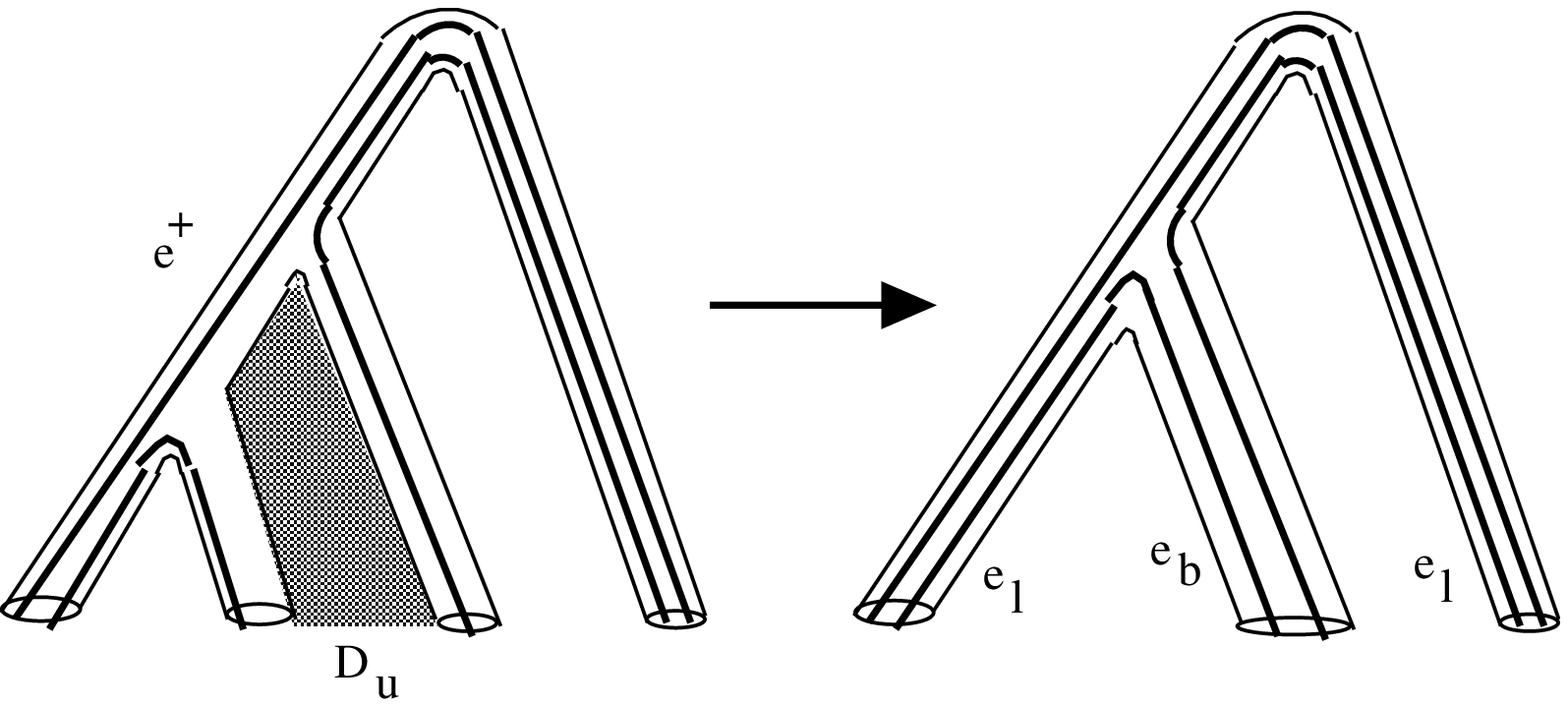}
\caption{} \label{fig:theta2eye}
\end{figure} 

Now suppose $\aaa_{l}$, like $\aaa_{u}$, runs from 
$\mlm$ to $\mrm$ (so, in particular, $e^{+} \cup e^{-}$ is vertical).  
We could construct $\bowtie$ as just described; the construction 
places $e_{r} \subset \bowtie$ level in $P$.  Unfortunately, when 
$e_{r}$ is tilted to restore genericity, $P$ is again a level sphere 
for $\bowtie$ and no thinning will have occurred.  So a different 
argument will be used, this one exploting the fact that $\aaa_{u}$ and 
$\aaa_{l}$ are disjoint from $K$.  Consider the situation once $D_{u}$ 
and $D_{l}$ have been used to put $e_{r}$ into the plane $P$.  Let 
$T_{r}$ be a thickened regular neighborhood of $e_{r}$ and consider 
the two longitudinal circles $\{ \lll_{1}, \lll_{2} \} = P \cap \bdd 
T_{r}$.  It is easy to see what they are: Before the edge $e_{r}$ is 
levelled it intersects $P$ in two points (actually the points $e^{-} 
\cap P$).  The $\lll_{i}$ are obtained by banding the corresponding 
pair of meridians to itself using both $\aaa_{u}$ and $\aaa_{l}$.  The 
important point for our purposes is that $|K \cap (\lll_{1} \cup 
\lll_{2})| = 2$.  Consider where these points lie.  If they both lie 
on the same longitude, then an arc of $K$ they cut off is inessential 
in the annulus component of $T_{r} - (\lll_{1} \cup 
\lll_{2})$ in which it lies, and so it can be removed by an isotopy.  
On the other hand, if one point lies on each longitude, consider the 
algebraic intersection of $K$ with the disk component $P_{1}$ of $P - T_{r}$ 
bounded by $\lll_{1}$, say.  One point is the point $\lll_{1} \cap 
K$.  All others come from intersections of $e_{l}$, i. e. 
intersections of the old $e^{\bot}$.  Each point $e_{l} \cap P_{1}$
contributes $p+1$ points to $K \cap P_{1}$ and they are all of the 
same sign.  So the total algebraic intersection of $K$ with $P_{1}$ 
is $\pm 1 \not\equiv 0 \, mod \, (p+1)$.  This contradicts the fact that 
$K$ bounds $F$ in $S^{3}- H$.  We are left with the conclusion that 
indeed $K$ can be isotoped off of the two longitudes, so $K$ only 
intersects the top of $\bdd T_{r}$.  But in that 
case, consider $F \cap (P_{1} \cup P_{2})$.  It's easy to see that any 
component of intersection that is inessential 
in $F$ can be removed (else $K$ could be thinned).  So every component 
of $P \cap F$ is essential.  Now simply attach the bottom annulus of 
$T_{r}$ to $P_{1}\cup P_{2}$ to obtain a sphere intersecting the 
original $\theta$ only in $e^{\bot}$.  Then Lemma \ref{lemma:hexagon2} 
shows $genus(K) \geq 2$.

\begin{figure}
\centering
\includegraphics[width=.8\textwidth]{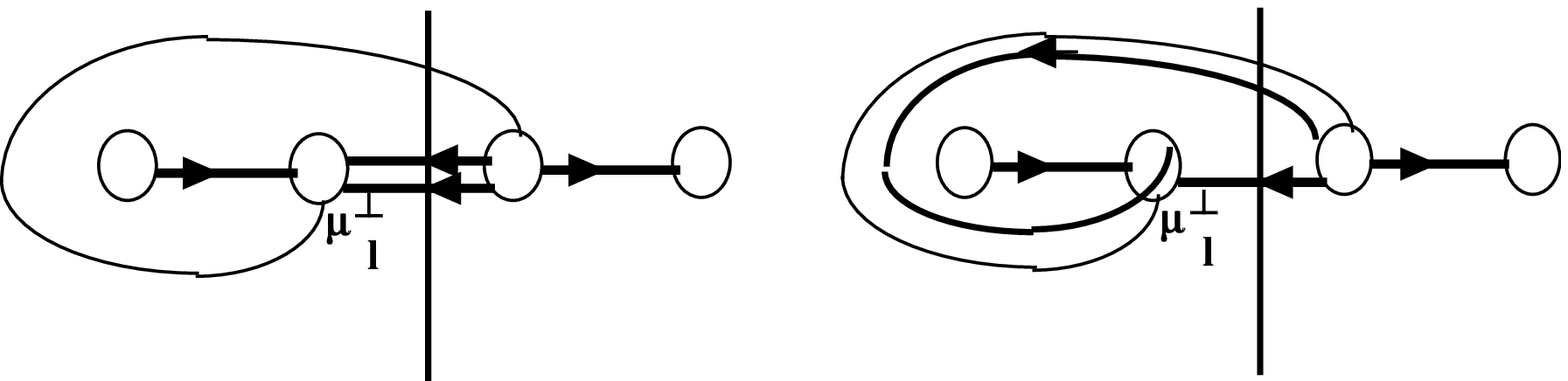}
\caption{} \label{fig:thintrick}
\end{figure}

The remaining case is if $\aaa_u$ has one end at each of the meridians 
$\mz_l$ and $\mz_r$.  Exclude any such arc that is parallel to a 
subarc of $K$, since if $\aaa_u$ were such an arc it would either 
violate the thinness of $K$ or (if $\aaa_l$ has ends at the same 
meridians) exhibit that $e^+ \cup e^{\bot}$ is vertical.  But the 
only way that $\aaa_u$ can be disjoint from $K$, connect $\mz_l$ to 
$\mz_r$, and not be parallel to a subarc of $K$ is if all $p$ arcs of 
$K$ with ends at these meridians are parallel and $\aaa_u$ is one of 
the other two paths connecting $\mz_l$ and $\mz_r$.  Although these 
paths are not parallel to a component of $K \cap \Ss$ in $\Ss$, they 
are sufficiently parallel in the component of $\eta(\theta) - P$ on 
whose boundary $\Ss$ lies to derive the same contradiction.  Here is 
the argument: If the ends of $\aaa_l$ are also at $\mz_l$ and $\mz_r$ 
then together $\aaa_r$ and $\aaa_l$ show that the cycle $e^{+} \cup 
e^{\bot}$ is vertical.  So suppose no end of $\aaa_l$ is at $\mz_l$, 
say.  One arc $\kkk$ of $K$ running from $\mz_l$ to $\mz_r$ is visibly 
isotopic to $\aaa_u$ in the $3$-ball component of $\theta - P$ on 
whose boundary $\Ss$ lies.  (See Figure \ref{fig:thintrick})  The isotopy moves the 
end of $\kkk$ across the meridian $\mz_{l}$ and so destroys the 
property that $K$ lies on $\bdd(\theta)$.  Nonetheless, once $\kkk$ 
is moved to $\aaa_{u}$ then there is no obstruction to pushing $\kkk$ 
below $P$ via $D_{u}$ while simultaneously pushing arcs of $K - P$ 
parallel to $\aaa_u$ above $P$ using $D_{l}$.  The result is a 
thinning of $K$, violating the hypothesis. \end{proof}  

\begin{lemma} \label{lemma:reguplowb}
Suppose there is a $\Th$-graph appropriate for $(K, F)$ and, for a 
thinnest one $\theta$, there is a dividing sphere $P$ that cuts off 
disjoint upper and lower disks $D_u$ and $D_l$ from $A$ so that the 
arcs $\aaa_u = D_u \cap \bdd\eta(\theta)$ and $\aaa_l = D_l \cap 
\bdd\eta(\theta)$ are regular curves disjoint from $K$.

If $e^{\bot}$ is disjoint from $P$ then either
\begin{enumerate}
	 
\item $genus(K) \geq 2$  

\item $e^+ \cup e^{\bot}$ is unknotted or

\item $p = 1$, and $e^- \cup e^{\bot}$ is unknotted.

\end{enumerate}
\end{lemma}

\begin{figure}
\centering
\includegraphics[width=.8\textwidth]{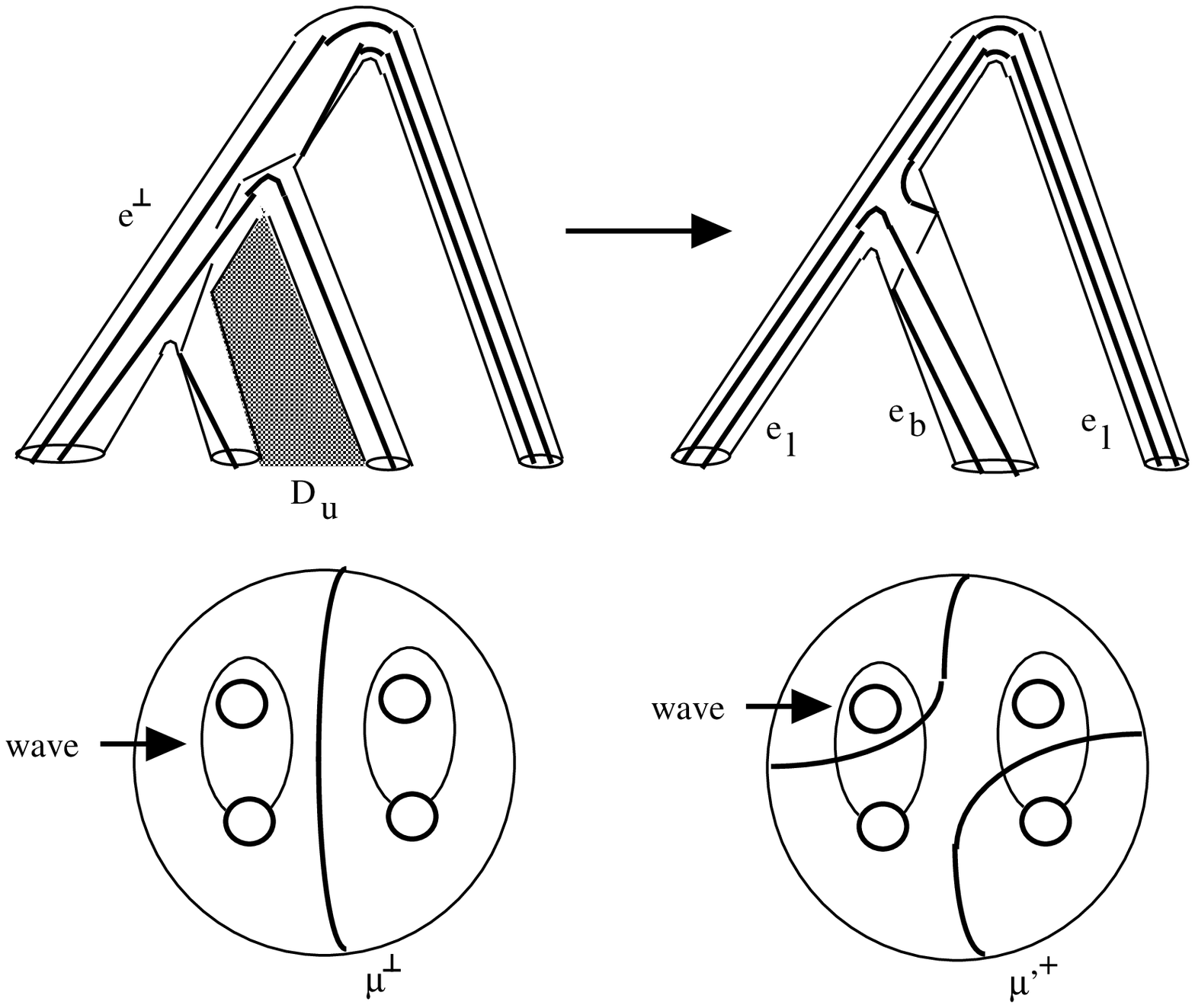}
\caption{} \label{fig:precable}
\end{figure}

\begin{proof} The proof is mostly similar to that of Lemma 
\ref{lemma:reguplow+}.  The relevant figure is modified as shown 
(Figure \ref{fig:reguplowb}), with $p-1$ arcs running between meridians 
$\mpp_{l}$ and $\mpp_r$.  There is only one regular path from 
$\mpp_{l}$ to $\mm_r$ (or from $\mpp_{r}$ to $\mm_l$) that is disjoint 
from $K$.  These paths do not cross $\mz$ and so, if $\aaa_{u}$ is 
such a path, we could use $D_u$ and $D_l$ to thin $\theta$ (without 
altering the wave condition), extending $e^{\bot}$ down to $P$.  The 
only regular path between $\mpp_{l}$ and $\mm_l$ (or $\mpp_{r}$ and 
$\mm_r$) is parallel to an arc of $K$, so $\aaa_{u}$ cannot be such a 
path.  If $p \geq 2$ and $\aaa_{u}$ runs between $\mpp_{l}$ and 
$\mpp_{r}$ we use the same argument as was used for paths from $\mz_l$ 
and $\mz_r$ previously.  

\begin{figure}
\centering
\includegraphics[width=.8\textwidth]{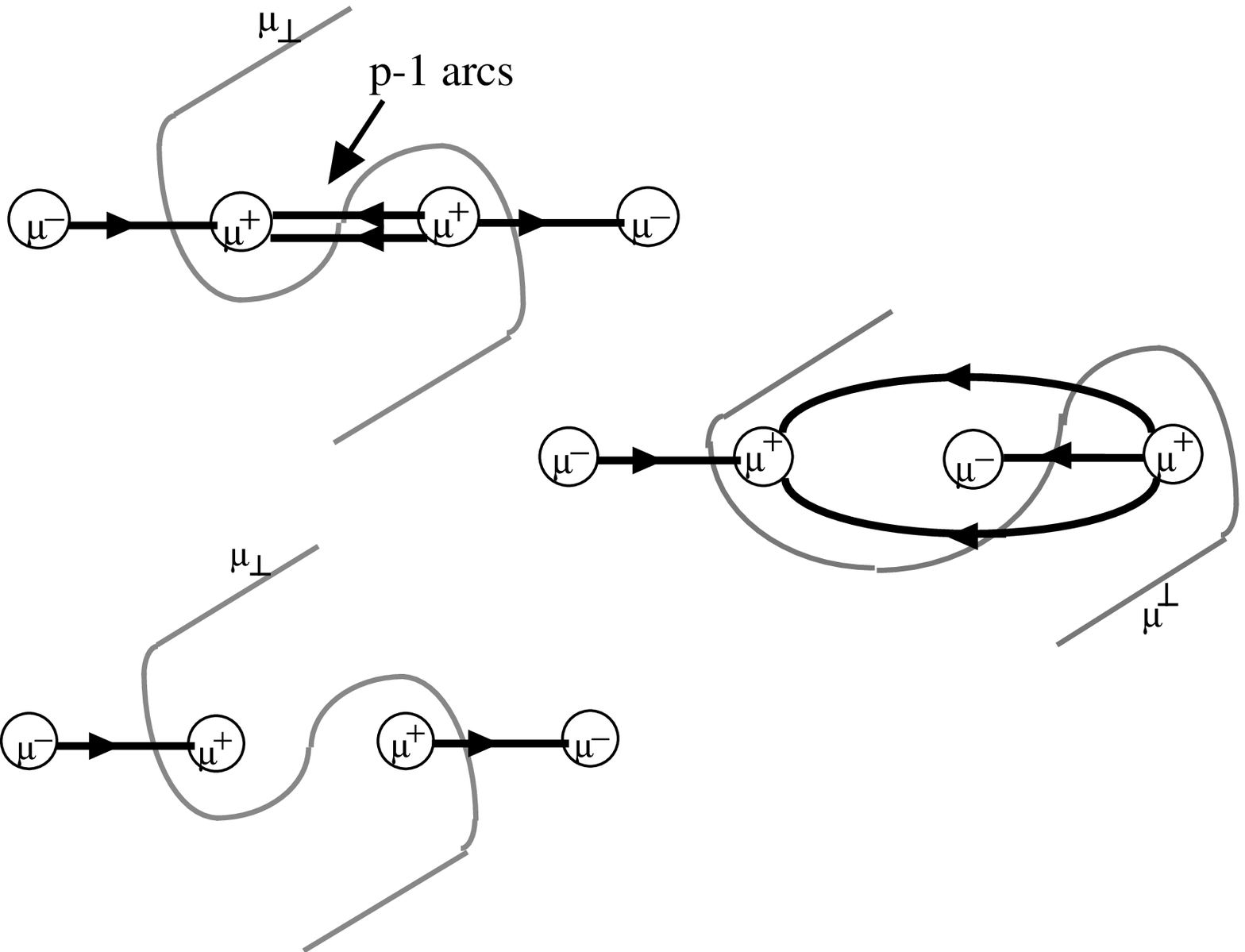}
\caption{} \label{fig:reguplowb}
\end{figure} 

Suppose finally that $\aaa_{u}$ runs between $\mm_{l}$ and $\mm_{r}$ 
(or, symmetrically, $p = 1$ and $\aaa_{u}$ runs between $\mpp_{l}$ and 
$\mpp_{r}$.)  We would like to use the same trick as was used 
previously, namely let $D_{u}$ describe a Whitney move that converts 
$\theta$ into an eyeglass graph.  There is a subtle complication, 
however.  Note that the eyeglass graph $\bowtie$ that is created by 
this move is in fact a $p$ eyeglass, not a $(p+1)$ eyeglass as before.  
In particular, the associated $\Theta$-graph $\theta'$ to $\bowtie$ is 
not $\theta$, which presented $K$ as a $(p, 1)$ quasi-cable.  Rather 
$\theta'$ presents $K$ as a $(p-1, 1)$ quasi-cable (or just as $K \cup 
\ggg$ if $p = 1$).  Nonetheless, we are still in a position to get the 
same contradiction with Corollary \ref{cor:eyetocable} (for, after 
all, $\theta$ was chosen to be thinnest among all appropriate 
$\Theta$-graphs and it was shown that such a graph is no thicker than 
$K \cup \ggg$) as long as we verify that $\theta'$ is still 
appropriate.  In other words, we need to verify that $\theta'$ still 
satisfies the wave condition.  The pre-cable disk is easy to identify: 
whereas in the previous argument it was (essentially) a thickened 
vertical arc in $D_{u}$ together with a meridian of $e^{+}$, here it 
is obtained from a thickened vertical arc in $D_{u}$ and a meridian of 
$e^{\bot}$ by adding a half-twist.  (See Figure \ref{fig:precable}.)  
Notice that $K$ intersects this meridian in $(p-1)$ points, so it is 
the meridian of $e^{+ \prime}$ (or $\ggg$ if $p = 1$) and, whether the 
wave was based at $e^{-}$ or at $e^{+}$, the extra half-twist 
guarantees that it is afterwards based at the new meridian, as 
required.  (The language when $p = 1$ is: the slope of the wave is still 
finite.)
\end{proof}

It will be useful to provide notation for arcs in the $4$-punctured 
sphere $\Ss$ discussed in the proof of Lemmas \ref{lemma:reguplow+}, 
\ref{lemma:reguplowb}.  As motivation for the notation we will use, 
suppose, as above, that the knot $K$ is thinly presented as a $(p, 1)$ 
quasi-cable on a $\Th$-graph $\theta$ in $S^3$.  Suppose further that 
$\theta$ is in bridge position with respect to a height function $h$, 
and there is a dividing sphere $P$ that is disjoint from one of the 
edges.  In particular, one of the components of $\theta - K$ is a tree 
$C$ with $4$-ends (whose regular neighborhood intersects 
$\bdd\eta(\theta)$ in the $4$-punctured sphere $\Ss$) and all the 
other components of $\theta - K$ are simple arcs.

Of course the fundamental group of $\theta$ is free on two generators.  
The natural generators of $\pi_1(\theta)$ are loops that traverse each 
of the two edges not disjoint from $P$ exactly once.  To be concrete, 
orient $K$ and suppose $e^{\bot}$ is disjoint from $P$; then a loop in 
$\theta$, based at a point in $C$ will give rise to a word in letters 
$a$ and $b$ (with inverses $\oa$ and $\ob$), where $a$ corresponds to 
traversing $e^-$ once, and $b$ to traversing $e^+$ once, each in the 
same direction as $K$.  Similarly, if one of $e^{\pm}$ is disjoint 
from $P$, then $a$ will correspond to traversing $e^{\mp}$ once and 
$b$ to traversing $e^{\bot}$ once.  We will only be interested in such 
presentations for regular simple closed curves on $\bdd\eta(\theta)$ 
that are disjoint from $K$.  In this case, the cyclic permutation 
class of the word in $a, b, \oa, \ob$ corresponding to the regular 
simple closed curve $\sss \subset \bdd\eta(\theta) - K$ determines the 
isotopy class of $\sss$ in $\bdd\eta(\theta)$ almost precisely.  The 
only ambiguity is in how $\sss$ intersects $\Ss$.

\begin{figure}
\centering
\includegraphics[width=.8\textwidth]{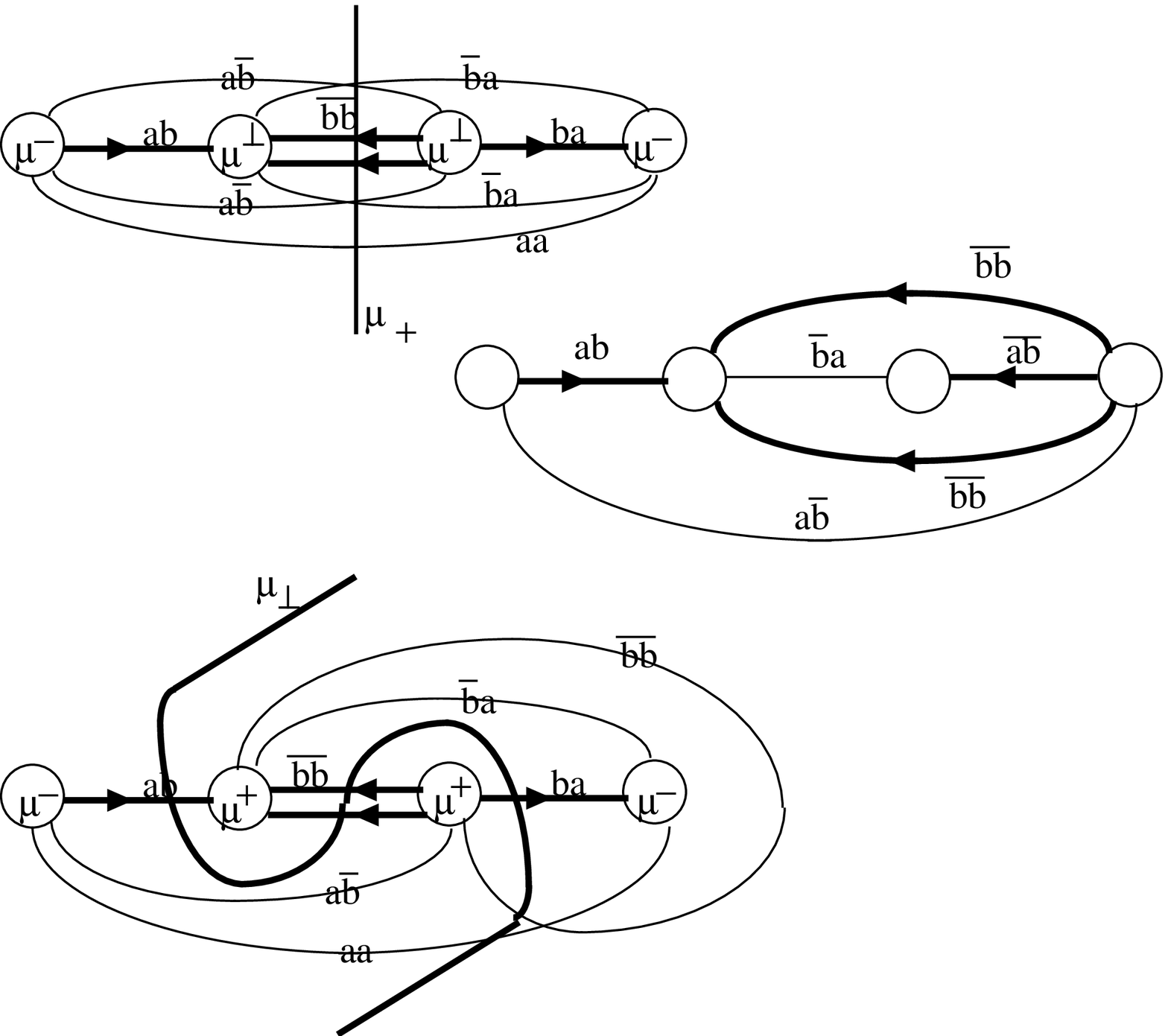}
\caption{} \label{fig:aandb}
\end{figure}

So consider isotopy classes of regular arcs in $\Ss$ that are disjoint 
from $K$.  Representative are those illustrated and labelled in Figure 
\ref{fig:aandb}.  Two of the figures show (via a heavy line from top 
to bottom) the meridian of the edge that's disjoint from $P$.  This 
meridian is relevant for determining that an arc is regular: a regular 
arc can cross the meridian at most once.  The meridians of the arcs 
intersecting $P$, i.  e.  the boundary components of $\Ss$, are shown 
as circles.  The arcs are labelled by where they would occur in a word 
in $a, b, \oa, \ob$; thus if the sequence $\ldots a\ob\ldots$ (or, 
inversely, $\ldots b\oa \ldots$) occurs in the word, the corresponding 
arc in $\Ss$ would be one labelled $a\ob$.  Notice, as one example, 
that there are two arcs labelled $\ob a$, indicating ambiguity 
in how such an arc may run through $\Ss$.  There is less ambiguity in 
the lowest figure (corresponding to the edge $e^{\bot}$ disjoint from 
$P$).  The heavy oriented horizontal curves correspond to arcs in $K$, 
and they are labelled in the same manner.  In the first picture, to 
avoid crowding only one arc labelled $\ob\ob$ (or, equivalently, $bb$) 
is shown; it's a subarc of $K$.  One of the two others is shown in the 
lowest figure; the other arc labelled $\ob\ob$, which is only relevant 
to the first figure, is obtained from the second by reflection through 
the vertical arc $\mu^+$.  The special case $p = 1$ and $e^{\bot}$ 
disjoint from $P$ is not shown.

With this labelling, another way of stating the first possibility in 
Lemma \ref{lemma:reguplow+} would then 
be: $\aaa_u$ is of type $a\ob$ or $\ob a.$

Under the hypotheses of the lemma, it is natural to define

\begin{defin}

The {\em geometric length} $l(w)$ of a word $w$ in $a, b, \oa, \ob$ 
representing a regular loop $\sss \subset \bdd\eta(\theta)$ is $|\sss 
\cap P|$, for $P$ a dividing sphere.
	
\end{defin}

Given $l(a)$ and $l(b)$ the geometric length of $w$ is clearly $ml(a)
+ nl(b)$, where $m$ and $n$ are the total number of occurences of,
respectively, $a$ and $\oa$, $b$ and $\ob$.  We will say that a word
$w$ is {\em positive in $a, b$} if neither $\oa$ nor $\ob$ occur in
$w$ (e. g. when $w = \emptyset$).  We then have

\begin{cor}  \label{cor:reguplow0}
Suppose there is a $\Th$-graph $\theta$ appropriate for $(K, F)$ and 
$\theta$ is a thinnest one.  Suppose further that

\begin{itemize}

\item $e^+$ is disjoint from a dividing sphere.  

\item there is a properly imbedded normal form annulus $A$ in $S^3 - 
\eta(\theta)$ whose boundary components $\bdd_{\pm} A$ are disjoint 
from $K$ and regular in 
$\bdd\eta(\theta)$ and

\item a dividing sphere $P$ cuts off upper and lower disks from $A$ and 
the component(s) of $\bdd A$ to which these disks are incident 
represents a word that is positive in $a, b$.

\end{itemize}

Then either $e^+ \cup e^{\bot}$ is unknotted or $genus(K) \geq 2$

\end{cor}

\begin{proof} The upper disk cannot be of type $a\ob$ or $\ob a$ 
because a word that is positive in $a, b$ contains neither $\ob$ nor 
$\oa$.  The result then follows from Lemma \ref{lemma:reguplow+}.
\end{proof}

Similarly, from Lemma \ref{lemma:reguplowb} we derive this corollary

\begin{cor}  \label{cor:reguplowb} Suppose there is a 
$\Th$-graph $\theta$ appropriate for $(K, F)$ and $\theta$ is a thinnest 
one. Suppose further that 

\begin{itemize}

\item $e^{\bot}$ is disjoint from a dividing sphere

\item there is a properly imbedded normal form annulus $A$ in $S^3 - 
\eta(\theta)$ whose boundary components $\bdd_{\pm} A$ are disjoint 
from $K$ and regular in $\bdd\eta(\theta)$ and

\item a dividing sphere $P$ cuts off upper and lower disks from $A$.

\end{itemize}

Then either $e^{\bot} \cup e^{+}$ is unknotted or $genus(K) \geq 2$ or 
$p = 1$ and $e^{\bot} \cup e^{-}$ is unknotted.
\end{cor}

\section{When $K$ has genus one}

Our interest in regular annuli comes from the following observation. 
Suppose $F$ is genus one.  Then an outermost disk cut off by $F \cap
E$ in $E$ provides a boundary-compression of $F$ to an essential
annulus $A \subset S^3 - \eta(\theta)$.  When viewed on
$\bdd\eta(\theta)$ the relation between $\bdd A$ and $K$ is this: $K$
is banded to itself via a subarc $\omega \subset \bdd E - K$ whose
ends are incident to the same side of $K$.  If $D \subset
\eta(\theta)$ satisfies the wave condition then any subarc of $\bdd
E = \bdd D$ that is disjoint from $K$ is regular.  So we are assured
that $\bdd A$ is regular in $\eta(\theta)$.

This leads to 

 \begin{theorem} \label{theor:eplus}
Suppose in a thinnest $\Th$-graph $\theta$ appropriate for $(K, F)$, 
the edge $e^{+}$ is disjoint from a dividing sphere.  Suppose further 
that $genus(F) = 1$ and that the wave for $\bdd D$ is based at $\mm$.  
Then either the cycle $e^+ \cup e^{\bot}$ is unknotted or $p = 1$ and the 
cycle $e^- \cup e^{\bot}$ is unknotted.
\end{theorem}

\begin{proof} That $q = 1$ follows from Proposition 
\ref{prop:oneonethin2}.  We will show that the annulus $A$ obtained 
from $\bdd$-compressing $F$ to $\eta(\theta)$, using the disk $E$ from 
the splitting sphere, gives rise to an annulus satisfying the 
conditions of Corollary \ref{cor:reguplow0}.  The result then follows 
from that corollary, possibly by way of Proposition 
\ref{prop:oneonethin2}.

Let $\omega \subset 
\bdd\eta(\theta) - K$ be the arc described above, which we may as well 
slide to minimize intersections with the meridians of $\eta(\theta)$.  
In particular, for $P$ a level sphere between the lowest maximum and 
the highest minimimum of $\theta$, the ends of $\omega$ will lie on 
the $4$-punctured sphere component $\Ss$ of $\bdd\eta(\theta) - P$.  
Let $w$ be the word in $a, b, \oa, \ob$ represented by $\omega$.  
Because the wave of $\bdd D$ is based at $\mm$

\begin{itemize}
	
	\item any occurence of the letter $b$ (resp.  $\ob$) in $w$ is 
	preceded and followed by the letter $a$ (resp. $\oa$).

	\item any occurence of the letter $a$ in $w$ is 
	followed by $a$ or $b$ and preceded by $a$ or $b$ and
	
	\item any occurence of the letter $\oa$ in $w$  is followed by $\oa$ or $\ob$
	and preceded by $\oa$ or $\ob$.
	
\end{itemize}

In particular, by a choice of orientation for $w$, we can assume that 
$w$ is positive (say) in $a$ and $b$ (e.  g.  perhaps $w = 
\emptyset$).  Exploiting these facts, together with the symmetries of 
the diagram, we have three essentially different ways in which the 
ends of $w$ can lie in $\Ss$.  These are shown in Figure \ref{fig:wend}.

\begin{figure}
\centering
\includegraphics[width=.5\textwidth]{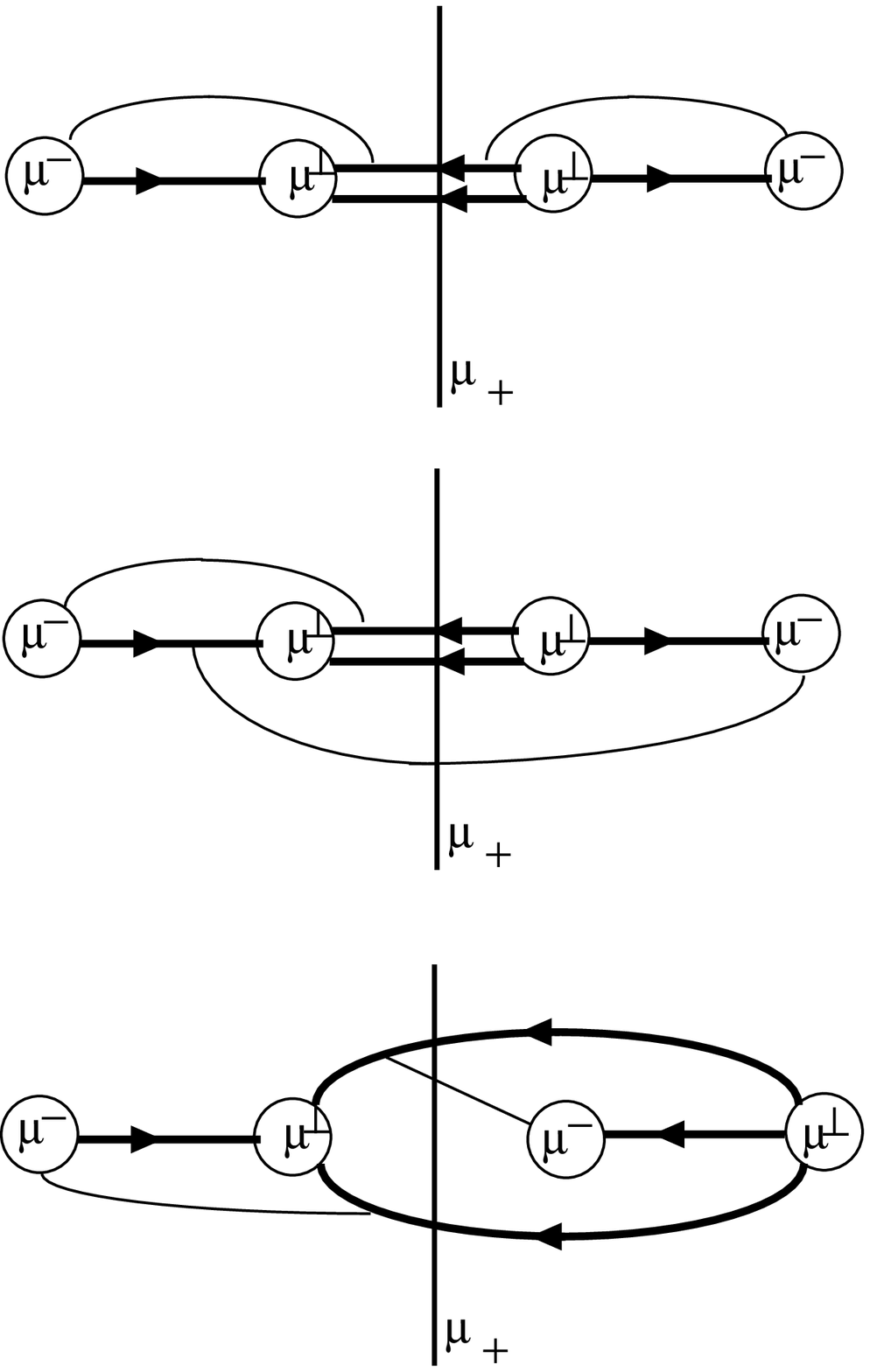}
\caption{} \label{fig:wend}
\end{figure} 

Now orient $w$ from left to right, and read off the words 
corresponding to $\bdd_{\pm}A$ (boundaries oriented to be parallel in 
$A$, not antiparallel).  All three cases can be expressed in one of 
the following two forms, with the details depending on where the ends 
of $\omega$ are incident to the word given by $K$, namely $ab^{p+1}$.  
Note that, in all cases, the word $w$ begins and ends with the letter 
$a$.  The choice of labelling of the components $\bdd_{\pm} A$ of 
$\bdd A$ is made so that the word corresponding to $\bdd_{+} A$ is 
positive in $a, b$.

\bigskip

{\bf Form 1:}

\begin{itemize}
	
	\item $\bdd_+ A \leftrightarrow wb^{i}$
	
	\item $\bdd_- A \leftrightarrow w\ob^{j}\oa\ob^{k}$
	
\end{itemize}

Here $j, k > 0$ and $i+j+k = p+1$. 

\bigskip

{\bf Form 2:} 

\begin{itemize}
	
	\item $\bdd_+A \leftrightarrow wb^{j}ab^{k}$
	
	\item $\bdd_-A \leftrightarrow w\ob^{i}$
	
\end{itemize}

Here $j > 0$, $i, k \geq 0$ and $i + j + k = p+1$.  

\bigskip

Note that the geometric length of $\bdd_{+} A$ is greater than that of 
$\bdd_{-} A$ exactly when, for annuli of the first form, $(i -j - 
k)l(b) > l(a)$ and, for annuli of the second form, $(i - j - k)l(b) < 
l(a)$.  In either case, $\bdd_{+} A$ contains an occurence of the 
letter $b$, so $\bdd_{+} A$ runs along $e^{\bot}$.  If $e^{\bot}$ 
intersects a dividing sphere $P$ just twice, then $e^{+} \cup 
e^{\bot}$ is vertical, and we are done.  If $|e^{\bot} \cap P| \geq 
4$, then a regular maximum and minimum (which we can take, 
respectively, to be the lowest maximum and the highest minimum) lie in 
the interior of $e^{\bot}$.  It follows, then, from Lemma 
\ref{lemma:unbalanced} and Corollary \ref{cor:reguplow0} that we are done if 
$(i -j - k)l(b) > l(a)$ for annuli of the first form, or $(i - j - 
k)l(b) < l(a)$ for annuli in the second form.  So we now consider only 
the alternative possibilities.
\bigskip

{\bf Case 1:} $(i - j - k)l(b) = l(a)$ 

In this case, it follows roughly from the same argument used in Lemma 
\ref{lemma:unbalanced} that there is a dividing sphere $P$ so that no 
arc of $P \cap A$ cuts off an outermost disk incident to $\bdd_+ A$.  
Indeed, as above, a high dividing sphere cuts off an upper disk 
incident to $\bdd_{+} A$ and a low dividing sphere cuts off a lower 
disk incident to $\bdd_{+} A$.  There can't be both an upper and a 
lower such disk incident to $\bdd_+ A$ by Corollary 
\ref{cor:reguplow0}.  So at some level there is a dividing sphere $P$ 
that cuts off no outermost disk incident to $\bdd_{+} A$.  But since 
the geometric lengths of $\bdd_{\pm} A$ are equal, this implies that 
$P$ intersects $A$ only in spanning arcs.

We will argue that this is impossible.  Suppose (with no loss) that 
$e^{+}$ and hence $\Ss$ lie above (not below) a dividing sphere.  The 
spanning arcs determine a correspondence between subintervals of 
$\bdd_-A$ and $\bdd_+A$.  To be precise, say that a component of 
$\bdd_{-}A - P$ and a component of $\bdd_{+}A - P$ are {\em opposite} 
if there is a (``square'') component of $A - P$ incident to both.  
More generally, a segment of $\bdd_+ A$ and a segment of $\bdd_{-} A$ 
are opposite each other if a spanning arc of $P \cap A$ runs between 
the beginning of each and another spanning arc runs between the end of 
each.  For example, we first observe that no segment of $\bdd_{-} A$ 
that is part of a $\ob$ interval (meaning that it comes from an 
occurence of a letter $\ob$, i.  e.  that it runs along $e^{\bot}$ 
with an orientation opposite to that of $K$) can lie opposite a 
segment of $\bdd_{+} A$ that is part of a $b$ interval.  For if this 
occured then it is easy to see that somewhere on the entire length of 
the $b$ and $\ob$ intervals there would be components of $\bdd A - P$ 
that are opposite to each other on $A$ but lie on the same component 
$\kkk$ of $e^{\bot} - P$.  Since they have opposite orientation in 
$\bdd_- A$ and $\bdd_+ A$ the square component of $A - P$ connecting 
them can be attached to the punctured solid torus $P \cup \eta(\kkk)$ 
to create a punctured Lens space $L(2, 1) \subset S^{3}$, which is 
absurd.  Similarly no segment of $\bdd_{-} A$ that is part of an $\oa$ 
interval can lie opposite a segment of $\bdd_{+} A$ that is part of an 
$a$ interval.  This immediately rules out the first form above (again, 
only under the assumption that $(i -j - k)l(b) = l(a)$) since, 
following these observations, the only possible segment opposite the 
transition segment from $\ob$ to $\oa$ in $\bdd_- A$ would be exactly 
a transition segment from $a$ to $b$ in $\bdd_+ A$, and that would 
lead to the same contradiction.

Ruling out the second form (in which no letter $\oa$ appears) is only 
slightly more complicated.  We will focus on the segments of $\bdd A - 
P$ that lie on the middle component of $e^{\bot} - P$; that is, on the 
arc component that is equidistant (measuring distance by intersection 
with $P$) from both ends of $e^{\bot}$.  Note that this segment of 
$e^{\bot} - P$ lies below $P$.  Label corresponding segments of $\bdd 
A - P$ by $\kkk$.  Note that none of these can be opposite to a 
segment of $\bdd A$ lying on $\Ss$ (e.  g.  those segments in $\bdd A$ 
that correspond to the transition between different letters) since the 
segments on $\Ss$ lie above $P$.  This remark allows us to be a bit 
casual about length arguments in the next few paragraphs, since it 
means that inequalities will usually imply strict inequalities.

The first observation is that no label $\kkk$ occurs opposite to any part 
of a $\ob^p$ interval in $\bdd_- A$, for this would allow us to display 
a Lens space $L(2,1)$ in $S^3$, as noted above.

Let $\sss \subset \bdd_- A$ be the segment between the first and last 
labels $\kappa$ in $\ob^i \subset \bdd_{-} A$.  Then $l(\sss) = 
(i-1)l(b) \geq l(a)$.  Since no label $\kkk$ lies opposite to $\ob^p$ 
it follows that opposite to $\sss$ is part of a segment in $\bdd_+ A$ 
corresponding to $a^q, q \geq 1$.  Notice that if a label $\kkk$ in 
$\bdd_- A$, say, is opposite to any part of a $b$ interval in $\bdd_+ 
A$, than the relation is reciprocal: the label $\kkk$ in the $b$ 
interval on $\bdd_+ A$ is opposite to the $b$ interval in $\bdd_- A$ 
containing the original label $\kkk$.  (This is not deep, just a 
reflection that we have taken $\kkk$ to lie half way along $e^{\bot}$ 
and so it appears half way along each $b$ or $\ob$ interval.)  Because 
$j + k > 0$ there is at least one more label $\kkk$ in $\bdd_+ A$ then 
there are $\kkk$ labels in $\bdd_- A$, not counting the labels $\kkk$ 
in $\ob^i$.  It follows that some label $\kkk$ in $\bdd_+ A$ is 
opposite a part of an $a$-interval in $\bdd_- A$, so an entire half of 
a $b$-interval in $\bdd_+ A$ is opposite to a subsegment $\ttt$ of a 
single $a$-interval in $\bdd_- A$, for $l(a) \geq l(b)$.  And, as 
we've seen, another copy of $\ttt$ lies opposite to a subsegment of 
$\sss$.  This works just as well for construcing a Lens space in $S^3$ 
as having the half of the $b$-segment itself opposite to $\sss$.  (Two 
rectangles are glued together along the boundary interval they share 
on a component of $e^{-} - P$ corresponding to part of the $a$ 
intervals.)  So we arrive at the same contradiction as previously.

\bigskip

{\bf Case 2:} For an annulus of the first form, $(i -j - k)l(b) < l(a)$ or, 
for one of the second form, $(i -j - k)l(b) > l(a)$.  

\begin{figure}
\centering
\includegraphics[width=.7\textwidth]{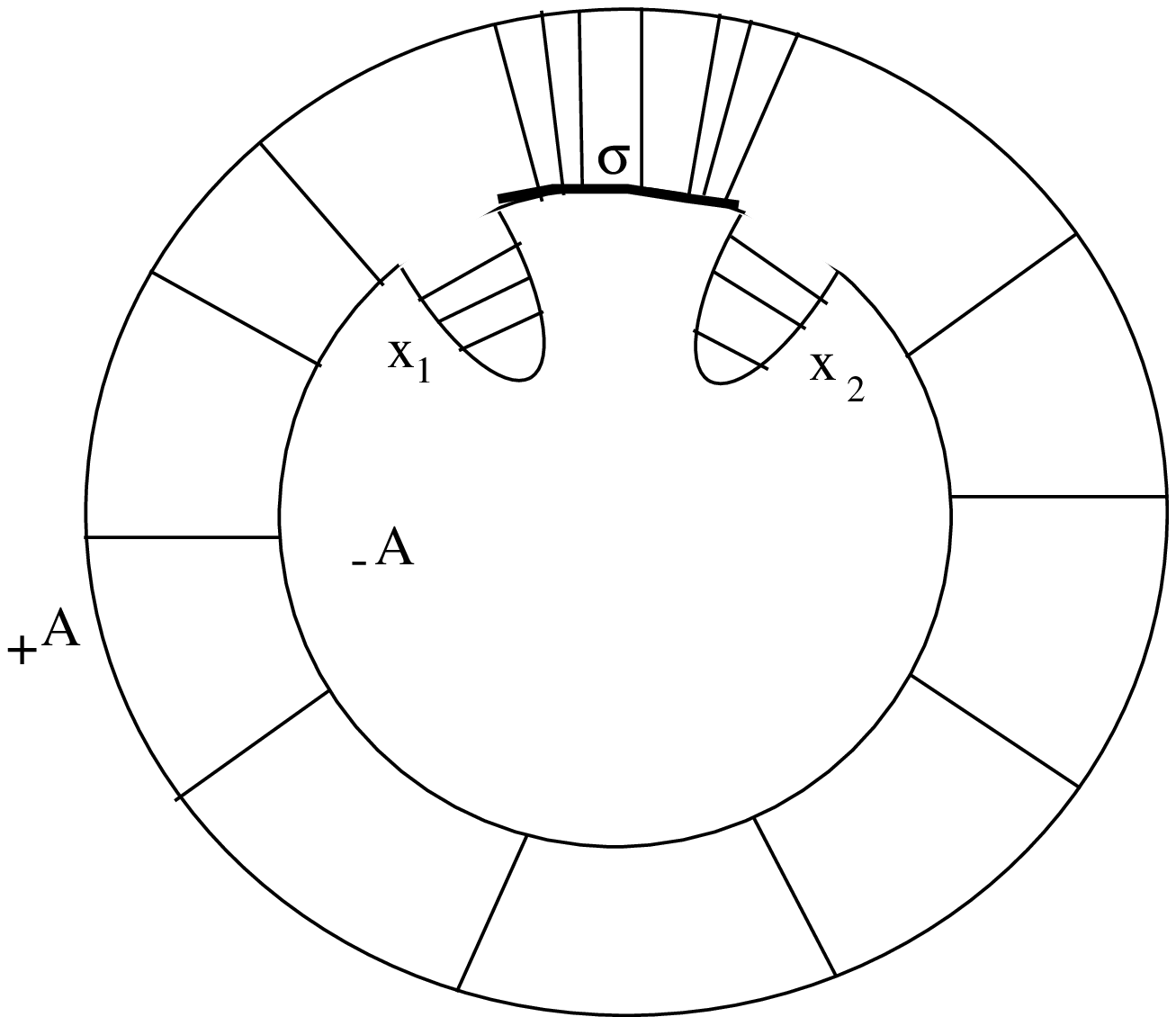}
\caption{} \label{fig:xannulus}
\end{figure} 

We will arrive at the same sort of contradiction, though the argument 
is a bit more complicated.  Again, with no loss, we assume that 
$e^{+}$ lies above a dividing sphere and that both the lowest maximum 
and the highest minimum of $\theta$ are regular critical points on 
$e^{\bot}$.  This implies that just below the lowest maximum (i.  e.  
at a high dividing sphere) the dividing sphere cuts off an upper disk 
from $A$ that is incident to a regular maximum (namely the lowest 
maximum).  Then no outermost disk cut off from $A$ by this high 
dividing sphere can be a lower disk, by thin position.  On the other 
hand, a low dividing sphere does cut off a lower disk from $A$.  So 
there is a height $y$ such that a dividing sphere just below $y$ cuts 
off a lower disk from $A$ but just above $y$ a dividing sphere cuts 
off no lower disk.  But any dividing sphere must cut off some 
outermost disk, since the geometric lengths of the words represented 
by the two boundary components of $A$ are different.  It follows from 
Corollary \ref{cor:reguplow0} then that either $e^{+} \cup e^{\bot}$ 
is unknotted or just above $y$ all outermost disks cut off of $A$ by 
the level sphere $P$ are incident to $\bdd_{-} A$ and, moreover, for 
each such disk $D_{u}$ the arc $\aaa_{u} = \bdd D_{u} \cap \bdd H$ is 
of the form $a\ob$.  That is, there are at most two outermost disks in 
$A$ and they are incident to the subarcs $\aaa_1$ and $\aaa_2$ 
labelled $a\ob$ and $\ob a$ of $\bdd_-A$.

In this position, the total number of non-spanning arcs in $A$, all of 
them incident to $\bdd_- A$ and each of them cutting off a disk 
containing either $\aaa_1$ or $\aaa_2$, is $x = |(i-j-k)l(b) - 
l(a)|/2$, since each arc has two ends.  Since $x$ is less than the 
distance between the ends of $\ob^{j}\oa\ob^{k}$ (first form) and less 
than the distance between the ends of $\ob^{i}$ (second form), each 
non-spanning arc cuts off a disk containing exactly one of $\aaa_1$ or 
$\aaa_2$ and so each arc is parallel to one of the $\aaa_i$.  Let 
$x_i$ denote the number parallel to $\aaa_i$, so $x_1 + x_2 = x$.  Let 
$\sss$ denote the segment in $\ob^p$ that is still incident to 
spanning arcs, and let $s$ be its length.  See Figure 
\ref{fig:xannulus}.  For obvious pictorial reasons, we'll refer to the 
part of $A$ containing the collection of arcs parallel to $\aaa_{i}$ 
as the {\em $x_{i}$ peninsula.}

\bigskip

{\bf Subcase 2a:}  $\bdd A$ is of the first form.

In this case note that $s + x = (j+k)l(b) + l(a)$, so $s = x + 
il(b)$.

\medskip

{\bf Subcase 2a.i:}  The entire $\oa$ interval of $\bdd_{-} A$ 
is disjoint from $\sss$.

Say the entire $\oa$ interval lies on the $x_1$ peninsula.  Then 
$\sss$ is entirely made up of powers of $\ob$ and its length is at 
least $x_1 + il(b) \geq l(a) + (i+j)l(b) \geq l(a) + 2l(b)$.  Since 
$\sss$ is made up of $\ob$-intervals, at most $l(b)$ of the length of 
the segment opposite to $\sss$ can lie in $b$ intervals (half at each 
end), so in particular, the segment opposite to $\sss$ contains an 
$a$-segment longer than $l(a)$.  In particular, if $\ttt$ denotes the 
terminal segment of $a$ at the end of the $x_1$ peninsula, there is 
also a copy of $\ttt$ lying opposite $\sss$.  Since across the ends of 
the $x_1$ peninsula the orientations of $a$ and $b$ coincide (that is, 
the orientations of $a$ and $\ob$ are reversed by the folding along 
$x$) whereas across from $\sss$ they disagree, we obtain the standard 
Lens space contradiction.  See Figure \ref{fig:sub2}i

\begin{figure}
\centering
\includegraphics[width=.8\textwidth]{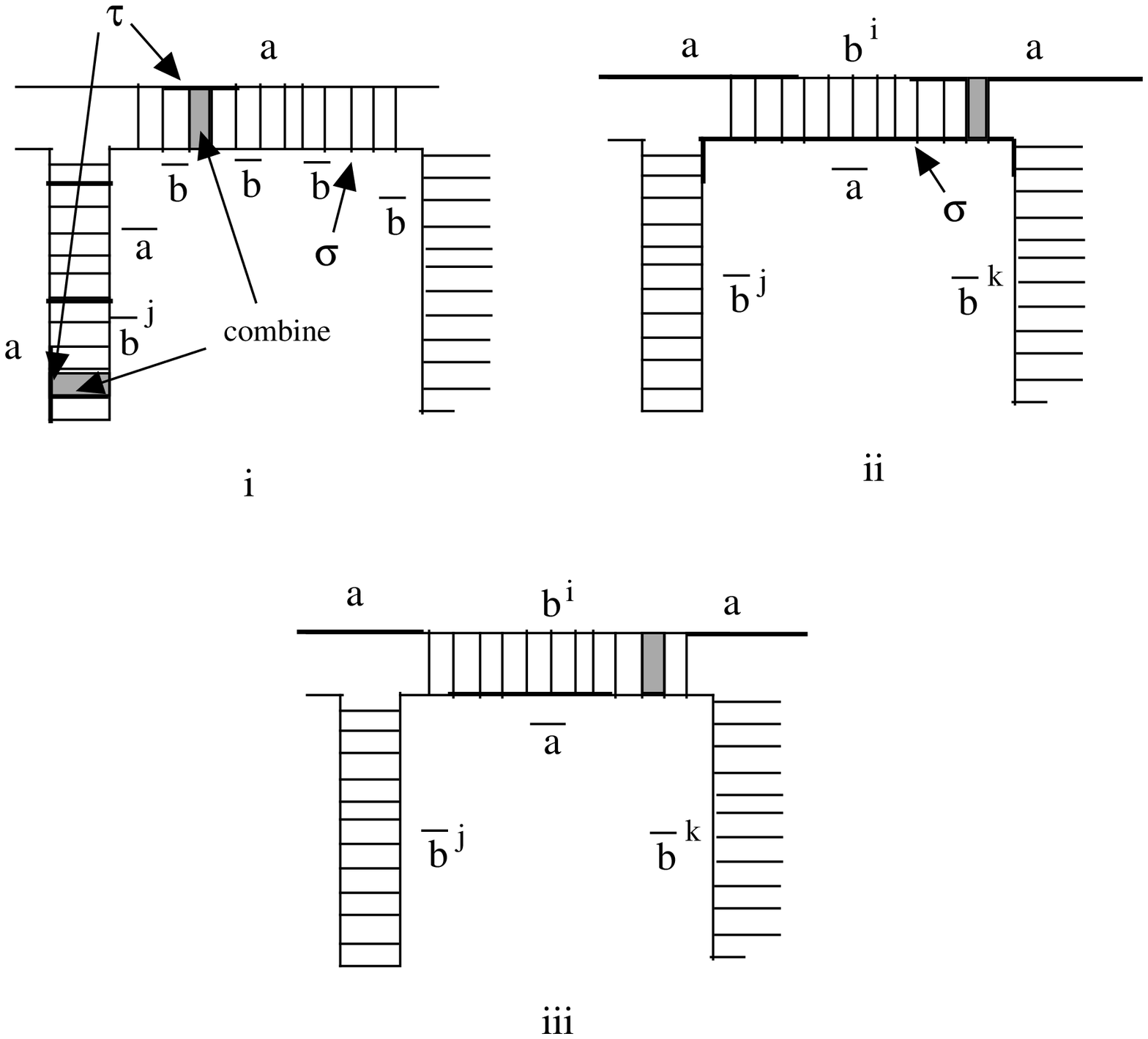}
\caption{} \label{fig:sub2}
\end{figure} 

\medskip

{\bf Subcase 2a.ii:} $\sss$ is completely contained in the $\oa$ 
interval.  

We have $s = x + il(b)$ and the longest $b$-segment in $\bdd_{+}A$ is 
of length $il(b)$, since the wave assumption ensures there are no 
proper powers of $b$ in $w$ (i.  e.  no powers greater than one).  
Hence the total length of $a$-segment(s) opposite $\sss$ is at least 
$x$.  In particular, the length of the $a$-segment(s) across from 
$\sss$ must be greater (by at least $(j+k)l(b)$) than the length of 
the ends of $a$ not contained in $\sss$.  This implies that some 
subsegment of $\sss$ is opposite a copy of itself, leading to the 
standard Lens space contradiction.  See Figure \ref{fig:sub2}ii.
 
 \medskip

{\bf Subcase 2a.iii:} The $\oa$ interval is completely contained in $\sss$. 

The argument of the previous subcase applies as long as $s < l(a) + 
l(b)$, so assume that $s \geq l(a) + l(b).$ To avoid the standard Lens 
space contradiction, across from the $\oa$ interval in $\sss$ is a 
segment comprised entirely of $b$-intervals.  The largest power 
of $b$ is $b^i$ and $\sss$ is even longer than that, so at least one 
end of $\oa$ is across from a copy of $b$ that lies completely in the 
segment across from $\sss$.  If the ends of $\oa$ and that copy of $b$ 
coincide, we get a Lens space contradiction via the terminal segment 
of $x$.  If the copy of $b$ extends out beyond $\oa \subset \sss$ then 
we get a Lens space contradiction with the end of $\ob$ adjacent to 
$\oa$ in $\sss$.  See Figure \ref{fig:sub2}iii.

\medskip

{\bf Subcase 2a.iv:} One end of $\oa$ is contained in $x_1$, say, and the 
other end is contained in $\sss$.  

Suppose first that the segment opposite the end of $\sss$ at $x_1$ is 
part of a $b$-interval.  Since the longest $b$-intervial comes from 
$b^{i}$ and $s = x +il(b)$ it follows that the $b$-interval ends 
somewhere in $\sss$ and is followed by an $a$-interval.  We get the 
standard Lens space contradiction with either $\oa$ or $\ob$, 
depending on whether the $b$-interval ends across from a point in 
$\oa$ or a point in $\ob$.  See \ref{fig:subiv}i.

Next suppose the segment opposite the end of $\sss$ at $x_1$ is part 
of an $a$-interval and let $\aaa$ denote that part of the {\em single} 
$a$-interval that lies across from $\sss$.  (So $\aaa$ is followed 
either by a $b$-interval or another $a$-interval.)  Abusing notation 
somewhat, let $\oa \cap x$ denote that part of the $\oa$ interval that 
lies on the $x_{1}$ peninsula.  If $\aaa$ is longer than $\oa \cap x$ 
we get a Lens space contradiction between $\aaa$ and $\oa \cap \sss$.  
See \ref{fig:subiv}ii.

If $\aaa$ is shorter than $\oa \cap x$ then, since $s = x + il(b)$ and 
$b^i$ is the highest power of $b$, there is more $a$-segment across 
from $\sss$ than just $\aaa$.  If there are some $b$ intervals between 
$\aaa$ and the additional $a$-segment, then the far (right-hand in the 
figure) end of the $b$-segment gives the same Lens space contradiction.  
So we conclude that $\aaa$ is immediately followed by another copy of 
$a$, which we'll call $a_1$.  

If $\aaa$ is shorter than $\oa \cap \sss$, as must happen if most of 
$\oa$ lies in $\sss$, we get a Lens space contradiction between 
$a_{1}$ and the end of $\oa \cap \sss$.  See \ref{fig:subiv}iii.  If 
$\aaa$ is longer than $\oa \cap \sss$ we get a Lens space 
contradiction, comparing the end of $\aaa$ across from a $\ob$ segment 
with the end of the $x_{1}$ peninsula.  See \ref{fig:subiv}iv.

\begin{figure}
\centering
\includegraphics[width=.8\textwidth]{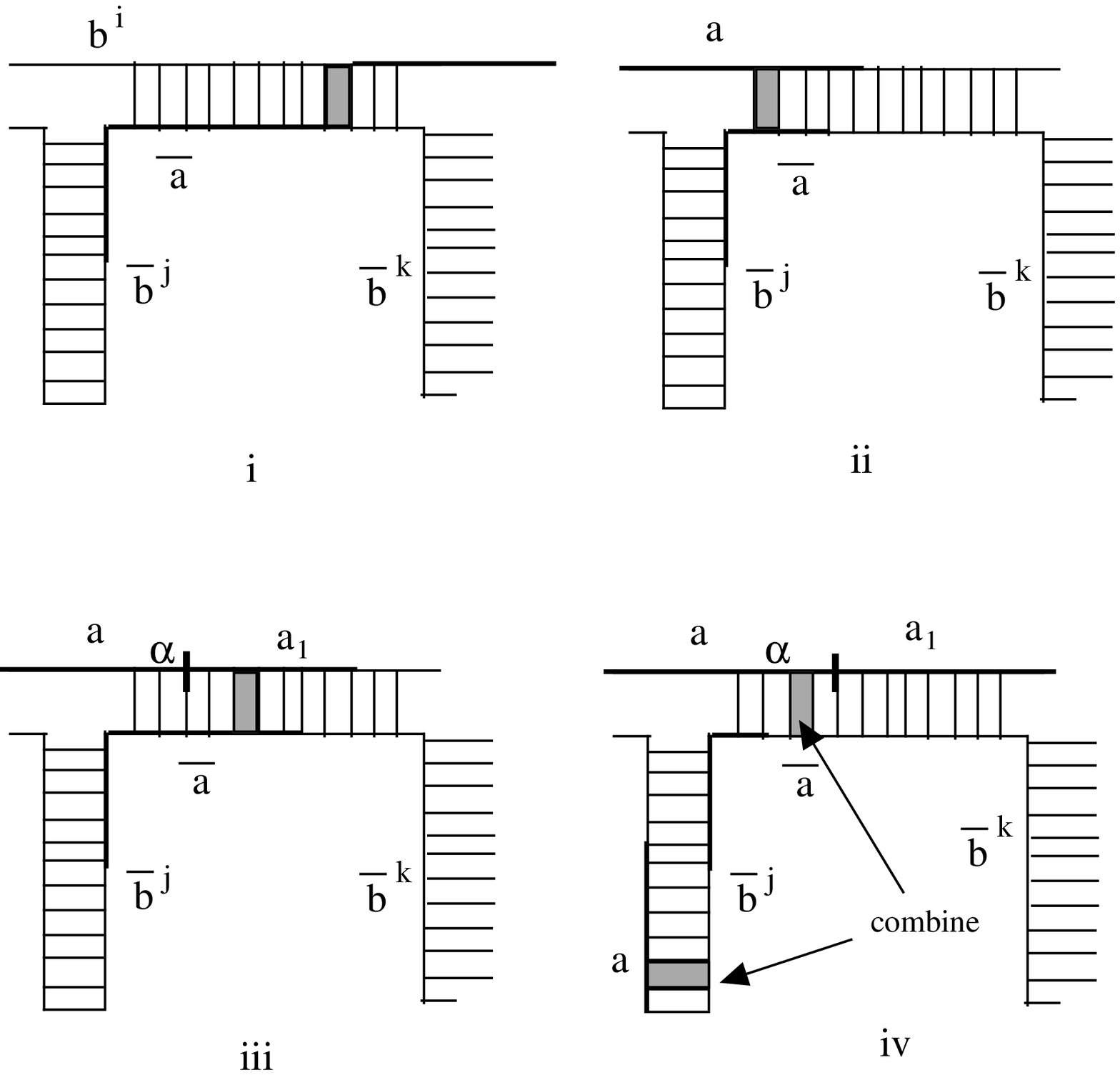}
\caption{} \label{fig:subiv}
\end{figure} 

\bigskip

{\bf Case 2b}  $\bdd A$ is of the second form.

Observe then that $s+x = il(b)$ so $s = x + (j+k)l(b) + l(a) > l(b) + 
l(a)$.  It follows that the distance between the outermost labels 
$\kkk$ (see Case 1) is greater than $l(a)$.  The rest now follows 
almost exactly as for the second form in Case 1.

\end{proof}

When the wave is based at $\mpp$ the result is less ambitious.  We 
need the following lemma:

\begin{lemma} \label{lemma:noshort}
If $|\bdd A \cap \mm| = 1$ then $e^{\bot} \cup e^+$ is unknotted.  
Symmetrically, if $|\bdd A \cap \mpp| = 1$ then $e^{\bot} \cup e^-$ is 
unknotted.
\end{lemma}

\begin{proof} Suppose $|\bdd A \cap \mm| = 1$.  Take two parallel 
copies of $\mm$ and band them together along the part of 
$\bdd_- A$ that does not lie between them.  The result is a disk $E 
\subset H$ that is disjoint from $\bdd A$ and separates $H$, leaving 
one of $\bdd_{\pm} A$ in the boundary of each of the solid tori 
components of $H - E$.  Label these solid tori (correspondingly) 
$L_{\pm}$ and denote by $L$ the link whose core circles are $L_- \cup 
L_+$.  Note that $\bdd_- A$ is a longitude of $L_-$ and $\bdd_+ A$ is 
a $(p, q)$ cable of $L_+$, some $q$.  $L$ is visibly a non-hyperbolic 
(because of $A$) tunnel number one link (the tunnel is dual to $E$).  
These have been classified (cf.  \cite{EU}): In particular, $L_+$ is 
the unknot.  But the core of $L_+$ is $e^{\bot} \cup e^+$, as required.

If $|\bdd A \cap \mpp| = 1$ (so $p = 1$) the argument is symmetric, 
interchanging $\mm$ and $\mpp$.
\end{proof}

\begin{prop} \label{prop:eplus2}  
Suppose in a thinnest $\Th$-graph $\theta$ appropriate for $(K, F)$, 
the edge $e^+$ is disjoint from a dividing sphere.  Suppose also 
that $F$ is of genus one and that the wave for $\bdd D$ is based at 
$\mpp$.  Then either $e^{\bot} \cup e^{+}$ is unknotted, or a Whitney 
move on $\theta$ changes it to an equally thin $\Th$-graph $\theta'$ 
that is appropriate for $(K, F)$.  

If $\theta$ presented $K$ as a 
$(p, 1)$ quasi-cable, then $\theta'$ presents it as a $(p+1, 1)$ 
quasi-cable.
\end{prop}

\begin{proof}
With no loss of generality, assume that $e^{+}$ lies above the 
dividing sphere and that $e^{+}$ is monotonic.

The proof now has the same features as the proofs of Theorems 
\ref{theor:eplus} and we use similar notation.  
Let $\omega \subset \bdd\eta(\theta) - K$ be the arc as previously, 
again slid to minimize intersections with the meridians of 
$\eta(\theta)$.  Let $\Ss$ again be the $4$-punctured 
sphere lying in $\bdd\eta(\theta) - P$, on the neighborhood of the 
component of $\theta - P$ that lies above $P$ and contains $e^{+}$. 

For the purposes of the argument, we will assume that all $p$ arcs of 
$K \cap \Ss$ that run between the two copies of $\mz$ in $\bdd \Ss$ 
are parallel.  If in fact there are two families of parallel arcs, 
the argument is essentially identical, except for one difference 
which is noted below.  

Then $K \cap \Ss$ consists of three families of arcs.  One family of 
$p$ arcs runs between the two copies of meridian $\mz$ in $\Ss$; two 
arcs each run from a copy of $\mz$ to a copy of $\mm$.  It is natural 
to parameterize slopes of proper arcs on $\Ss$ using these arcs of 
$K$.  Indeed, the discussion will now, in some sense, be parallel to 
that of \cite{ST1}.  We declare the family of $p$ arcs to have slope 
$0$ and the second pair to have slope $\infty$.  An outermost disk of 
$D \subset \eta(\theta)$ cut off by the pair of meridians $\mm, \mz$ 
defines a wave in $\Ss$; the wave assumption guarantees that such an 
outermost disk also intersects $\mpp$, so we conclude that the wave 
has finite slope $u/v$ in the coordinates just defined by the arcs $K 
\cap \Ss$.  Moreover $u$ is odd since a wave in $\Ss$ will be based at 
each copy of a single meridian (either $\mm$ or $\mz$) (see \cite{ST1} 
for details).  An argument will now show that either $e^{\bot} \cup 
e^{+}$ is unknotted or $u/v = \pm 1$.

The arc $\omega$ is disjoint from the wave.  Suppose to begin that 
$\omega$ intersects both meridians $\mm$ and $\mz$.  Then some arc 
component $\bbb$ of $\omega \cap \Ss$ has one end on a copy of each of 
$\mm$ and $\mz$ in $\bdd \Ss$.  Then the slope $r/s$ of $\bbb$ is odd 
and can't differ from $u/v$, the slope of the waves, since if it did 
its ends would have to run between the base of both waves, i.  e.  
different copies of the same meridian.  On the other hand, since 
$\bbb$ is disjoint from $K$, which has one parallel family of arcs of 
slope $0$ and two non-parallel arcs of slope $\infty$, we have $|r| 
\leq 2$ (hence $r = \pm 1$) and $|s| \leq 1$.  Since we are given that 
$u/v \neq \infty$ it follows that $u/v = r/s = \pm 1$ as claimed.

Next suppose that $\omega$ intersects $\mm$ but never $\mz$.  Then any 
component $\bbb$ of $\omega \cap \Ss$ that has both ends on copies of 
$\mm$ in $\Ss$ will have slope $0$ (since it's disjoint from $K$).  
The two terminal segments of $\omega$ in $\Ss$ will then each have one 
end on different copies of $\mm$.  But then they can't have their 
other end (i.  e.  the end points of $\omega$) on the same side of 
$K$.  For if they did, then either the arcs cross in the ``square'' 
component of $\Ss - (K \cup \bbb)$ in which they lie or one must be 
part of a segment of $\bdd D \cap \Ss$ of slope $\geq 1$ and the other 
of slope $\leq -1$.  See Figure \ref{fig:pillow}.  (If not all $p$ arcs 
of $K \cap \Ss$ that run between the two copies of $\mz$ in $\bdd \Ss$ 
are parallel, the slopes of both these arcs could be $\pm 1$, still 
sufficient to deduce that this is the slope of the wave.)

\begin{figure}
\centering
\includegraphics[width=.8\textwidth]{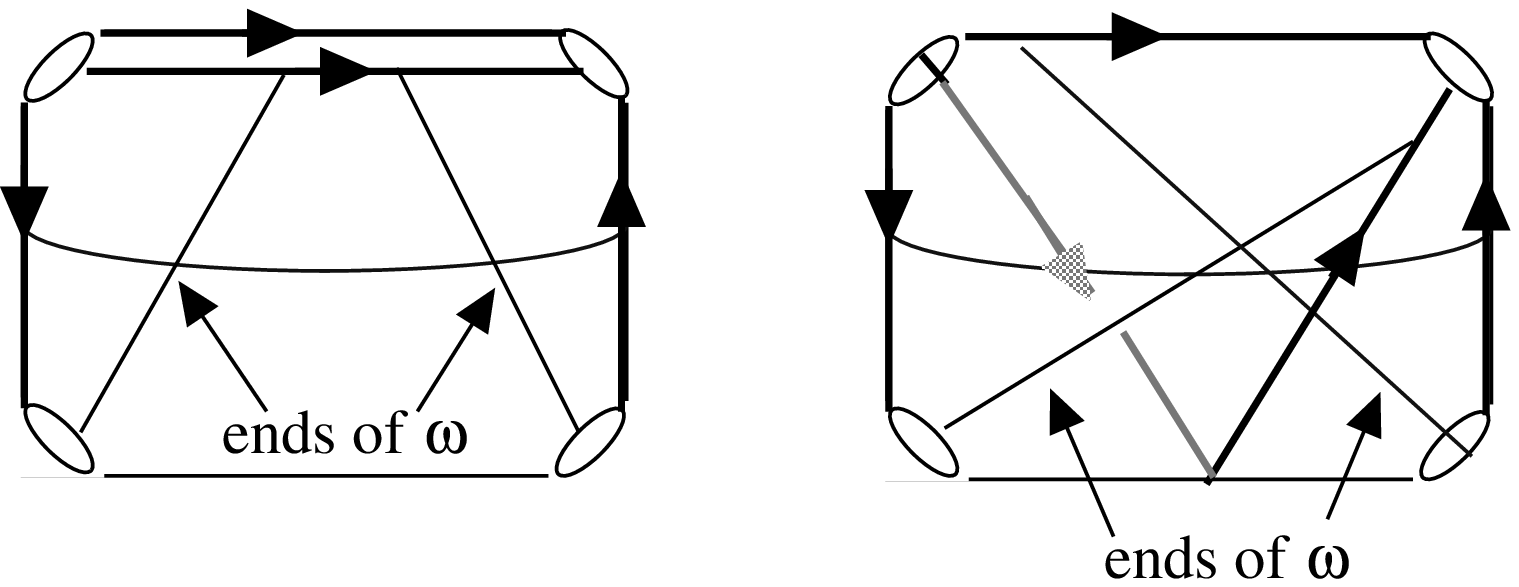}
\caption{} \label{fig:pillow}
\end{figure} 

The only remaining possibility (to avoid the conclusion that 
$u/v = 1$) would be that $\omega$ is disjoint from $\mm$.  But 
that would imply that the boundary of the annulus $A$ obtained by 
$\bdd$-compressing $F$ to $\bdd \eta(\theta)$ intersects $\mm$ in a 
single point.  Then by Lemma \ref{lemma:noshort} $e^{\bot} \cup 
e^{+}$ is unknotted.

So we continue, assuming that the slope of the waves in $\Ss$ is $\pm 
1$.  Now apply a Whitney move, replacing the meridian of $e^{+}$ 
(slope $\infty$) with the disk whose boundary has slope $u/v$.  This 
redefines the $\Theta$ curve as $\theta'$, presenting $K$ as a $(p+1, 
1)$ quasi-cable.  Moreover it is no thicker than $\theta$ and it 
satisfies the wave condition, since the meridian of the new edge 
$e^{\prime\bot}$ has been chosen to be disjoint from the wave.  As 
usual, we can ensure that the Whitney disk has no effect on the bridge 
structure of $K$, so $K$ remains in thin position.  Recall Figure 
\ref{fig:descend2}.
\end{proof}

\begin{prop} \label{prop:eminus}
Suppose in a thinnest $\Th$-graph $\theta$ appropriate for $(K, F)$ 
the edge $e^{-}$ is disjoint from a dividing sphere and suppose 
$genus(F) = 1$.  Then $p = 1$ and either $e^{\bot} \cup e^{-}$ is 
unknotted, or a Whitney move on $\theta$ changes it to an equally thin 
$\Th$-graph $\theta'$ that is appropriate for $(K, F)$.  Moreover 
$\theta'$ presents $K$ as a $(2, 1)$ quasi-cable.
\end{prop}

\begin{proof} Without loss we assume $e^{-}$ lies above a dividing 
sphere and $e^{-}$ is monotonic.  Suppose $p \geq 2$.  Then $K$ has a 
maximum at the lowest vertex.  Find a level sphere as in Lemma 
\ref{lemma:essential}.  The result contradicts Lemma 
\ref{lemma:hexagon1}.  Having established that $p = 1$, switch the 
labels of $e^{\pm}$ and apply Proposition \ref{prop:eplus2}.  \end{proof}

Propositions \ref{prop:eplus2} and \ref{prop:eminus} focus attention 
on the single remaining case to consider: when $e^{\bot}$ is disjoint 
from a dividing sphere.

\begin{theorem} \label{theorem:ebot}
Suppose in a thinnest $\Th$-graph $\theta$ appropriate for $(K, F)$, 
the edge $e^{\bot}$ is disjoint from a dividing sphere.  Suppose 
further 
that $genus(F) = 1$.

Then either the cycle $e^+ \cup e^{\bot}$ is unknotted or $p = 1$ and 
the cycle $e^- \cup e^{\bot}$ is unknotted.
\end{theorem}

\begin{proof} That $q = 1$ follows from Proposition 
\ref{prop:oneonethin2}.  As previously, let $A$ be the annulus 
obtained from $\bdd$-compressing $F$ to $\eta(\theta)$ using the disk 
$E$ from the splitting sphere.  Without loss of generality, assume 
$e^{\bot}$ lies above the dividing sphere.  If there are no regular 
maxima of $\theta$ then $e^{+} \cup e^{\bot}$ is vertical, and we are 
done.  If there is a regular maximum of $\theta$, we can assume it's 
the lowest maximum.  In that case, a level sphere just below the 
lowest maximum cuts off an upper disk from $A$ and a level sphere just 
above the highest minimum cuts off a lower disk from $A$.  So either 
some dividing sphere cuts off both an upper disk and a lower disk or 
some dividing sphere $P$ intersects $A$ only in essential arcs.  In 
the former case Corollary \ref{cor:reguplowb}
finishes the proof.  

The rest of the proof is an extended proof by contradiction.  We will 
show that it is impossible for a dividing sphere to intersect $A$ only 
in spanning arcs.   

Let $\omega \subset \bdd\eta(\theta) - K$ be the arc as previously, 
again slid to minimize intersections with the meridians of 
$\eta(\theta)$.  Let $P$ again be a level sphere between the highest 
minimum and the lowest maximum, and $\Ss$ again be the $4$-punctured 
sphere lying in $\bdd\eta(\theta) - P$, again supposing (with no loss) 
that both vertices are $\lll$ vertices, so $\Ss$ lies above $P$.  Let 
$w$ be the word in $a, b, \oa, \ob$ represented by $\omega$.  The wave 
condition now guarantees that (with the right choice of direction for 
$w$) $w$ is positive in $a$ and $\ob$ (including as usual the 
possibility that $w$ is the empty word).  If the wave is based at 
$\mm$ then no proper power of $\ob$ occurs in $w$ (i.  e.  no power 
greater than one); if it's based at $\mpp$ then no proper power of $a$ 
occurs in $w$.  Exploiting these facts, together with the rotational 
symmetry of the diagram, we have several essentially different ways in 
which the ends of $w$ can lie in $\Ss$.  Representative samples 
indicating that $w$ can begin or end on any letters are shown in 
Figure \ref{fig:wend6}.  We have oriented $w$ from left to right.  
Symmetric figures in which $w = \ob\ldots a$ are not shown.  The waves 
themselves are also not shown, but they are described (except for 
details of how their ends lie near $\mu^{\pm}$) by the requirement 
that they are disjoint from $\mu^{\bot}$, which is shown.  Note also, 
that the number of arcs of $K$ connecting the copies of $\mu^+$ is 
now $p - 1$.

The case $p = 1$ is special, since in this case there are no arcs 
connecting the copies of $\mu^+$.  Variants that arise in this 
case are shown in Figure \ref{fig:wend4}.

\begin{figure}
\centering
\includegraphics[width=1.0\textwidth]{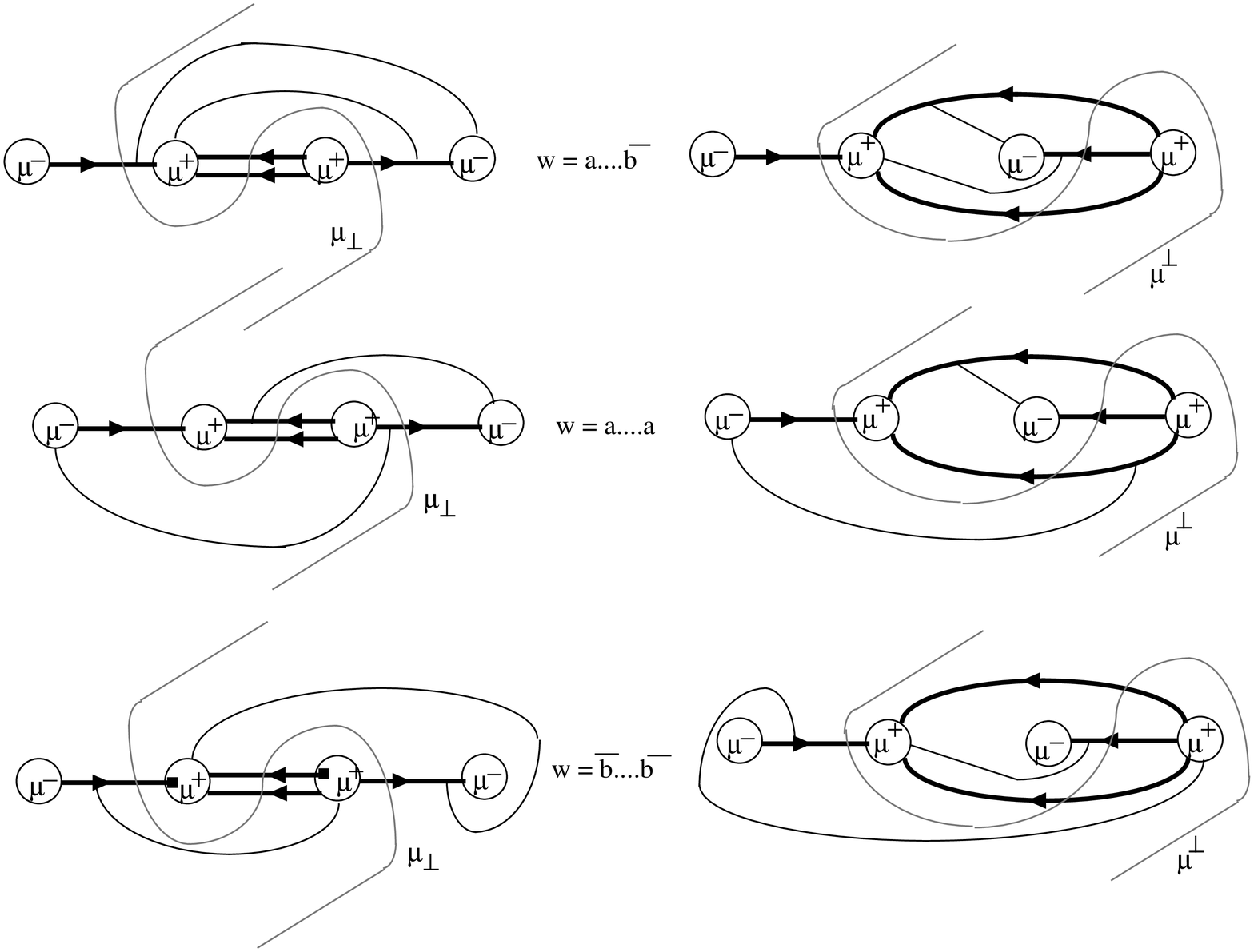}
\caption{} \label{fig:wend6}
\end{figure} 

\begin{figure}
\centering
\includegraphics[width=.8\textwidth]{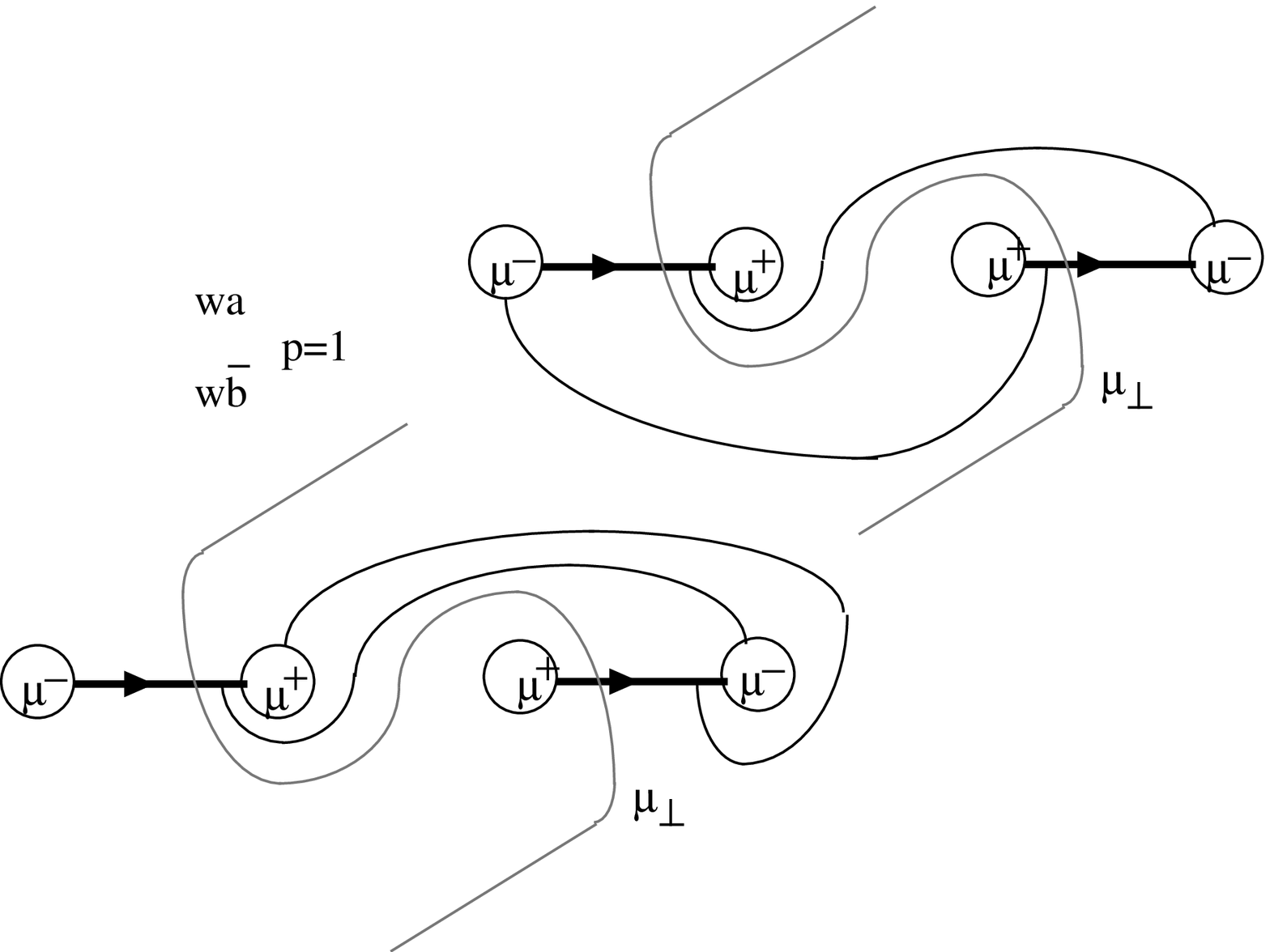}
\caption{} \label{fig:wend4}
\end{figure} 

The resulting words corresponding to $\bdd_-A$ and $\bdd_+A$ can be 
put in two forms:

\begin{enumerate}
	
\item \begin{itemize}
	
	\item $\bdd_+ A \leftrightarrow wb^{i}ab^{j}$
	
	\item $\bdd_- A \leftrightarrow w\ob^{k}$
	
\end{itemize}

\item \begin{itemize}
	
	\item $\bdd_+ A \leftrightarrow w\ob^{i}\oa\ob^{j}$
	
	\item $\bdd_- A \leftrightarrow wb^{k}$
	
\end{itemize}

\end{enumerate}

Here $i, j, k \geq 0, i+j+k = p.$

Now apply thin position.  It follows immediately that there is a level 
sphere $P$ so that no arc of $P \cap A$ cuts off an outermost disk 
incident to $\bdd A$, since there can't be both an upper and a lower 
such disk incident to $\bdd A$ simultaneously by Lemma 
\ref{lemma:reguplowb}.  In particular, $P \cap A$ consists only of 
spanning arcs, so $l(\bdd_-A) = l(\bdd_+A).$ This immediately implies 
that $(k - i - j)l(b) = l(a)$.  Thus it also forces $k > 0$.  In the 
second form above observe then that by 
regularity of $\bdd_{-} A$, $w$ begins and ends with the letter $a$ so by 
regularity of $\bdd_{+} A$, $i, j > 0$.

In fact we will show that the second category above 
does not arise and the first is limited to the case $i = j = 0, k = p$.  

\begin{lemma}  \label{lemma:simple}
The letters $b$ and $\oa$ do not occur in the words determined by 
$\bdd A$.
\end{lemma}

\begin{proof} Consider first the form
\begin{itemize}
	
	\item $\bdd_+ A \leftrightarrow wb^{i}ab^{j}$
	
	\item $\bdd_- A \leftrightarrow w\ob^{k}$.
	
\end{itemize}
We will show that $i = j = 0$ so $k = p$.  Suppose, with no loss of 
generality, that $i > 0$ so the length of the segment $b^{i}ab^{j}$ 
(the same as the length of the segment $\ob^{k}$) is at least $l(a) + 
l(b)$.  The easy case is when the wave is based at $\mm$ so there are 
no repeating $\ob$'s in $w$.  Then opposite (in $\bdd_+ A$) to the 
segment $\ob^{k}$ is a complete collection of $a$-labels (perhaps 
separated by a single letter $\ob$).  On the other hand, opposite (in 
$\bdd_- A$) to $b^{i}$ must be part of an $a$-segment.  Combining the 
two easily gives a Lens space contradiction.

It seems to be harder to establish a Lens space contradiction in the 
case when the wave is based at $\mpp$, so there are no repeating $a$'s 
in $w$.  For the first time we need to use the graph $G$ in the 
dividing sphere $P$ whose vertices are the points of intersection of 
$e^{\pm}$ with $P$ (we will call these points the $a$ and $b$-vertices 
in $G$) and whose edges are the arcs $P \cap A$, viewed both as 
spanning arcs of $A$ and as edges in $G$.  Ends of edges at the same 
vertex in $G$ will be said (usually in $\bdd A$) to have the same 
label.  The label will be an $a$-label, $b$-label, $\oa$-label or 
$\ob$ label depending on whether $\bdd A$ at that point (as oriented 
so that $w$ is positive in $a$ and $\ob$ as above) is passing through 
an $a$- or $b$-vertex in the direction that the $\theta$-graph is 
oriented or in the opposite direction.  Our use of the graph $G$ in 
this lemma will be modest, mostly as a book-keeping device.  Once the 
lemma is established we will need to examine $G$ much more seriously.

Orient all edges in $G$ to point from $\bdd_{+} A$ to $\bdd_{-} A$.  
We first claim there is an oriented path that begins with an edge 
incident to a $b$-label and ends with an edge incident to a 
$\ob$-label.  To see this, remove the $b$-vertices from $G$ (but not 
their incident edges), and let $G_{b}$ denote those edges, and the 
$a$-vertices they pass through, that are part of an oriented path 
beginning with an edge incident to a $b$-label.  If there are $m$ 
$a$-vertices in $G_b$ and there are $q$ occurences of $a$ in the word 
$w$ then $G_b$ has $(q+1)m + (i + j)l(b)$ ends of edges from $\bdd_{+} 
A$.  Yet only $qm$ ends of edges in $\bdd_{-} A$ could be in $G_b$ but 
not at $\ob$-labels.  So some ends of the edges must be at $\ob$-labels, 
as claimed.  Consider then the shortest oriented paths beginning with 
a $b$-label and ending with a $\ob$-label.  Among all 
shortest such paths, pick a path $\rrr$ whose ends are closest in 
$e^{+}$ as measured by the number of components of $e^{+} - P$ that 
lie between them.  We claim that that number is one; i.  e.  $\rrr$ 
begins and ends at the opposite ends of the same interval of $e^{+} - 
P$.

First note that the ends of $\rrr$ can't be at the same $b$-vertex, 
because $\rrr$ would then be a loop in $P$ whose normal $I$-bundle, as 
pieced together from neighborhoods of the edges in $A$, would not be 
oriented.  So call the initial $b$-vertex $b_{i}$ and the terminal 
$b$-vertex $b_t$.  To be concrete suppose that, in a single $b$-letter 
of $\bdd_{+} A$, $b_{t}$ precedes $b_{i}$ (we'll say that $b_{t}$ lies 
to the left of $b_{i}$ in the oriented $\bdd A$).  Unless the labels 
are precisely adjacent in $\bdd_{+} A$ (which is our claim), we can 
construct a better path $\rrr'$ as follows: Start at the $b$-label 
just to the left of the origin of $\rrr$ and construct a path by 
always using the edge that is one to the left (in $A$) of the edge in 
$\rrr$.

Notice first of all that the collection of edges $\rrr'$ we have just 
described is indeed a path in $G$: Suppose $\aaa_1$ and $\aaa_2$ are 
successive edges of $\rrr$ and the edges to their left in $A$ are 
$\aaa'_{1}$ and $\aaa'_2$.  We need to show that the end of 
$\aaa'_{1}$ in $\bdd_{-} A$ is at the same vertex as the end of 
$\aaa'_2$ in $\bdd_{+}A$.  (See Figure \ref{fig:path}.)  This is 
obvious unless the end of $\aaa_1$ at $\bdd_- A$ (and so the end of 
$\aaa_2$ in $\bdd_{+}A$) is the first label of an $a$-segment.  But if 
it were then, since there are no repeating $a$'s in $w$, the end of 
$\aaa'_{1}$ would in fact be a $\ob$ label and we would have found a 
shorter path.  Having established that $\rrr'$ is in fact a path in 
$G$, notice that it ends just to the left of the label $b_t$ in $\ob$ 
hence to the right of the label $b_t$ in $b$.  Hence we have found a 
path of equal length but with ends closer together.

\begin{figure}
\centering
\includegraphics[width=.5\textwidth]{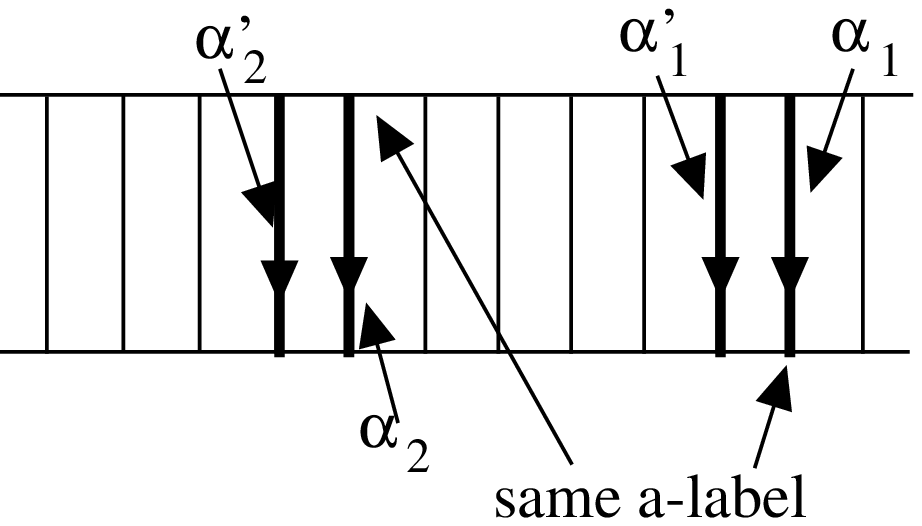}
\caption{} \label{fig:path}
\end{figure}

Having established that the ends of $\rrr$ represent adjacent 
intersections of $e^{+}$ with $P$, carry through the above 
construction of $\rrr'$.  Consider the sequence of squares that lie 
between the two paths.  When glued together along the components of 
$e^{+} - P$ that their edges in $\bdd A$ represent, the result is a M\"obius 
band whose boundary lies on the level sphere $P$.  This is impossible 
(and represents a different way of viewing a Lens space 
contradiction).

\bigskip

Next consider the form

\item \begin{itemize}
	
	\item $\bdd_+ A \leftrightarrow w\ob^{i}\oa\ob^{j}$
	
	\item $\bdd_- A \leftrightarrow wb^{k}$
	
\end{itemize}

\bigskip

Since $kl(b) = l(a) + (i +j)l(b) \geq l(a) + 2l(b)$ and, by the Lens 
space argument, no part of the $b^{k}$ interval can lie across from 
any part of an $\ob$-interval, it follows that no part of the $b^{k}$ 
interval can lie across from any part of the $\ob^{i}\oa\ob^{j}$ 
segment.  In particular, across from the $b^{k}$ interval must be an 
$a$-segment longer than $l(a)$.  This means that there are repeated 
$a$'s in $w$, hence no repeated $\ob$'s.  What then lies across from 
the $\oa$ interval?  None of it can by part of an $a$-segment, by the 
Lens space argument, so it must be part of a single $\ob$-interval.  
Then $l(b) \geq l(a)$.  On the other hand, $l(a) = (k - i - j)l(b) 
\geq l(b)$.  We conclude that $l(a) = l(b)$ and immediately across 
from $\oa$ is precisely a single $\ob$.  But if this is the case, then 
each segment corresponding to a letter in one boundary component will 
lie exactly opposite a segment corresponding to a single letter in the 
other boundary component.

This is clearly a very special case, and resolving it will begin the 
process of understanding how to use the graph $G$ effectively.  We've 
already identified (across from $b^{k}$) a term in $\bdd_{+}A$ of the 
form $a^{n}, n \geq 3$.  Across from that same term $a^{n}$ in $\bdd_- 
A$ must be three letters, at least one of which is also an $a$.  There 
are then two letters $a$ exactly aligned opposite each other, 
exhibiting that each $a$-vertex is the base of a loop in $G$.  Ask 
then what lies in $\bdd_{-} A$ exactly opposite $\ob^{i}$.  If any of 
it is an $a$ interval then this fact, together with the established 
fact that across from the $\oa$ interval is a $\ob$ interval, gives a 
Lens space contradiction.  If any of it is a $\ob$ interval then we 
will have exhibited that every $b$-vertex in $G$ is also the base of a 
loop.  Then, since every vertex is the base of a loop, some such loop 
will contain no vertices in its interior.  The following lemma shows 
that this is impossible.
\end{proof}

\begin{lemma}  \label{lemma:noloop}
Any loop in the graph $G$ must have vertices in both disks into which the 
loop divides $P$.
\end{lemma}

\begin{proof} A loop without such vertices in its interior would give 
a problematic $\bdd$-compression of $A$.  To see the problem, consider 
the base of an innermost such loop, that is, the subarc $\sss$ of the 
meridian to which $A$ is $\bdd$-compressed, an arc in 
$\bdd(\eta(\theta)) - K$.  The arc $\sss$ could not be a simple cocore 
of the band along $\omega$, otherwise $F$ would have been 
compressible.  On the other hand, if $\sss$ were incident to a copy of 
$\omega$ in $A$ from the side in $\bdd\eta(\theta)$ opposite to the band 
$\omega$, then the $\bdd$-compression would turn $A$ into an essential 
disk in $S^3 - \eta(\theta)$ whose boundary is disjoint from $K$, 
which is also absurd, for if we attach a neighborhood of the disk to 
$\eta(\theta)$, $K$ would lie on the resulting unknotted torus and $F$ 
would be a Seifert surface in the solid torus complement, forcing $K$ 
to be trivial.  The only remaining possibility is that $\sss$ has both 
ends incident to the parts of $\bdd A$ that come from $K$.  But $K$ 
crosses each meridian always in the same direction, so two such 
crossings can't be the ends of a $\bdd$-compression. 
\end{proof}

Lemma \ref{lemma:noloop} completes the proof of Lemma 
\ref{lemma:simple}.  We conclude that the words corresponding to the 
boundary components of the annulus $A$ are exactly

\bigskip

\begin{itemize}
	
	\item $\bdd_- A \leftrightarrow wa$
	
	\item $\bdd_+ A \leftrightarrow w\ob^{p}$.
	
\end{itemize}

\bigskip

In fact more can be said.  The fact that there is no occurence of the 
letter $b$ in $\bdd_- A$ or $\bdd_{+} A$ means that neither end of $w$ 
can lie on the segments of $K\cap \Ss$ connecting the two copies of $\mpp$ 
(when $p \geq 2$).  Then 

\begin{cor}  \label{cor:endb}
Suppose $p \geq 2$.  If the $p-1$ segments $\{ \kkk_{1}, \ldots, 
\kkk_{p-1} \}$ of $K \cap \Ss$ that connect 
the two copies of $\mpp$ are not all 
parallel in $\Ss$ then $w$ begins and ends with the letter $\ob$.  If 
they are all parallel then $w$ begins or ends (perhaps both) with the 
letter $\ob$.
\end{cor}

\begin{proof} If some $\kkk_{i}$ were incident to both ends of 
$\omega$ then one component of $\bdd A$ would represent the word $w$ 
and the other one $w\ob^{i}a\ob^{j}$.  This is impossible since these 
words have different lengths.  If some $\kkk_i$ contained a single end 
of $\omega$ then one component of $\bdd A$ would contain an occurence 
of the letter $b$, contradicting Lemma \ref{lemma:simple}.  So the 
ends of $w$ lie, one each, on the two components of $K \cap \Ss$ that 
are not among the $\kkk_{i}$.  The result follows easily (see Figure 
\ref{fig:wend6}).
\end{proof}

The argument now proceeds by considering every possible type of word 
$w$.  We begin by considering short words, then long words, then words 
of intermediate length.

\begin{lemma} \label{lemma:noshort2}
If $w = \ob^m, m \geq 0$ (e.  g.  $w = \emptyset$) then 
$e^{\bot} \cup e^+$ is uknotted.  

If $w = a^{m}, m \geq 1$ then $p = 1$ and 
$e^{\bot} \cup e^-$ is unknotted.
\end{lemma}

\begin{proof} If $w = \ob^{m}$ (or is empty) then $\bdd A$ intersects 
the meridian $\mm$ in exactly one point, a point in $\bdd_- A$.  If $w 
= a^{m}, m \geq 1$ then $p = 1$ by Corollary \ref{cor:endb}.  Then 
$\bdd A$ intersects the meridian $\mpp$ in exactly one point, a point 
in $\bdd_+ A$.  In both cases the result follows from Lemma 
\ref{lemma:noshort}.
\end{proof}

In view of Lemma \ref{lemma:noshort2} we can and will restrict our 
attention only to words that contain both letters $a$ and $\ob$.

To deal with longer words it will be useful to generalize Lemma 
\ref{lemma:noloop}.  To appreciate how, we examine the local structure 
of $G$.  The key to organizing the information is to orient each edge 
of $G$, as was done briefly above, so that the edge, when viewed in 
$A$, points from $\bdd_+ A$ to $\bdd_- A$.  This has the obvious 
consequence that any $a$-vertex has at least one edge pointing into 
it, since the word $wa$ contains the letter $a$, and each $b$-vertex 
has at least $p$ edges pointing out, since the word $w\ob^p$ contains 
at least $p$ occurences of the letter $\ob$.  Beyond these ends of 
edges, though, is one more pair at each $a$-vertex (resp.  $b$-vertex) 
for each occurence of $a$ (resp.  $\ob$) in $w$.  (One of the pair is 
identified with the occurence of the letter in $w \subset \bdd_{+} A$ 
and the other with the occurence in $w \subset \bdd_{-}A$.)  Call 
these ends of edges the $w$-ends.  At any $a$-vertex there is a single 
non-$w$ edge pointing in plus a sequence of $w$-ends alternating 
between pointing in and pointing out.  At any $b$-vertex there are $p$ 
non-$w$ edge pointing out and a collection of $w$-ends, 
the latter coming in pairs of adjacent ends, one pointing in and one 
pointing out.  See Figure \ref{fig:valence}.

\begin{figure}
\centering
\includegraphics[width=1.0\textwidth]{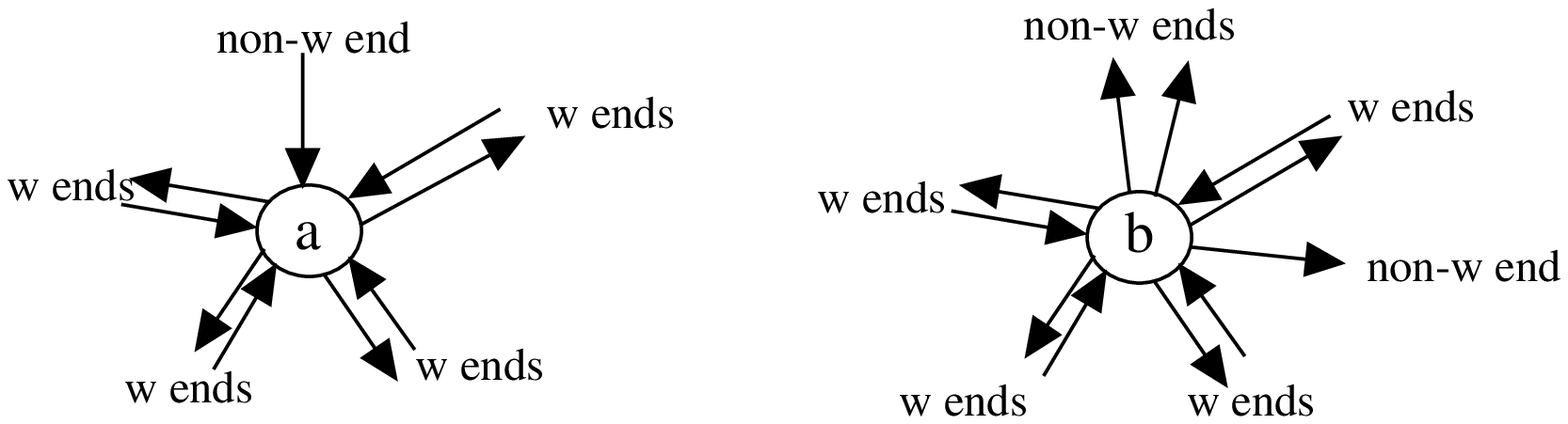}
\caption{} \label{fig:valence}
\end{figure}

Between each pair of $w$-ends in $\bdd \mm$ or $\bdd \mpp$ is an arc 
$\omega_{\bot}$ that is a cocore of the band along $\omega$.  Put 
another way, banding together the boundary components of $\bdd A$ along 
$\omega_{\bot}$ would recover $F$ from $A$.  Call that side the 
$w$-side of the $w$-end.  Since occurences of $a$ or $\ob$ in $w$ 
occur always with the same sign, $w$ crosses $\mm$ always in the same 
direction, and similarly for $\mpp$.  It follows that, at any given 
vertex of $G$, the $w$-side of any (oriented) $w$-end is always to the 
right (or always to the left) of the end, as the end is oriented by 
the orientation of its edge.  Moreover, since the ends of 
$\omega_{\bot}$ are incident to the same side of $A$, that normal 
direction is well-defined, so if both ends of the same edge in $G$ are 
$w$-ends, then the $w$-side is the same side at both ends.  Combining 
these facts we discover that, throughout any one component $G_0$ of 
$G$, the $w$-side of any $w$-end at any vertex lies always to the same 
side (say always to the right as the edge is oriented) of the $w$-end.  
Under these circumstances, note that a cycle in $G_0$ that has no 
vertices or edges in its interior, and which moves clockwise around 
its interior (we call it a {\em clockwise face}), must have corners 
that are always on the $w$-side of $w$-edges.  In particular, if such 
a cycle can be found, then its interior would, before the 
$\bdd$-compression that changed $F$ to $A$, correspond to a 
compressing disk $E$ for $F$.  (We know that $\bdd E$ would be 
essential in $F$, since it crosses a proper arc in $F$, namely the one 
$\bdd$-compressed to $\omega$, always in the same direction.)  Such a 
compression, of course, violates our assumption that $F$ is 
incompressible.

Symmetrically, in a component of $G$ for which the $w$-side of any 
$w$-end is to the left of the oriented end, there could be no 
counterclockwise cycle whose interior is empty (i.  e.  no 
counterclockwise face).

Although it might not be easy to see if a given component of the graph 
$G$ in $P$ is ``right-handed'' or ``left-handed'' in this sense, it is 
possible to use the extreme regularity of $G$ (guaranteed by the fact 
that all edges in $G$ are parallel in $A$) to identify circumstances 
in which there are {\em both} clockwise and counterclockwise cycles in 
the same component of $G$ with no vertices in the interior of either.  
That, then, forces one or the other to define a compression of $F$, a 
contradiction.  For example, we have

\begin{lemma}  \label{lemma:nofacepair}
No two faces of $G$ can be adjacent and have boundaries that are 
cycles.
\end{lemma}

\begin{proof}
Since the cycles are adjacent, one is clockwise and one is 
counterclockwise.
\end{proof}

\begin{defin}  \label{defin:face}
 A disk component of $P - G$ will be called a {\em face}.  If the 
 boundary of the face is a cycle we call it a {\em face cycle}.  A 
 clockwise (resp.  counterclockwise) face cycle will be called a {\em 
 clockwise (resp.  counterclockwise) face}.  A face incident only to 
 $a$-vertices (resp.  $b$-vertices) will be called an {\em $a$-face} 
 (resp.  $b$-face).
\end{defin}

\begin{lemma}  \label{lemma:afacecycle}
Any $a$-face is a face cycle.  If $p = 1$ any $b$-face is a face cycle.
\end{lemma}

\begin{proof}  At any $a$-vertex exactly one end of an edge is not a 
$w$-end and it points into the vertex.  Hence there cannot be two 
adjacent ends of edges pointing out, as there would be in an $a$-face 
that is not a cycle.
\end{proof}

\begin{lemma} \label{lemma:noadja}
No two $a$-faces can be adjacent.  If $p = 1$ no two $b$-faces can be 
adjacent.
\end{lemma}

\begin{proof}  Combine lemmas \ref{lemma:nofacepair} and \ref{lemma:afacecycle}
\end{proof}

It is natural to seek features of the graph which guarantee the 
existence of cycles. The following lemma suggests a possibility.

\begin{lemma} \label{lemma:ebigon}
 No distinct 
$b$-vertices can have edges pointing to (resp. from) the same vertex.
\end{lemma}

\begin{proof} Suppose $\aaa_{1}, \aaa_{2}$ are two edges in $G$ with their 
heads, say, at the same vertex of $G$.  Then the ends of the $\aaa_{i}$ in 
$\bdd_{-} A$ are some multiple of $l(b)$ apart.  (Recall that $l(a) = 
pl(b)$.)  The ends of the $\aaa_{i}$ in $\bdd_{+} A$ are the same 
distance apart.  So if both ends on $\bdd_{+} A$ are at $b$-labels, 
they must be at the same $b$-label.
\end{proof}

It seems from this lemma that bigons may be prevalent.  To be 
precise, define a {\em bigon} in $G$ to be a pair of 
edges, each running between the same pair of vertices.  A {\em 
parallel} bigon will be a bigon in which both edges of the bigon are 
oriented in the same direction.  An {\em anti-parallel} bigon will be 
one in which the edges are oriented in the opposite direction, forming 
a cycle in $G$ of length two.  In the case of a parallel bigon we will 
denote the vertex from which the edges of the bigon point out as $v_+$ 
and the vertex into which the edges point as $v_{-}$.

\begin{lemma}  \label{lemma:nobigonc} Suppose there is a 
parallel bigon in $G$ and let $B$ be a disk in $P$ that it bounds.  
Suppose in $B$ there is an oriented path from $v_{-}$ to $v_{+}$.  
Then in the interior of $B$ there is a cycle that is disjoint from 
both $v_{\pm}$.
\end{lemma}

\begin{proof} With no loss we may assume that $B$ contains no other 
parallel bigon, else we would focus on an innermost one.  Since $w$ 
contains both letters $a$ and $\ob$, any vertex is incident to an edge 
pointing out and an edge pointing in; no vertex is a sink or source.  
Hence any component of $G$ contains a cycle, so we may as well assume 
that every vertex in the disk belongs to the same component $G_0$ of 
$G$ as the bigon.  Suppose, with no loss (as explained above), $G_0$ 
contains no clockwise face.

If there were an oriented path that runs from $v_-$ to $v_+$ inside 
$B$, then the closed disk would contain both a clockwise and a 
counterclockwise cycle.  Consider an innermost clockwise cycle 
(perhaps passing more than once through the same vertex, but not 
crossing at such a vertex) and the disk $B'$ that it bounds.  Suppose 
$B'$ contains a vertex.  That vertex must be part of an oriented path 
in the interior of $B'$.  If that path forms a cycle completely in the 
interior of $B'$ we are done.  If not (e.  g.  the ends of the path 
are at the same vertex of the cycle $\bdd B'$) the path would cut off 
a clockwise cycle that would be even further in, a contradiction.  
(See Figure \ref{fig:ccycle}.)  So there is no vertex in the interior of $B'$.  
Similarly, if there were an edge in the interior of $B'$ there would 
be a further in clockwise cycle.  We conclude that $B'$ would have to 
be a clockwise face, which is impossible.  
\end{proof}

\begin{figure}
\centering
\includegraphics[width=.6\textwidth]{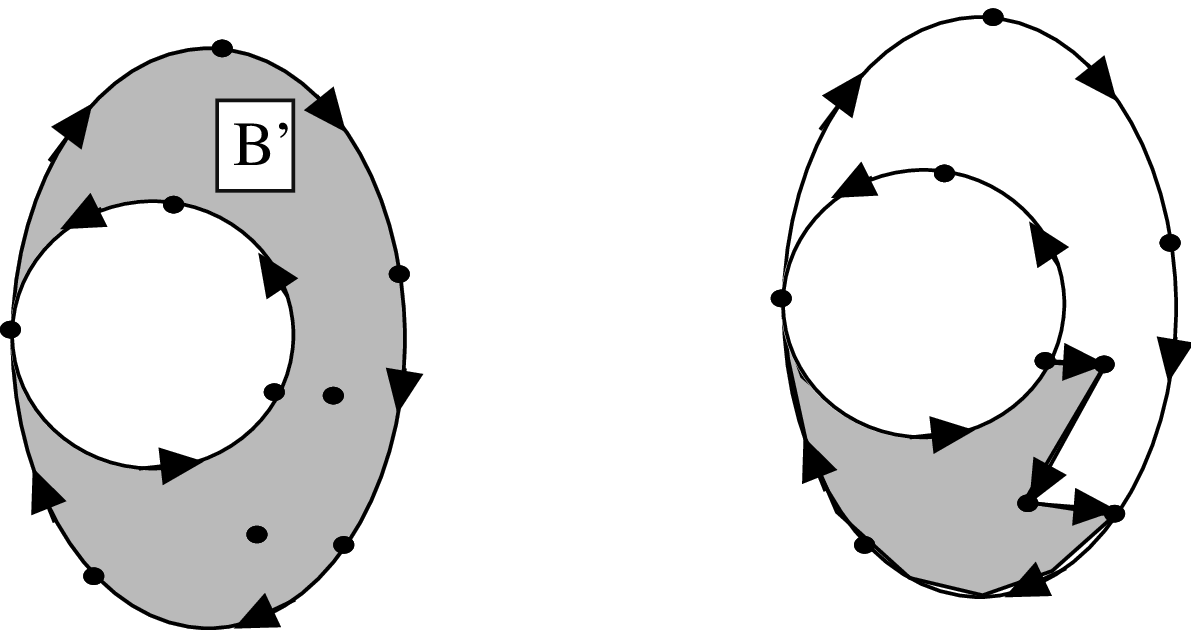}
\caption{} \label{fig:ccycle}
\end{figure}

\begin{lemma} \label{lemma:nobigon} Suppose there is a 
parallel bigon in $G$ and let $B$ be a disk 
in $P$ that it bounds.  Then in the closure of $B$ (i.  e.  including 
the vertices $v_{\pm}$) there is a cycle that includes at most one of 
$v_{\pm}$.
\end{lemma}

\begin{proof} As in the proof of Lemma \ref{lemma:nobigonc} we can 
assume, with no loss, that $B$ contains no other parallel bigon, no 
vertex is a sink or source, that every vertex in the disk belongs to 
the same component $G_0$ of $G$ as the bigon, and that $G_0$ contains 
no clockwise face.

If any end of $G$ in $B$ points out from $v_-$, then it is part of an 
oriented path (since no vertex is a source or sink).  If the path ends 
in $v_{+}$ we are done by Lemma \ref{lemma:nobigonc}.  If not, it must 
contain a cycle not incident to $v_+$, and we are done.

If the only ends of $G$ in the bigon that are incident to $v_-$ point 
into $v_-$ then there can be no edges other than the bigon itself, 
since even an $a$-vertex can have at most two adjacent ends pointing 
in.  (Recall that $w$-ends alternate between pointing in and pointing 
out.)  So we may as well assume that no end of $G$ in the interior of 
$B$ is incident to $v_-$.  In that case, any end lying in the the 
interior of $B$ and incident to $v_+$ must be part of a cycle in the 
bigon incident only to $v_+$ and we are done.

The possibility remains that the interior of the bigon is empty.  In 
that case, at least one end of the edges of the bigon at each vertex 
is not a $w$-end (since $w$-ends alternate between pointing in and 
pointing out) so $v_-$ is an $a$-vertex and $v_+$ is a $b$-vertex.  If 
neither end at $v_+$ is a $w$-end then, considering how non-$w$-ends 
arise, necessarily $p \geq 2$ and the edges of the bigon, when viewed 
in $A$, are some $kl(a)/p, 1 \leq k < p$ apart.  But then they can't 
have their other end at the same $a$-vertex, since two ends in 
$\bdd_{-} A$ with the same $a$-label are at least $l(a)$ apart in $A$.  
Hence we conclude that exactly one end of the bigon at each of $v_+$ 
and $v_-$ is a $w$-end.  Necessarily their $w$-side is the same side 
and not the side in the bigon.  Hence exactly one edge $\aaa_l$ in the 
bigon (the left one, say) has both of its ends $w$-ends and {\em 
neither} of the ends of the other edge $\aaa_r$ are $w$-ends.  In 
particular, $\aaa_r$ connects, in the annulus $A$, a point in $\bdd_- 
A$ corresponding to a point in the last letter of $wa$, to a point in 
$\bdd_+ A$ that lies in the final syllable $\ob^p$ of $w\ob^p$.

\begin{figure}
\centering
\includegraphics[width=1.0\textwidth]{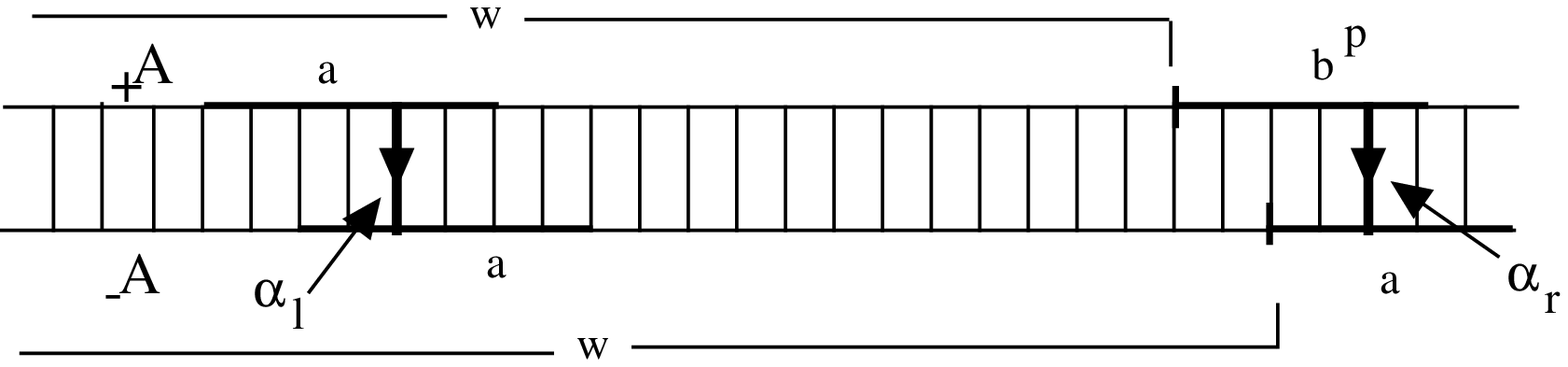}
\caption{} \label{fig:fence}
\end{figure}

Now $l(a) = pl(b)$; suppose the end of $\aaa_r$ in $\bdd_{+} A$ is the 
$i$th end in the final syllable $\ob^{p}$ of $w\ob^{p}$ and the end in 
$\bdd_{-}A$ is the $j$th end in the final letter $a$ of $wa$.  Suppose 
for concreteness that $i \geq j$.  Now consider $\aaa_{l}$.  Since 
it's ends are the same as those of $\aaa_r$, the end of $\aaa_{l}$ in 
$\bdd_{-} A$ is the $j$th end in some letter $a$ in $w$.  But the 
length of the terminal segment of $w$ is of course the same in 
$\bdd_{+} A$ as in $\bdd_{-}A$ and so the end of $\aaa_{l}$ in 
$\bdd_{+} A$ is the point exactly $i - j$ later than the end is in 
$\bdd_{-} A$.  That is, it is still part of the same letter $a$ of 
$w$.  (See Figure \ref{fig:fence}.)  This is a contradiction, since 
the other end of $\aaa_l$ is at a $b$-vertex.  If $i \leq j$ we get 
the same contradiction, using the initial segment of $w$ instead of 
the terminal segment.
\end{proof}

With a little determination, we have more:

\begin{lemma} \label{lemma:nobigona} The interior of any 
parallel bigon contains a $b$-vertex.

\end{lemma}

\begin{proof} With no loss we can assume that the parallel bigon is 
innermost.  Suppose the disk $B$ contained only $a$-vertices in its 
interior.  By Lemma \ref{lemma:nobigon}, $B$ contains a cycle that 
passes through some vertices in its interior.  Since the valence of 
each vertex is greater than two, there are adjacent faces in $B$.  
Then Lemma \ref{lemma:noadja} shows that $v_{\pm}$ can't both be 
$a$-vertices.  Similarly, all vertices in $B$ are in the same 
component of $G$ as the bigon.

Let $\Ggg$ denote the subgraph of $G$, lying in the interior of $B$, 
obtained by deleting all edges incident to $v_{\pm}$. 

{\bf Claim:} $\Ggg$ contains a cycle.  

Otherwise consider an oriented path from a source vertex $a_+$ to a 
sink vertex $a_{-}$.  All edges (of $G$) pointing into $a_+$ must have 
their other ends at $v_{\pm}$ and not the same vertex, since this 
would exhibit a further in parallel bigon.  There must be at least two 
such edges, since by Lemma \ref{lemma:noshort2}, $w$ contains at least 
one occurence of $a$ so the word $wa$ contains at least two.  This 
implies that there are exactly two edges pointing into $a_{+}$ and one 
edge comes from each of $v_{\pm}$.  Then Lemma \ref{lemma:ebigon} 
implies that $v_{\pm}$ can't both be $b$-vertices.  Suppose that 
$v_{+}$ were an $a$-vertex.  We have already identified three edges 
pointing out from $v_{+}$: the edges of the bigon and an edge to 
$a_{+}$.  Then three edges point out from $a_{-}$ and at least two 
would have to go to the same vertex in the pair $v_{\pm}$, creating a 
further in bigon.  So $v_+$ is a $b$-vertex and $v_-$ is an 
$a$-vertex.  Now consider $a_-$.  If an edge in $G$ pointing out from 
$a_{-}$ goes to $v_-$, then $v_{-}$, hence every $a$-vertex, has three 
edges pointing in.  This would force $a_{+}$ to be part of a further 
in bigon, a contradiction.  So the edge pointing out from $a_{-}$ goes 
to $v_+$ and this edge, together with the path from $a_+$ to $a_{-}$ 
together with the edge from $v_{-}$ pointing into $a_+$ give an 
oriented path from $v_{-}$ to $v_{+}$.  Then Lemma 
\ref{lemma:nobigonc} provides a cycle in $\Ggg$, as claimed. 

Having established the claim we continue with the proof of Lemma 
\ref{lemma:nobigona}. Note that, once we 
have a cycle in $\Ggg$ we know that $\Ggg$ contains a face hence, 
following Lemma \ref{lemma:afacecycle}, an $a$-face cycle.

\medskip

{\bf Case 1}:  $w$ contains only one occurence of 
$a$.  

In this case each $a$-vertex has valence three -- two ends of edges 
pointing in and one pointing out.  (So, for example, $v_{+}$ must be a 
$b$-vertex.)  We have seen that $\Ggg$ contains an $a$-face cycle 
$\sss$.  Now notice that there are only two possible ``corners'' of 
cycles at each $a$-vertex, since the valence is three; one set of 
corners occurs only in clockwise cycles and one in counterclockwise 
cycles.  As a result, all the corners of $\sss$ come from a single 
occurence of the letter $a$ in $\bdd_- A$ and a single occurence of 
the letter $a$ in $\bdd_{+} A$.  We claim this is impossible: Of all 
the vertices of $\sss$, let $a_{l}$ be the one that is first 
encountered when passing along the oriented edge $e^{-}$, so, in any 
occurence of the letter $a$ in $\bdd A$, the label corresponding to 
$a_{l}$ lies most to the left among all labels coming from vertices of 
$\sss$.  Similarly, define $a_r$ to be the last vertex of $\sss$ that 
is encountered along the oriented edge $e^{-}$.  We repeat: the edges 
pointing out from $a_{l}$ and from $a_{r}$, when viewed in $A$, leave 
from the same $a$-interval $a_{+}$in $\bdd_{+} A$ and end in the same 
$a$-interval $a_{-}$ in $\bdd_{-} A$.  But, by definition of $a_{r}$ 
and $a_{l}$, the edge pointing out from $a_{r}$ in $a_{+} \subset 
bdd_{+} A$, goes to a label to the left of $a_r$ in $a_{-} \subset 
\bdd_{-}A$ whereas the edge pointing out from the label $a_l$ in 
$a_{+}$ goes to a label to the right of $a_{l}$ in $a_{-}$.  This 
presents a clear contradiction: there are more ends of edges between 
the two edges in $\bdd_{+} A$ than there are in $\bdd_{-}A$.  See 
Figure \ref{fig:acycle}.

\begin{figure}
\centering
\includegraphics[width=.6\textwidth]{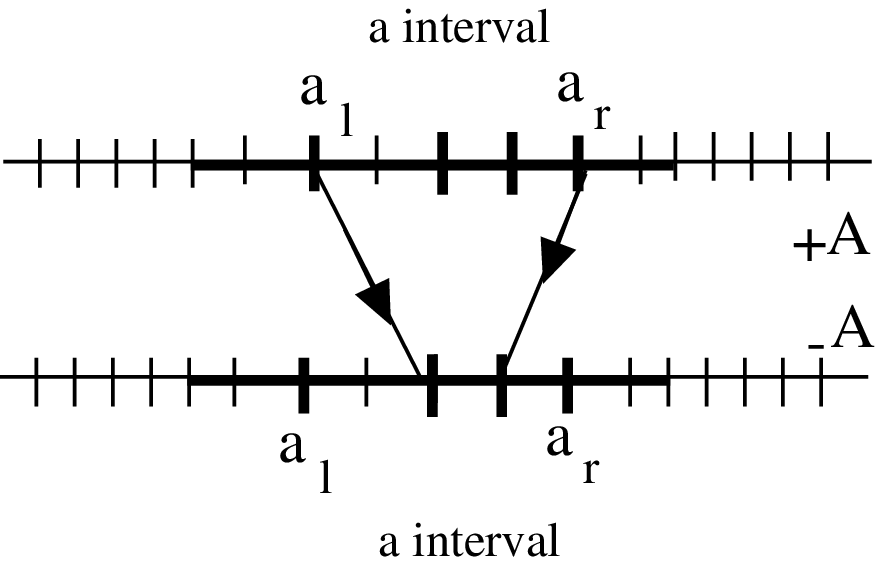}
\caption{} \label{fig:acycle}
\end{figure}

\medskip

{\bf Case 2}:  $w$ contains three or more occurences of $a$.

In this case we will show that there are two adjacent $a$-faces in 
$\Ggg$, contradicting Lemma \ref{lemma:noadja}.  Let $\Ggg_{0}$ be a 
component of $\Ggg$ and let $k$ denote the number of vertices in 
$\Ggg_{0}$.  If $v_{\pm}$ are both $b$-vertices then at each 
$a$-vertex in $B$ there can be at most two edges also incident to one 
of the $v_{\pm}$, one pointing in and one pointing out, since two 
different $b$-vertices can't have edges pointing toward (or away from) 
the same $a$-vertex (Lemma \ref{lemma:ebigon}) and if one $b$-vertex 
had two edges pointing toward (or away from) the same $a$-vertex it 
would be a further in bigon.  On the other hand, if either vertex (say 
$v_-$) were an $a$-vertex then at most two edges could be incident to 
both $\Ggg_{0}$ and $v_-$, for otherwise $G$ would have adjacent 
$a$-faces.  So clearly $k \geq 2$ (the valence of each vertex is at 
least $7$) and at most $2k+2 \leq 3k$ edges connect $\Ggg_{0}$ to 
$v_{\pm}$ in $G$.  It follows that $\Ggg_{0}$ contains at least $7k/2 
- 3k/2 = 2k$ edges.  The proof that this would provide two adjacent 
$a$-faces now follows from this simple observation:

\medskip

{\bf Claim:} If $\Ggg$ is a connected graph in the plane with $k$ 
vertices and at least $2k-1$ edges then either $\Ggg$ contains a 
trivial loop or two faces of $\Ggg$ (i. e. compact complementary 
components, necessarily disks) are adjacent.

Here's the proof of the claim.  Let $e \geq 2k-1$ be the number of 
edges.  If any face is a monogon, we are done.  So suppose every face 
has at least two edges.  Either some edge is incident to two faces, 
and we are done, or the number of faces $f \leq e/2$.  Then consider 
the Euler characteristic of $\Ggg$ and all its faces: $1 = v - e + f 
\leq v - e/2 \leq k - (k - 1/2) = 1/2$, a contradiction.

\medskip

{\bf Case 3}:  $w$ contains exactly two occurences of $a$.

Suppose first that $v_{\pm}$ are both $b$-vertices.  As usual, for 
each of the vertices in the interior of $B$, there is at most one edge 
pointing from $v_{\pm}$ to the interior vertex by Lemma 
\ref{lemma:ebigon} and our assumption that $B$ is an innermost 
bigon.  Now consider a component $\Ggg_{0}$ of $\Ggg$ with, say, $k$ 
vertices.  
We have just shown that there are at least $2$ edges in $\Ggg$ 
pointing into each $a$-vertex (for at most one edge pointing into the 
$a$-vertex comes from a $b$-vertex).  Hence there are at least $2k$ 
edges in $\Ggg_{0}$.  The proof now follows as in the previous case.

Suppose finally that one of $v_{\pm}$ is a $b$-vertex $v_{b}$ and one 
is an $a$-vertex $v_{a}$.  If, for some component $\Ggg_{0}$ of $\Ggg$ 
there are no edges (in $G$) running from $\Ggg_{0}$ to the $b$-vertex, 
we are done much as before.  Similarly, if there is at most one edge 
$\aaa$ in $G$ pointing from $v_{a}$ to $\Ggg_{0}$ we are done: In 
$\Ggg_{0}$ there are still at least $2$ edges pointing into every 
vertex, except for the single vertex at the end of $\aaa$.  Hence 
there are at least $2k-1$ edges in $\Ggg_{0}$ and we can still apply the 
combinatorial claim above.  If there are at least two edges, say 
$\aaa_1$ and $\aaa_2$, in $G$ pointing from $v_a$ to $\Ggg_{0}$ then, 
since no adjacent ends at an $a$-vertex point out, there is another 
end of an edge $\aaa_3$ between the ends of the $\aaa_i$ at $v_a$.  If 
$\aaa_3$ also goes from $v_a$ to $\Ggg_{0}$ then on either side of it 
are adjacent $a$-faces, contradicting Lemma \ref{lemma:noadja}.  If 
instead it goes to another component of $\Ggg$ then that component is 
cut off from $v_b$ by $\Ggg \cup \aaa_1 \cup \aaa_2$ so no edge in $G$ 
connects it to $v_{b}$, a case we have already established.
\end{proof}

Lemma \ref{lemma:nobigon} immediately eliminates the possibility 
that $w$ is a long word.  Explicitly, we have:

\begin{lemma} \label{lemma:oneb} The letter $\ob$ occurs at most 
once in $w$.
\end{lemma}

\begin{proof}
Suppose $\aaa$ and $\aaa'$ are two edges in $G$ that point out from 
the same $b$-vertex.  Then (echoing the argument of Lemma 
\ref{lemma:ebigon}) the distance between $\aaa$ and $\aaa'$ as 
measured along either of $\bdd_{\pm} A$ is some multiple of $l(b) = 
l(a)/p$.  In particular, there are at most $p$ candidates for 
$a$-vertices the other ends of $\aaa$ and $\aaa'$ might be incident 
to, plus a $b$-vertex.  If the letter $\ob$ occurs more than once, 
then in $\bdd_{+}A$ there are at least $p+2$ occurences of each 
$b$-label.  And, for each $b$-label, we have just argued that there 
are at most $p+1$ possible labels in $\bdd_{-}A$ to which they can 
point, $p$ of them $a$-labels and one a $b$ label.  Hence at least two 
of the edges point to the same label.  This shows that every 
$b$-vertex is part of a parallel bigon.  An innermost parallel bigon 
then would have to contain only $a$-vertices, contradicting Lemma 
\ref{lemma:nobigona}.
\end{proof} 

\begin{lemma} \label{lemma:onea} If $\ob$ occurs in $w$ then the 
letter $a$ occurs at most once.
\end{lemma}

\begin{proof} Following Lemmas \ref{lemma:noshort2} and 
\ref{lemma:oneb} we can restrict to the case in which $\ob$ occurs 
exactly once in $w$.  We will assume that $a$ occurs $m \geq 2$ times 
and derive a contradiction.  The structure of the proof depends on 
whether $p = 1$ or $p \geq 2$.

If $p = 1$ then any two occurences of the same $a$-label or the same 
$b$-label in $\bdd_{\pm} A$ occur a multiple of $l(a) = l(b)$ apart.  
It follows that at least two edges pointing into any given $a$-vertex 
in $G$ have their other ends at the same vertex.  In other words, each 
$a$-vertex is contained in some parallel bigon.  Thus an innermost 
parallel bigon contains only $b$-vertices.  Any $b$-vertex has at 
least one edge pointing out that goes to an $a$-vertex (since there is 
only one occurence of the letter $\ob$ in $\bdd_{-}A$ but two in 
$\bdd_{+} A$) and edges that point from $b$-vertices to a given 
$a$-vertex must all come from the same $b$-vertex by Lemma 
\ref{lemma:ebigon}.  It follows that there are at most two 
$b$-vertices in the the interior of the bigon.  If there were only 
one, then the two edges coming out from it can't go to the same 
vertex, for that would form another parallel bigon, so one goes to 
each of the vertices forming the bigon.  This would force $v_+$ to be 
a $b$-vertex and $v_-$ to be an $a$-vertex.  Furthermore, $b$-vertices 
are of valence $3$, and the edge pointing into the interior $b$-vertex 
has nowhere to come from but $v_-$.  The two adjacent face cycles 
contradict Lemma \ref{lemma:nofacepair}.  See Figure \ref{fig:onea}i.

\begin{figure}
\centering
\includegraphics[width=.5\textwidth]{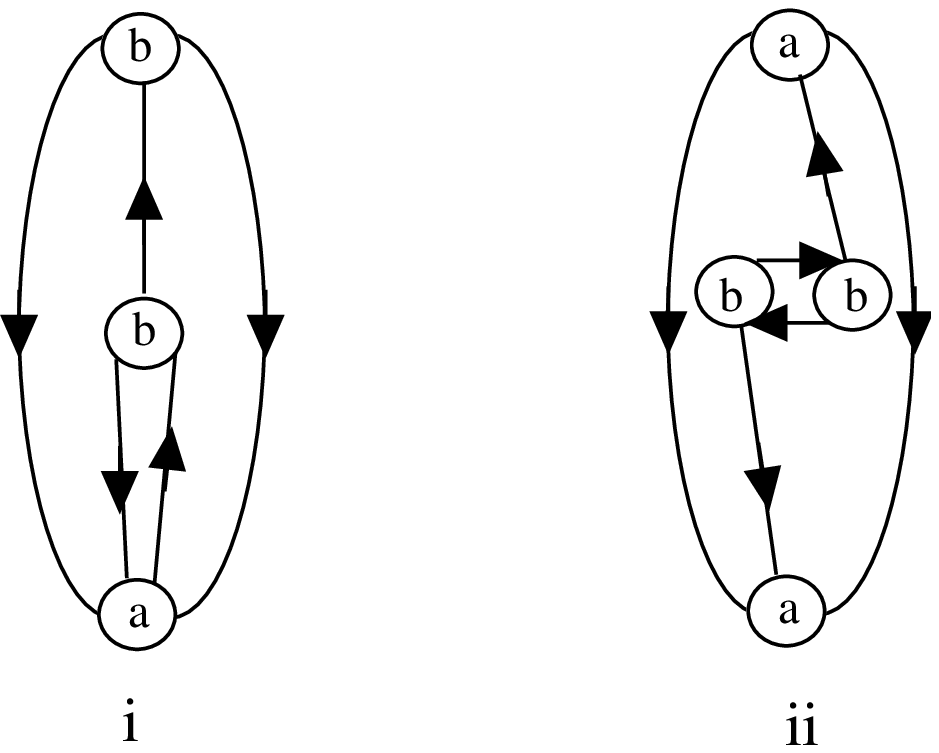}
\caption{} \label{fig:onea}
\end{figure}

So there are exactly two $b$-vertices inside the bigon.  As noted 
above, at least one edge pointing out from each of these $b$-vertices 
must go to an $a$-vertex, and edges can't point from different 
$b$-vertices to the same $a$-vertex.  It follows that both $v_{\pm}$ 
are $a$-vertices.  If both of the edges pointing out from a $b$-vertex 
go to $a$-vertices then, since $p = 1$, they would in fact go to the 
same $a$-vertex, contradicting our assumption that the parallel bigon 
is innermost.  So each $b$-vertex has an edge pointing from it to the 
other $b$-vertex.  In other words, the two $b$-vertices in the bigon 
are the vertices of a $2$-cycle, necessarily a $b$-face cycle.  (See 
Figure \ref{fig:onea}ii.)  But this leads to the same contradiction as 
in the proof of Case 1 of Lemma \ref{lemma:nobigona}.

\medskip

If $p \geq 2$ then we know from Corollary \ref{cor:endb} that $w$ 
either begins or ends in $\ob$, so the words corresponding to 
$\bdd_{\pm} A$ are (up to cyclic rotation) $a^m \ob^{p+1}$ and 
$a^{m+1} \ob$ respectively.  Here $l(a) = pl(b)$ as usual.  Suppose 
first that no edge in $G$ runs from one $b$-vertex to another 
$b$-vertex.  Then the $p+1$ occurences of any $b$ vertex in $\bdd_+ A$ 
have their other ends in only the $p$ possible $a$-vertices, so there 
is a parallel bigon at each $b$-vertex.  An innermost parallel bigon 
then contains only $a$-vertices.  This contradicts Lemma 
\ref{lemma:nobigona}.

Suppose then that some edge $\bbb$ in $G$ has both ends at $b$-vertices.  
Consider the distance in $A$ between copies of the same $a$-label, 
counting distance (i. e. intersection with $P$) along the arcs of 
$\bdd A$ that {\em don't} intersect $\bbb$.  Measured on this side, 
the distance between any two copies of the same $a$-label in either 
component of $\bdd A$ is a multiple of $l(a)$.  It follows that the 
$m+1 \geq 3$ copies of the same label in $\bdd_{-} A$ have their other 
ends at at most two labels in $\bdd_{+}A$, one an $a$-label and one a 
$b$-label.  Thus in this case every $a$-vertex is part of a parallel 
bigon.  Then an innermost bigon contains only $b$-vertices.  The proof 
now follows as for the case $p = 1$ but is easier, since any 
$b$-vertex has $p+1 \geq 3$ edges pointing out.
\end{proof}

At this point there are only two remaining words to consider: $w = 
a\ob$ and $w = \ob a$.  Eliminating these two requires a bit more 
detailed argument.

\begin{lemma} \label{lemma:noab}  $w \neq a\ob$ or $\ob a$.
\end{lemma}

\begin{proof} The cases are symmetric, so without loss of generality 
suppose $w = a\ob$.  Then $\bdd_- A$ is represented by the word $a\ob 
a$ and $\bdd_+ A$ is represented by the word $a\ob^{p+1}$.  Let $n = 
l(a) = pl(b), p \geq 1$.  The type of contradiction depends on how the 
copies of $\ob$ in $\bdd A$ are aligned with each other.  There are 
three cases.  See Figure \ref{fig:noab}, where the orientation of $w$ is meant 
to be clockwise around $A$. 

\begin{figure}
\centering
\includegraphics[width=0.65\textwidth]{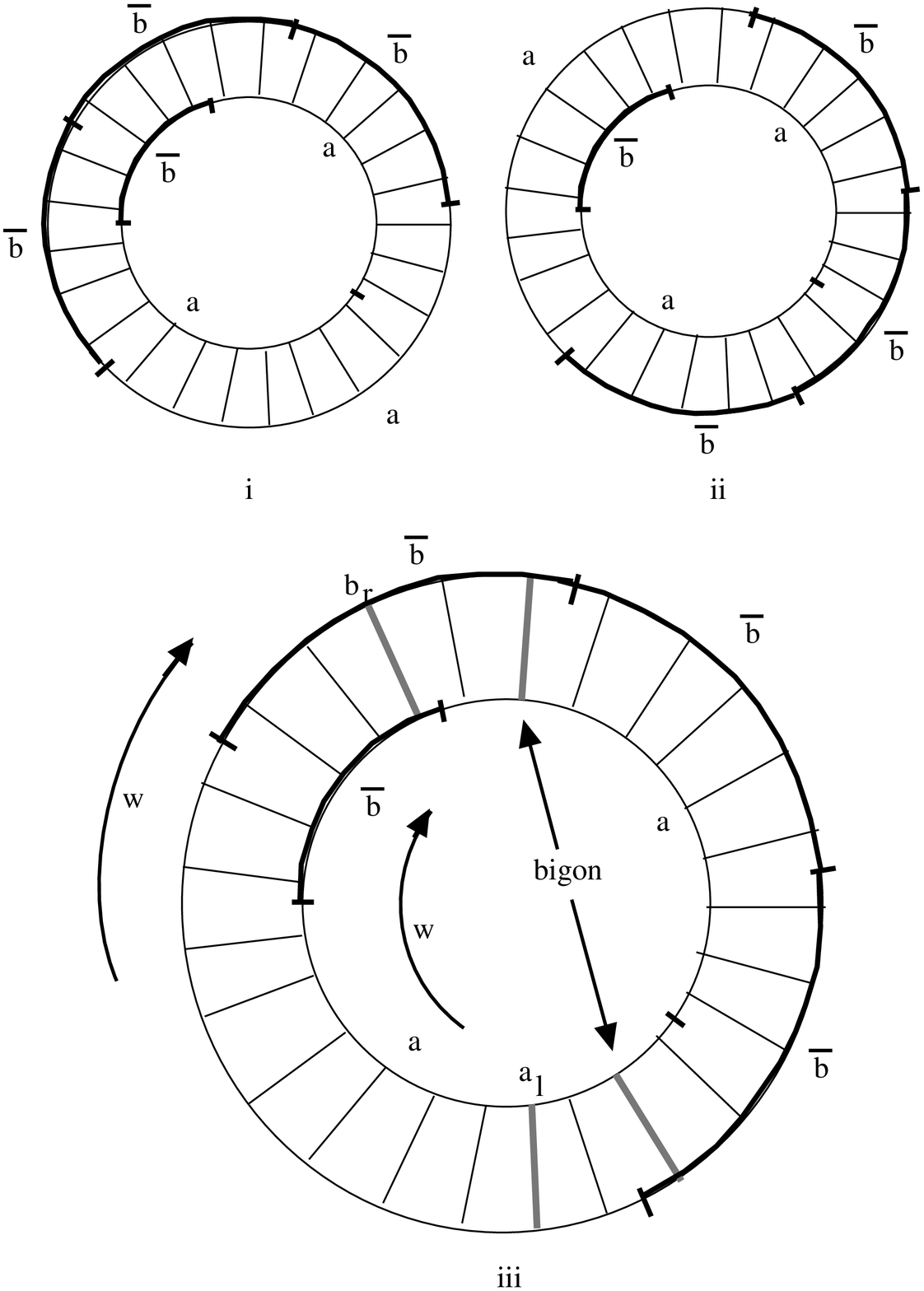}
\caption{}  \label{fig:noab}
\end{figure}

\medskip

{\bf Case 1} The $\ob$-segment in $\bdd_{-} A$ lies completely 
opposite a subsegment of the $\ob^{p+1}$ segment in $\bdd_{+}A$.  (Figure 
\ref{lemma:nobigona}i)

Then symmetrically, the $a$ segment in $\bdd_{+} A$ lies completely 
opposite a subsegment of the $a^{2}$ segment in $\bdd_{-}A$.  Consider 
the collection of edges incident to the copy of the $a$-segment in 
$\bdd_{+} A$.  Among those edges, every $a$-label occurs exactly once 
in $\bdd_{+}A$ and once in $\bdd_- A$.  It follows that these edges, 
when viewed in $G \subset P$, form a collection of $a$-cycles 
containing every $a$-vertex.  Similarly there is a collection of 
$b$-cycles containing every $b$-vertex.  An innermost pure $a$- or pure 
$b$-cycle can then contain no vertices in its interior and so must be 
an $a$-face or a $b$-face.  Either way, the argument presented in Case 
1 of Lemma \ref{lemma:nobigona} presents a contradiction.

\medskip

{\bf Case 2:} The $\ob$-segment in $\bdd_{-} A$ lies completely 
opposite a subsegment of the $a$ segment in $\bdd_{+}A$.  (Figure 
\ref{lemma:nobigona}ii)

Then dually the entire $\ob^{p+1}$-segment in $\bdd_{+} A$ lies 
completely opposite a subsegment of the $a^{2}$ segment in 
$\bdd_{-}A$.  This means that for each $b$-vertex, each of the $p+1$ 
edges in $A$ with that label in $\bdd_{+} A$ can go to at most $p$ 
different $a$-labels in $\bdd_{-} A$.  It follows that every 
$b$-vertex is part of a parallel bigon.  An innermost one can contain 
only $a$-vertices, contradicting Lemma \ref{lemma:nobigona}.

{\bf Case 3:} The $\ob$-segment in $\bdd_{-} A$ lies partly 
opposite an end of the $a$ segment in $\bdd_{+}A$ and partly 
opposite an end of the $\ob^{p}$ segment. (Figure 
\ref{lemma:nobigona}iii)

The argument in this case is a kind of degenerate variant of the 
argument in Lemma \ref{lemma:nobigon}.  Suppose, with no loss and as 
shown in the figure, that part of the $\ob$ segment in 
$\bdd_{-} A$ is opposite the beginning end of $\ob^{p}$, overlapping 
say on $j < l(b)$ edges.  ($j = 3, l(b) = 5$ in the figure).  Now 
consider any of the last $n-j$ labels in the first occurence of $\ob$ 
in $\bdd_{+} A$ and the corresponding label in the last occurence of 
$\ob$.  The distance between them is $l(a)$ in $\bdd_+ A$.  It follows 
that each of the corresponding $(n - j)$ $b$-vertices is part of a 
parallel bigon in $P$, each with an $a$-vertex for $v_{-}$.

Consider an innermost parallel bigon in $P$ and the disk $B$ that it 
bounds.  The interior of $B$ must contain $b$-vertices, by Lemma 
\ref{lemma:nobigona}, and it must also contain $a$-vertices, by the 
argument of Case 1 of the proof of that Lemma.  We claim that there is 
an oriented edge pointing from $v_{-}$ into $B$ and an oriented edge 
pointing out from $B$ into $v_{+}$.  To see the former, consider the 
$a$-vertex $a_{l}$ which, among all $a$-vertices lying in $B$, is the 
first encountered by $e^{-}$.  Then the corresponding label in the 
second copy of $a$ in $a^2 \subset \bdd_-A$ lies across from the 
label of an earlier vertex in $e^{-}$, hence the label is that of $v_{-}$, 
since there is no alternative.  So the edge in $P$ between them must 
connect $a_{l}$ to $v_{-}$, pointing toward $a_{l}$. A symmetric 
argument, using the last $b$-vertex $b_r$ in the interior of $B$ 
encountered by $e^{+}$, shows that there is an edge pointing from 
$b_{r}$ into $v_+$.  

Since each $a$-vertex has valence $3$ we have now accounted for all 
edges incident to $v_{-}$.  In particular one of the two edges of the 
bigon has a $w$-end at $v_{-}$ with the $w$ side lying within the 
bigon.  What's important here is not that one of the edges of the 
bigon has a $w$-end at $v_{-}$ -- that fact can be seen simply because 
one of the ends of the edges of the bigon, viewed in $A$, lies in the 
second occurence of $a$ in $a^2$ (see Figure \ref{fig:noab}iii), hence in 
the $a\ob$ section of the word $wa = a\ob a$.  What is important is 
that the $w$ side of this edge lies in the interior of the bigon.  But 
examining the figure again, we see that the {\em other} edge of the 
bigon lies in the first occurence of $\ob$ in $\ob^{p+1}$, hence in 
the $w$ section of the word $w\ob^{p} = a\ob \ob^{p}$.  That is, the 
{\em other} edge has a $w$ end at $v_+$.  So the $w$-side of that edge 
lies {\em outside} of $B$.  On the other hand, we've shown that some 
edge at $v_+$ points into $v_+$ from the interior of $B$ and, since 
any $b$-vertex has only one edge pointing into it, that edge must be 
adjacent to the (only) $w$-corner at $v_+$, so that corner must be 
{\em inside} B. This contradiction proves the Lemma, hence the 
theorem.
\end{proof}

Having eliminated every possible word for $w$, we deduce that no 
dividing sphere can intersect $A$ only in essential arcs, completing 
the proof of Theorem \ref{theorem:ebot}
\end{proof}

\section{The Goda-Teragaito Conjecture}

\begin{theorem}  \label{theorem:GodaT} 
Suppose $K$ is a tunnel number one knot of genus one 
and $\ggg$ is an unknotting tunnel.  Then either there is a genus one 
Seifert surface $F$ for $K$ that contains $\ggg$ or $\ggg$ can be slid 
and isotoped until it is an unknotted loop.
\end{theorem}

\begin{proof}  According to Proposition \ref{prop:oneonethin} if 
neither of the outcomes above occurs, then there is an appropriate 
$\Theta$-graph for $(K, F)$, thinly presenting it, say, as a $(p, q)$ 
quasi-cable.  Let $\theta$ be a thinnest appropriate $\Theta$-graph 
for $(K, F)$ and, among all such possibilities, choose one with $p$ 
maximal.  According to Propositions \ref{prop:oneonethin} and
\ref{prop:oneonethin2} $q$ = 1 and $\theta$ is in bridge position.  

We claim that if $p \geq 2$ the cycle $e^{\bot} \cup e^{+}$ is 
unknotted and, if $p = 1$, one of the two cycles $e^{\bot} \cup 
e^{\pm}$ is unknotted.  This follows immediately from Corollary 
\ref{cor:edgeabove} unless a dividing sphere is a critical sphere that 
is disjoint from some edge.  Consider the possibilities for such an 
edge: If the disjoint edge is $e^{\bot}$ then the claim is established 
by Theorem \ref{theorem:ebot}.  If $e^{-}$ is the disjoint edge then, 
since $\theta$ has been chosen to have maximal $p$, it follows from 
Proposition \ref{prop:eminus} that $p = 1$ and $e^{\bot} \cup e^{-}$ 
is unknotted, establishing the claim.  Similarly, if $e^{+}$ is the 
disjoint edge then it follows from 
Proposition \ref{prop:eplus2} that the wave is based at $\mm$.  Then 
Theorem \ref{theor:eplus} establishes the claim.  So the claim is 
established in all cases.

Now let $L$ be the unknotted solid torus neighborhood of $e^+ \cup 
e^{\bot}$ in $H$.  Since $q = 1$ we can apply the ``vacuum cleaner 
trick'': slide the ends of the $1$-handle corresponding to $e^-$ along 
the arc $K \cap L$ until $K$ has been made disjoint from a meridian of 
$L$.  At that point, $L$ has become a tunnel for $K$ and remains 
unknotted.  \end{proof}

\begin{cor} (Goda-Teragaito Conjecture) \label{cor:GodaT} 
Suppose $K$ is a tunnel number one knot of genus one that is not a 
satellite knot.  Then $K$ is 
$2$-bridge.

\end{cor}

\begin{proof} Let $\ggg$ be an unknotting tunnel for $K$.  If $\ggg$ 
can be slid and isotoped to lie on a genus one Seifert surface $F$ 
then $K$ is necessarily $2$-bridge (see Corollary 5.4 of \cite{ST1}).  
If not, then according to Theorem \ref{theorem:GodaT} $\ggg$ can be 
slid and isotoped until it is an unknotted loop.  The following 
argument (shown to me by Abby Thompson) shows that then $K$ is 
$1$-bridge on an unknotted torus.  Let $W$ denote the 
solid torus neighborhood of the loop, containing a short, 
$\bdd$-parallel arc of $K$.  Let $K_{-}$ denote the arc of $K$ that 
lies outside of $W$.  Since $\eta(K \cup \ggg)$ is an unknotted 
handlebody, it follows that the $1$-handle with $K_{-}$ at its core 
constitutes a genus two Heegaard splitting of the solid torus $S^3 - 
W$.  Any non-trivial splitting of a handlebody (e.  g.  of $S^{3} - 
W$) is stabilized \cite{ST3}, so in fact $K_-$ is also parallel to $\bdd W$.  
This shows that $K$ is $1$-bridge with respect to the unknotted 
torus $\bdd W$. 

Matsuda  \cite{Ma} has proven the statement for this 
class of knots.
\end{proof}

\end{document}